\documentclass[a4paper,german,final,twoside,notitlepage,12pt]{article}

\newcommand{\secondproof}{Second proof}

\usepackage{amssymb}
\usepackage{amsthm}
\usepackage{amsmath}
\usepackage[all]{xy}
\usepackage{graphicx}
\usepackage{mathrsfs}
\usepackage{afterpage}
\usepackage{float}
\usepackage[bottom]{footmisc}
\usepackage[margin=2cm,font=normalsize,labelfont=bf]{caption}
\usepackage[normalem]{ulem}
\usepackage{xstring}
\usepackage[section]{placeins}
\usepackage{verbatim}
\usepackage{amstext}

\makeatletter

\gdef\@ntitle{\@title}
\def\adress#1{\gdef\@adress{#1}}
\def\@adress{}
\def\preprint#1{\gdef\@preprint{#1}}
\def\@preprint{}
\def\keywords#1{\gdef\@keywords{#1}}
\def\@keywords{}
\def\email#1{\gdef\@email{#1}}
\def\@email{}

\def\refname{References}

\newlength{\myparskip}
\def\href#1#2{#2}

\def\kohyp{
  \usepackage{hyperref}
  \hypersetup{
    linktocpage = true,
    pdftitle = {\@title},
    pdfauthor = {\@author},
    pdfkeywords = {\@keywords},    
    bookmarksopen = true,
    bookmarksopenlevel = 1
  }}  
\def\showkeywords{\begin{flushleft}\footnotesize\textbf{Keywords}: \@keywords.\end{flushleft}}

\newcounter{mythm}[section]

\numberwithin{equation}{section}

\newtheorem{theorem}[mythm]{Theorem}
\newtheorem{definition}[mythm]{Definition}
\newtheorem{rem}[mythm]{Remark}
\newtheorem{corollary}[mythm]{Corollary}
\newtheorem{exa}[mythm]{Example}
\newtheorem{proposition}[mythm]{Proposition}
\newtheorem{lemma}[mythm]{Lemma}
\newtheorem{con}[mythm]{Construction}

\newenvironment{remark}{\begin{rem}\normalfont\setlength{\parskip}{0cm}}{\setlength{\parskip}{\myparskip}\end{rem}}

\def\N {\mathbb{N}}

\def\R {\mathbb{R}}
\def\C {\mathbb{C}}

\def\id{\mathrm{id}}

\def\lli#1{\,\,_{#1}\!}

\renewcommand{\varepsilon}{\epsilon}

\newcommand{\alxy}[1]{\begin{aligned}\xymatrix{#1}\end{aligned}}
\newcommand{\alxydim}[2]{\begin{aligned}\xymatrix#1{#2}\end{aligned}}

\renewenvironment{proof}{Proof.\setlength{\parskip}{0cm}
}{\hfill{$\square$}\\\setlength{\parskip}{\myparskip}
}
\newenvironment{proofblank}[1]{#1.\ }{\hfill{$\square$}\\}

\newcommand\erf[1]{(\ref{#1})}

\renewcommand{\emph}[1]{\def\reserved@a{it}\ifx\f@shape\reserved@a\uline{#1}\else\textit{#1}\fi}
\newlength{\myl}
\newcommand\sheaf[1]{\unitlength 0.1mm
  \settowidth{\myl}{$#1$}
  \addtolength{\myl}{-0.8mm}
  \begin{picture}(0,0)(0,0)
  \put(3,6){\text{\uline{\hspace{\myl}}}}
  \end{picture}#1\hspace{-0.15mm}}

\def\trivlin{\mathbf{I}}

\makeatother

\renewcommand{\to}{\!\xymatrix@C=0.5cm{\ar[r] &}}
\renewcommand{\mapsto}{\!\xymatrix@C=0.5cm{\ar@{|->}[r] &}\!}
\renewcommand{\Rightarrow}{\!\xymatrix@C=0.5cm{\ar@{=>}[r] &}\!}
\newcommand{\incl}{\!\xymatrix@C=0.5cm{\ar@{^(->}[r] &}\!}
\renewcommand{\hookrightarrow}{\incl}
\def\Leftrightarrow{\!\xymatrix@C=0.5cm{\ar@{<=>}[r] &}\!}

\makeatletter

\def\quot#1{``#1''}

\def\mytitle{}
\def\zmptitle{
  \begin{tabular}{cc}
    \begin{minipage}[c]{0.4\textwidth}
      \begin{flushleft}
        \includegraphics[width=110pt]{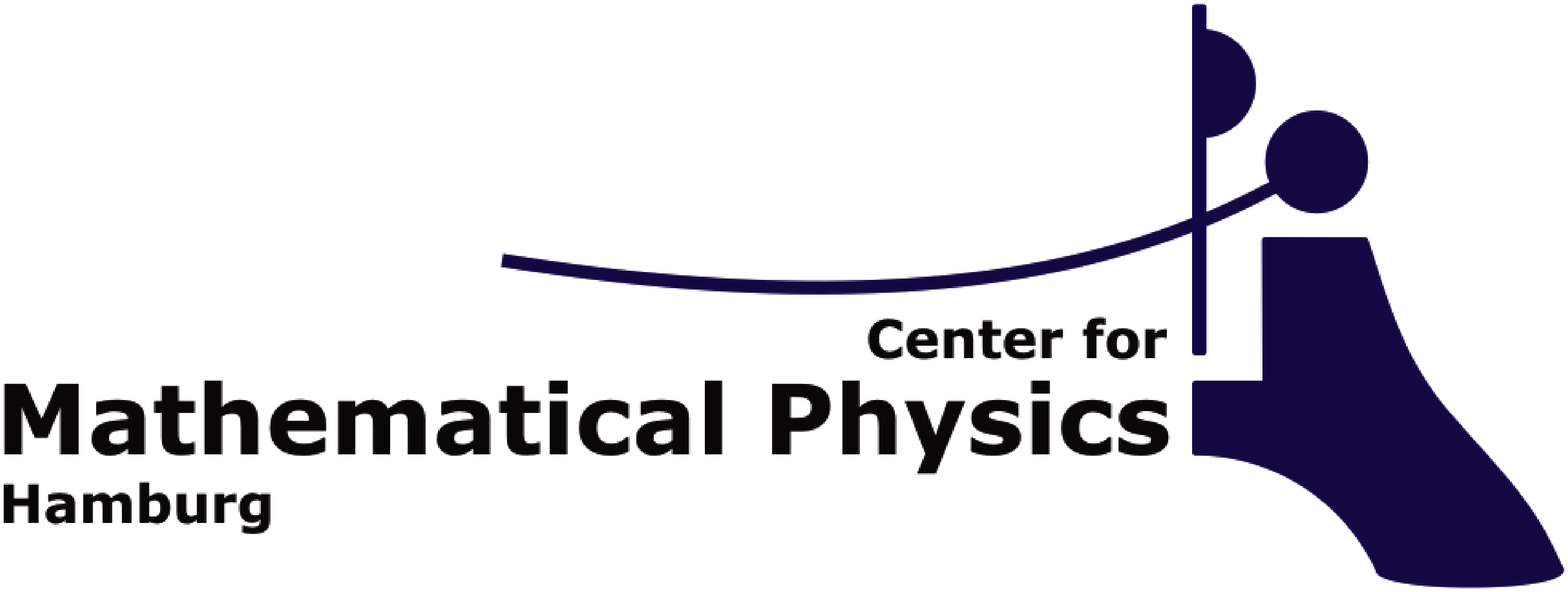}
      \end{flushleft}  
    \end{minipage}&
    \begin{minipage}[c]{0.55\textwidth}
      \begin{flushright}
      {\small\sf\@preprint}
      \end{flushright}
    \end{minipage}
  \end{tabular}
  \vskip 2cm}

\def\maketitle{
  \setlength{\parskip}{\myparskip}
  \newpage
  \noindent
  \mytitle
  \begin{center}
    \LARGE\@ntitle
    \if!\@author!\else \vskip 0.5cm \large\@author\fi
    \if!\@adress!\else \vskip 0.2cm \normalsize\@adress\fi
    \if!\@email!\else \vskip 0.2cm \normalsize\textit{\@email}\fi
  \end{center}
  \vskip 2cm\thispagestyle{empty}}

\def\kobiburl#1{
   \IfBeginWith
     {#1}
     {http://arxiv.org/abs/}
     {\kobibarxiv{#1}}
     {\kobiblink{#1}}}
\def\kobibarxiv#1{\href{#1}{\texttt{[arxiv:\StrGobbleLeft{#1}{21}]}}}
\def\kobiblink#1{
  \StrSubstitute{#1}{_}{\underline{\;\;}}[\mylink]
  \StrSubstitute{\mylink}{&}{\&}[\mylink]
  Available as: \mbox{\;}
  \href{#1}{\texttt{\mylink}}}

\def\kobib#1{

}

\def\showcomments{ -- Comments suppressed}

\newif\if@fewtab\@fewtabtrue{
  \count255=\time\divide\count255 by 60
  \xdef\hourmin{\number\count255}
  \multiply\count255 by-60\advance\count255 by\time
  \xdef\hourmin{\hourmin:\ifnum\count255<10 0\fi\the\count255}}
\def\ps@draft{
  \let\@mkboth\@gobbletwo
  \def\@oddfoot{
    \hbox to 7 cm{\tiny \versionno\hfil}
    \hskip -7cm\hfil\rm\thepage\hfil{\tiny\draftdate}}
  \def\@oddhead{}
  \def\@evenhead{}
  \let\@evenfoot\@oddfoot}
\def\draftdate{\number\month/\number\day/\number\year\ \ \ \hourmin }
\newcommand\version[1]{
  \typeout{}\typeout{#1}\typeout{}
  \vskip-1.7cm \centerline{\fbox{{\normalsize\tt DRAFT -- #1 -- 
  \draftdate\showcomments}}} \vskip0.92cm}
\def\draft#1{
  \def\versionno{#1}
  \pagestyle{draft}\thispagestyle{draft}
  \gdef\@ntitle{\version\versionno \@title}
  \global\def\draftcontrol{1}}
\global\def\draftcontrol{0}
\def\todo#1{\ifnum\draftcontrol>0 \\\bigskip {\tt #1} \\ \fi}

\makeatother

\usepackage[latin1]{inputenc}

\newcommand{\loctrivfunct}[3]{\mathrm{Triv}^{#3}_{\pi}(#1)}
\newcommand{\loctrivfunctsmooth}[4]{\mathrm{Triv}^{#4}_{#3}(#1)^{\infty}}
\newcommand{\trans}[3]{\mathfrak{Des}_{\pi}^{#3}(#1)}
\newcommand{\ex}[1]{\mathrm{Ex}_{#1}}
\newcommand{\transsmooth}[3]{\mathfrak{Des}^{#2}_{\pi}(#1)^{\infty}}
\newcommand{\transport}[4]{\mathrm{Trans}^{#2}_{#3}(M,#4)}
\newcommand{\diffco}[3]{Z^{#2}_{#3}(#1)^{\infty}}
\newcommand{\gconn}[1]{\diffco{G}{1}{#1}}
\newcommand{\gbun}{\mathfrak{Bun}^{\nabla}_G(M)}
\newcommand{\upp}[2]{\mathcal{P}^{#2}_1(#1)}
\newcommand{\uppp}{p^{\pi}}
\newcommand{\smsp}{D^{\infty}}
\newcommand{\set}{\mathfrak{Set}}
\newcommand{\sm}{C^{\infty}}
\newcommand{\fu}{\mathcal{P}}
\newcommand{\fo}{\mathcal{D}}
\newcommand{\funct}{\mathrm{Funct}}

\title{Parallel Transport and Functors}
\author{Urs Schreiber and Konrad Waldorf}
\adress{Organisationseinheit Mathematik\\Schwerpunkt Algebra und Zahlentheorie\\Universit\"at
Hamburg\\Bundesstra\ss e 55\\D--20146 Hamburg
}
\preprint{arxiv:0705.0452\\Hamburger Beitr\"age zur Mathematik Nr. 269\\ZMP-HH/07-5}

\kohyp
\usepackage{amstext}

\begin{document}


\def\mytitle{\zmptitle}
\maketitle 

\begin{abstract}
Parallel transport of a connection in a smooth fibre  bundle yields  a  functor from the path groupoid of the base
manifold 
into a  category that describes the  fibres of the bundle. We characterize functors obtained like this by   two notions we introduce: local trivializations and smooth descent data. This provides a way to substitute
  categories of
functors for categories of smooth fibre bundles with connection. We indicate that this concept can be generalized to connections in categorified bundles, and how this generalization improves the understanding of higher dimensional parallel transport. 
\end{abstract}

\newpage

{
\footnotesize
\tableofcontents
\thispagestyle{empty}
}

\setcounter{page}{1}

\section{Introduction}

\label{sec8}

Higher dimensional parallel transport generalizes parallel transport along
curves to parallel transport along higher dimensional objects, for instance  surfaces. One motivation
to consider  parallel transport along surfaces comes from two-dimensional
conformal field theories, where so-called Wess-Zumino terms have been recognized as surface holonomies  \cite{gawedzki3,carey2,schreiber1}.

Several mathematical objects
have have been used to define higher dimensional parallel transport, among
them classes in Deligne cohomology \cite{deligne1}, bundle gerbes with connection and curving \cite{murray}, or 2-bundles with 2-connections \cite{baez,baez2}. The development of such definitions often occurs in two steps:
 an appropriate definition of parallel transport along curves,
followed by a generalization  to higher dimensions. For instance,
bundle gerbes with connection can be obtained as a generalization of principal bundles
with connection. However,  in the case of both bundle gerbes and Deligne classes one encounters the obstruction that the structure group has to be abelian.  It is hence desirable to find  a reformulation of fibre bundles with connection, that brings along a natural generalization for arbitrary structure group.

A  candidate for such a  reformulation are holonomy maps
\cite{barret1,caetano}. These are   group homomorphisms
\begin{equation*}
\mathcal{H}:\pi_1^1(M,*) \to G
\end{equation*}
from the group of thin homotopy classes of based loops in a smooth manifold  $M$ into a   Lie group $G$. Any principal $G$-bundle with connection over $M$ defines
a group homomorphism $\mathcal{H}$, but the crucial point is to  distinguish 
those   from arbitrary ones. By imposing a certain smoothness
condition on $\mathcal{H}$, these holonomy maps correspond -- for connected manifolds -- bijectively
to principal $G$-bundles with connection \cite{barret1,caetano}.  On the other hand, they have a natural generalization from loops to surfaces. However,
the obstruction for $M$ being connected  becomes  even stronger: only if the manifold $M$ is connected and simply-connected, holonomy maps generalized to surfaces capture all aspects of surface holonomy \cite{mackaay1}. Especially the  second obstruction  erases one of the most interesting of these
aspects, see, for example, \cite{gawedzki1}. 

\medskip

In order to obtain  a formulation
of parallel transport  along curves
without topological assumptions on the base manifold $M$, one considers functors
\begin{equation*}
F: \mathcal{P}_1(M) \to T
\end{equation*}
from the path groupoid $\mathcal{P}_1(M)$
of  $M$ into another category $T$ \cite{mackenzie,mackaay1}. 
The set of objects of the path groupoid $\mathcal{P}_1(M)$ is the manifold $M$ itself, and the set of morphisms between two points
$x$ and $y$ is the set of thin homotopy classes
of curves starting at $x$ and ending at $y$. A functor $F: \mathcal{P}_1(M) \to T$ is a generalization
of a group homomorphism $\mathcal{H}:\pi_1^1(M,*) \to G$,
but it is not clear how the  smoothness condition for holonomy maps has to
be generalized to these functors. 

Let us first review how  a functor $F: \mathcal{P}_1(M) \to T$ arises from parallel transport in  a, say, principal $G$-bundle
$P$ with connection. In this case, the category $T$ is the category $G\text{-}\mathrm{Tor}$ of
smooth manifolds with smooth, free and transitive  $G$-action from the right, and smooth equivariant maps
between those.
Now, the
connection on $P$ associates to any smooth curve $\gamma: [0,1] \to M$ and
any element in the fibre $P_{\gamma(0)}$ over the starting
point,  a unique horizontal lift $\tilde \gamma:[0,1] \to P$.
Evaluating this lift at its endpoint defines a smooth map
\begin{equation*}
\tau_{\gamma}: P_{\gamma(0)} \to P_{\gamma(1)}\text{,}
\end{equation*}
the parallel transport in $P$ along the curve $\gamma$. It is $G$-equivariant
with respect to the
$G$-action on the fibres of $P$, and
it is  invariant under  thin homotopies. Moreover, it satisfies
\begin{equation*}
\tau_{\id_x}=\id_{P_x}
\quad\text{ and }\quad
\tau_{\gamma'
\circ \gamma} = \tau_{\gamma'} \circ \tau_{\gamma}\text{,}
\end{equation*} 
where $\id_x$ is the constant curve and $\gamma$ and $\gamma'$ are smoothly
composable curves. These are the axioms of a functor
\begin{equation*}
\label{2}
\mathrm{tra}_P: \mathcal{P}_1(M) \to G\text{-}\mathrm{Tor}
\end{equation*}
which sends an object $x$ of $\mathcal{P}_1(M)$ to the object $P_x$ of
$G\text{-}\mathrm{Tor}$ and a morphism $\gamma$ of $\mathcal{P}_1(M)$ to
the morphism $\tau_{\gamma}$ of $G\text{-}\mathrm{Tor}$. Summarizing, every principal $G$-bundle with connection over $M$ defines a functor
$\mathrm{tra}_P$.  Now the crucial point is to characterize these functors among
all functors from $\mathcal{P}_1(M)$ to $G\text{-}\mathrm{Tor}$. 

\medskip

In this article we describe such a characterization. For this purpose, we introduce, for general target categories $T$, the notion  of a
transport functor. These are certain functors
\begin{equation*}
\mathrm{tra}: \mathcal{P}_1(M) \to T\text{,}
\end{equation*}
such that the category they form is -- in the case of $T=G\text{-}\mathrm{Tor}$ -- equivalent to the category of principal $G$-bundles with connection. 

The defining properties of a transport functor capture two
important concepts: the existence of local trivializations and the smoothness  of  associated descent data.
Just as for fibre bundles, local
trivializations are specified with respect to an open cover of the base manifold $M$ and   to a choice of a typical
fibre. Here, we  represent an open cover by a surjective
submersion $\pi:Y \to M$, and encode the typical fibre in the notion of a structure groupoid: this is a Lie groupoid $\mathrm{Gr}$ together with a functor 
\begin{equation*}
i: \mathrm{Gr} \to T\text{.}
\end{equation*}
Now, a $\pi$-local $i$-trivialization of a functor $F:\mathcal{P}_1(M) \to T$
is another functor 
\begin{equation*}
\mathrm{triv}: \mathcal{P}_1(Y) \to \mathrm{Gr}
\end{equation*}
together with
a natural equivalence
\begin{equation*}
t: F \circ \pi_{*} \to i \circ \mathrm{triv}\text{,}
\end{equation*}
where $\pi_{*}: \mathcal{P}_1(Y) \to \mathcal{P}_1(M)$ is the induced functor between
 path groupoids. In detail, the natural equivalence $t$ gives for every point $y\in Y$ an isomorphism $F(\pi(y))
\cong i(\mathrm{triv}(y))$ that identifies the \quot{fibre} $F(\pi(y))$ of $F$ over $\pi(y)$ with the image
 of a \quot{typical fibre}
$\mathrm{triv}(y)$ under the functor $i$.
In other words, a functor is $\pi$-locally $i$-trivializable, if its pullback to the cover $Y$ factors through the functor
$i$ up to a natural equivalence. Functors with a chosen $\pi$-local $i$-trivialization $(\mathrm{triv},t)$
form a category $\loctrivfunct{i}{\mathrm{Gr}}{1}$.

The second concept we introduce is that of smooth descent data. Descent data is specified with respect to a surjective submersion $\pi$ and a structure groupoid
$i:\mathrm{Gr} \to T$.
While  descent data for a fibre bundle with connection is a collection
of transition
 functions and local connection 1-forms,  descent  data for a functor 
 $F:\mathcal{P}_1(M) \to T$ is a pair $(\mathrm{triv},g)$
consisting of a functor $\mathrm{triv}:\mathcal{P}_1(Y) \to \mathrm{Gr}$
 like the one from a local trivializations and of a certain natural equivalence $g$ that compares $\mathrm{triv}$ on the two-fold fibre product of $Y$ with itself. Such pairs  define  a descent category
$\trans{i}{\mathrm{Gr}}{1}$. The first result of this article (Theorem \ref{th3}) is to prove the
descent property: extracting descent data and, conversely, reconstructing
a functor from descent data, are equivalences of categories
\begin{equation*}
\loctrivfunct{i}{\mathrm{Gr}}{1} \cong \trans{i}{\mathrm{Gr}}{1}\text{.}
\end{equation*}
We introduce  descent data because one can  precisely decide whether
a pair $(\mathrm{triv},g)$ is smooth or
not (Definition \ref{def1}).  The smoothness conditions we introduce can be expressed
in basic terms of smooth maps between smooth manifolds, and arises from the theory of  diffeological spaces \cite{chen1}. The concept of smooth descent data is our generalization of the smoothness
condition for holonomy maps to functors.

Combining both concepts we have introduced, we  call  a functor that allows -- for some surjective submersion $\pi$ -- a $\pi$-local $i$-trivialization
whose corresponding descend data
is smooth, a transport functor on $M$
in $T$ with $\mathrm{Gr}$-structure. The category formed by these transport functors is denoted by $\transport{i}{1}{\mathrm{Gr}}{T}$. 

\medskip

Let us return  to the particular target category $T=G\text{-}\mathrm{Tor}$.
As described above, one obtains a functor $\mathrm{tra}_P:\mathcal{P}_1(M) \to G\text{-}\mathrm{Tor}$
from any principal $G$-bundle $P$ with connection. We consider
the  Lie groupoid $\mathrm{Gr}=\mathcal{B}
G$, which has only one object, and where every group element $g \in G$ is
an automorphism of this object. The notation indicates the fact that the geometric realization of the nerve of this category yields the classifying space $BG$ of the group $G$. The Lie groupoid $\mathcal{B}G$ can be embedded in the
category $G\text{-}\mathrm{Tor}$  via the functor
$i_G: \mathcal{B}G \to G\text{-}\mathrm{Tor}$ which sends the object of $\mathcal{B}
G$ to
the group $G$ regarded as a $G$-space, and a morphism $g\in G$ to the equivariant
smooth map  which   multiplies with $g$ from the left.

The descent category $\trans{i_G}{}{1}$ for the structure groupoid $\mathcal{B}G$ and some surjective submersion $\pi$ is closely related to differential
geometric objects: we  derive a one-to-one correspondence between
smooth functors $\mathrm{triv}:\mathcal{P}_1(Y) \to \mathcal{B}G$, which are part of the objects of  $\trans{i_G}{\mathcal{B}G}{1}$, and 1-forms $A$ on $Y$ with values in the Lie algebra
of $G$ (Proposition \ref{th2}). The correspondence can be symbolically expressed as the
path-ordered exponential
\begin{equation*}
\mathrm{triv}(\gamma) = \mathcal{P}\exp \left ( \int_{\gamma} A \right)
\end{equation*}
for a path $\gamma$. Using this relation between smooth
functors and differential forms, we show that a functor
$\mathrm{tra}_P:\mathcal{P}_1(M) \to G\text{-}\mathrm{Tor}$ obtained
from a principal $G$-bundle with
connection, is a transport
functor on $M$ in $G\text{-}\mathrm{Tor}$ with $\mathcal{B}G$-structure. 
The main   result of this article (Theorem \ref{th1}) is that this  establishes
an equivalence of categories
\begin{equation*}
\gbun \cong \transport{i}{1}{\mathcal{B}G}{G\text{-}\mathrm{Tor}}
\end{equation*} 
between the category of  principal $G$-bundles with connection over $M$ and
the category of transport functors on $M$ in $G\text{-}\mathrm{Tor}$ with $\mathcal{B}
G$-structure. In other words, these transport functors  provide a proper reformulation of principal bundles
with connection, emphasizing the aspect of parallel transport.

\medskip

This article is organized as follows. In Section
\ref{sec2} we review the path groupoid of a smooth manifold and describe some properties of functors
defined on it. We introduce local trivializations for functors and the descent category
$\trans{i}{\mathrm{Gr}}{1}$.
In Section \ref{sec3} we  define the category $\transport{i}{1}{\mathrm{Gr}}{T}$ of transport functors on $M$ in $T$ with $\mathrm{Gr}$-structure and discuss several properties. In Section \ref{sec7} we derive the result that relates the descent category
$\trans{i_G}{\mathrm{Gr}}{1}$ for the particular functor $i_G:\mathcal{B}G \to G\text{-}\mathrm{Tor}$
to differential forms. In  Section \ref{sec4} we provide examples that show that the
theory of transport functors  applies well to several situations: we prove our main result concerning principal $G$-bundles with connection, show a similar statement for 
vector bundles with connection, and also discuss holonomy maps. In Section \ref{sec15} we discuss principal groupoid bundles and show how transport functors can be used to derive the definition of a connection on such groupoid bundles. Section \ref{sec1} contains various directions in which the concept of transport functors can be generalized. In particular, we outline a possible generalization of transport functors
to transport $n$-functors
\begin{equation*}
\mathrm{tra}:\mathcal{P}_n(M) \to T\text{,}
\end{equation*}
which provide an implementation for higher dimensional parallel transport.
The discussion of the interesting  case $n=2$ is the subject of a separate publication \cite{schreiber2}.

\section{Functors and local Trivializations}

\label{sec2}

We give the definition of the path groupoid of a smooth manifold and describe
functors defined on it. We introduce local trivializations and descent
data of such functors. 

\subsection{The Path Groupoid of a smooth Manifold}

\label{sec10}

We start by setting up the basic definitions around the path groupoid of
a smooth manifold $M$. We use the conventions of \cite{caetano,mackaay1},
generalized from loops to paths. 

\begin{definition}
\label{def9}
A \emph{path} $\gamma:x \to y$
between two points $x,y\in M$  is a smooth map $\gamma:[0,1] \to M$ which has a sitting instant:
a number $0 < \epsilon < \frac{1}{2}$ such that $\gamma(t)=x$
for $0 \leq t < \epsilon$ and $\gamma(t)=y$ for $1-\epsilon<t \leq 1$. 
\end{definition}

Let us denote the set of such paths by $PM$. For example, for any point $x\in M$ there is the constant path $\id_x$ defined
by $\id_x(t):=x$. Given a path $\gamma_1:x \to y$ and another path $\gamma_2:y
\to z$ we define their composition to be the path $\gamma_2 \circ \gamma_1:
x \to z$ defined by
\begin{equation*}
(\gamma_2 \circ \gamma_1)(t) := \begin{cases} \gamma_1(2t) & \text{for }
0 \leq t \leq \frac{1}{2} \\
\gamma_2(2t-1) & \text{for } \frac{1}{2}\leq t \leq 1\text{.} \\
\end{cases}
\end{equation*}
This gives a smooth map since $\gamma_1$ and $\gamma_2$ are both constant
near the gluing point, due to their
sitting instants. 
We also define the inverse $\gamma^{-1}:y
\to x$ of a path
$\gamma:x \to y$ by $\gamma^{-1}(t) := \gamma(1-t)$. 

\begin{definition}
\label{def10}
Two paths $\gamma_1:x \to y$ and $\gamma_2:x \to y$ are called \emph{thin homotopy equivalent},
if there exists a smooth map $h: [0,1] \times [0,1] \to M$ such that
\begin{enumerate}
\item
there exists a number $0 < \epsilon < \frac{1}{2}$ with 
\begin{enumerate}
\item 
$h(s,t)=x$ for $0 \leq t <\epsilon$ and $h(s,t)=y$ for $1-\epsilon< t
\leq 1$.
\item
$h(s,t)=\gamma_1(t)$ for $0 \leq s < \epsilon$ and $h(s,t)=\gamma_2(t)$ for
$1-\epsilon<s\leq 1$. 
\end{enumerate}
\item
the differential of $h$ has at most rank 1 everywhere, i.e.
\begin{equation*}
\mathrm{rank}(\mathrm{d}h|_{(s,t)}) \leq 1
\end{equation*}
for all $(s,t) \in [0,1] \times [0,1]$. 
\end{enumerate}
\end{definition}
Due to condition (b), thin homotopy  defines  an equivalence relation on $PM$. The set of thin homotopy classes of paths is denoted by $P^1M$, and
the projection to classes is denoted by
\begin{equation*}
\mathrm{pr}: PM \to P^1M\text{.}
\end{equation*}
We
denote a thin homotopy class of a path $\gamma :x \to y$ by $\overline{\gamma}:x \to
y$. Notice
that thin homotopies include the following type of reparameterizations: let
$\beta: [0,1] \to [0,1]$ be a path $\beta: 0 \to 1$, in particular with $\beta(0)=0$ and
$\beta(1)=1$.  Then, for any path $\gamma:x \to y$, also $\gamma \circ\beta:x \to y$
is a path and
\begin{equation*}
h(s,t) := \gamma(t \beta(1-s) + \beta(t)\beta(s))
\end{equation*}
defines a thin homotopy between them.

The composition
of paths defined above on $PM$ descends to $P^1M$ due to condition (a), which
admits a smooth composition of smooth homotopies. The composition of thin homotopy
classes of paths obeys the following rules:

\begin{lemma}
\label{lem5}
For any path $\gamma:x \to y$, 
\begin{enumerate}
\item[a)] 
$\overline{\gamma} \circ \overline{\id_x} = \overline{\gamma} = \overline{\id_y} \circ \overline{\gamma}$, 
\item[b)]
For further paths $\gamma':y \to z$ and $\gamma'':
z \to w$,
\begin{equation*}
(\overline{\gamma''} \circ \overline{\gamma'}) \circ \overline{\gamma} = \overline{\gamma''} \circ (\overline{\gamma'} \circ \overline{\gamma})\text{.}
\end{equation*}
\item[c)]
 $\overline{\gamma} \circ \overline{\gamma^{-1}}
= \overline{\id_y}$ and $\overline{\gamma^{-1}} \circ \overline{\gamma} = \overline{\id_x}$.
\end{enumerate}
\end{lemma}

These three properties lead us to the following

\begin{definition}
\label{def8}
For a smooth manifold $M$, we consider the category whose set of objects is $M$,
whose set of morphisms is $P^1M$, where a class $\overline{\gamma}:x
\to y$ is a morphism from $x$ to $y$, and the composition is
as described above. Lemma \ref{lem5} a) and b) are the axioms of a category
and c) says that every morphism is invertible. 
 Hence we have defined a groupoid, called the \emph{path groupoid} of $M$, and denoted
 by $\mathcal{P}_1(M)$. 
 
\end{definition}
For a smooth map $f:M \to N$, we denote by
\begin{equation*}
f_{*}: \mathcal{P}_1(M) \to \mathcal{P}_1(N)
\end{equation*}
the functor with $f_{*}(x)=f(x)$ and $(f_{*})(\overline{\gamma}):=\overline{f \circ \gamma}$. The latter is well-defined, since a thin homotopy $h$ between
paths $\gamma$ and $\gamma'$ induces a thin homotopy $f \circ h$ between
$f \circ \gamma$ and $f \circ \gamma'$. 

\medskip

In the following we consider functors
\begin{equation}
\label{44}
F: \mathcal{P}_1(M) \to T
\end{equation}
for some arbitrary category $T$. Such a functor sends each point $p\in M$ to an object $F(p)$ in
$T$, and each thin homotopy class $\overline{\gamma}:x \to y$ of paths to a morphism
$F(\overline{\gamma}): F(x) \to F(y)$ in $T$. We use the following notation: we call $M$ the \emph{base space} of the functor $F$, and the object $F(p)$ the \emph{fibre} of $F$
over $p$. In the remainder of this section we give examples of natural constructions
with functors (\ref{44}).

\paragraph*{Additional Structure on $T$.}

Any additional structure for the category $T$ can be
applied pointwise to functors into $T$, for instance,
\begin{itemize}
\item[a)]
if $T$ has direct sums, we can take the direct sum $F_1 \oplus F_2$ of two
functors. 
\item[b)] 
if $T$ is a monoidal category, we can take tensor products $F_1 \otimes F_2$
of functors. 
\item[c)]
if $T$ is monoidal and has a duality regarded as  a functor $d:T
\to T^{\mathrm{op}}$, we can form the dual $F^{*}:=d \circ F$ of a functor
$F$. 
\end{itemize}

\paragraph*{Pullback.}

If $f: M \to N$ is a smooth map and $F:\mathcal{P}_1(N) \to
T$ is a functor, we define
\begin{equation*}
f^{*}F := F \circ f_{*}: \mathcal{P}_1(M) \to T
\end{equation*}
to be the pullback of $F$ along $f$. 

\paragraph*{Flat Functors.}

\label{p1}
Instead of the path groupoid,  one can also consider the fundamental groupoid $\Pi_1(M)$ of a smooth
manifold $M$, whose objects are points in $M$, just like for $\mathcal{P}_1(M)$,
but whose morphisms are smooth homotopy classes of paths (whose differential
may have arbitrary rank). The projection
from thin homotopy classes to smooth homotopy classes provides a functor
\begin{equation*}
p: \mathcal{P}_1(M) \to \Pi_1(M)\text{.}
\end{equation*}
We  call a functor $F: \mathcal{P}_1(M) \to T$  \emph{flat}, if there exists a functor
$\tilde F: \Pi_1(M) \to T$ with $F \cong \tilde F \circ p$.
This is motivated by parallel transport in
principal $G$-bundles: while it is invariant under \textit{thin} homotopy, it is only homotopy invariant if the bundle is flat, i.e. has vanishing
curvature. However, aside from Section \ref{sec6_2} we will not discuss the flat case any further in this article.

\paragraph*{Restriction to Paths between fixed Points.}

Finally, let us consider the restriction
of a functor $F:\mathcal{P}_1(M) \to T$ to paths between two fixed
points. This yields a map
\begin{equation*}
F_{x,y}: \mathrm{Mor}_{\mathcal{P}_1(M)}(x,y)
\to \mathrm{Mor}_T(F(x),F(y))\text{.}
\end{equation*} 
Of particular interest is the case
$x=y$, in which $\mathrm{Mor}_{\mathcal{P}_1(M)}(x,x)$
forms a group under composition,
which is called the thin homotopy
group of $M$ at $x$, and is denoted
by $\pi_1^1(M,x)$ \cite{caetano,mackaay1}.
Even more particular, we consider
the target category $G\text{-}\mathrm{Tor}$: by choosing a diffeomorphism
$F(x)\cong G$, we obtain an
identification 
\begin{equation*}
\mathrm{Mor}_{G\text{-}\mathrm{Tor}}(F(x),F(x))= G\text{,}
\end{equation*}
and the restriction $F_{x,x}$ of a functor $F: \mathcal{P}_1(M) \to
G\text{-}\mathrm{Tor}$  to the thin homotopy
group of $M$ at $x$ gives a group
homomorphism
\begin{equation*}
F_{x,x}: \pi_1^1(M,x) \to G\text{.}
\end{equation*}
This way one obtains the setup of 
\cite{barret1,caetano} and \cite{mackaay1}
for the case $G=U(1)$ as a particular case of our setup.
A further question is, whether  the group homomorphism $F_{x,x}$ is smooth
in the sense used in \cite{barret1,caetano,mackaay1}. An answer is given in Section \ref{sec9}.

\subsection{Extracting Descent Data from a Functor}

\label{sec6}

To define local trivializations
of a functor $F: \mathcal{P}_1(M)
\to T$, we fix three attributes:
\begin{enumerate}
\item 
A surjective
submersion $\pi:Y \to M$. Compared to local trivializations
of fibre bundles, the surjective
submersion replaces an open cover
of the manifold. Indeed, given an open cover $\lbrace
U_{\alpha} \rbrace_{\alpha \in A}$ of $M$, one obtains  a surjective
submersion by taking $Y$ to be
the disjoint union of the $U_{\alpha}$
and $\pi:Y \to M$ to be the union of  the inclusions $U_{\alpha} \hookrightarrow M$. 
\item
A Lie
groupoid $\mathrm{Gr}$, i.e. a groupoid
whose sets of objects and morphisms
are smooth manifolds, whose source and target maps
\begin{equation*}
s,t: \mathrm{Mor}(\mathrm{Gr})
\to \mathrm{Obj}(\mathrm{Gr})
\end{equation*}
are surjective submersions, and  whose composition  
\begin{equation*}
\circ: \mathrm{Mor}(\mathrm{Gr})
\;_s\!\!\times_{t}\;\mathrm{Mor}(\mathrm{Gr})
\to \mathrm{Mor}(\mathrm{Gr})
\end{equation*}
and the identity $\id: \mathrm{Obj}(\mathrm{Gr})
\to \mathrm{Mor}(\mathrm{Gr})$
are smooth maps. The Lie groupoid $\mathrm{Gr}$ plays the role of the typical
fibre of the functor $F$.
\item
A  functor $i:\mathrm{Gr} \to T$, which relates the typical fibre $\mathrm{Gr}$
to the target category $T$ of the functor $F$. In  all of our examples, $i$ will be an equivalence
of categories. This is important for some results derived in Section \ref{sec11}.
\end{enumerate}

\begin{definition}
\label{def5}
Given a Lie groupoid $\mathrm{Gr}$, a functor $i: \mathrm{Gr} \to T$ and a surjective submersion
$\pi:Y \to M$, a \emph{$\pi$-local $i$-trivialization} of a functor 
\begin{equation*}
F: \mathcal{P}_1(M)
\to T
\end{equation*}
is a pair $(\mathrm{triv},t)$ of a functor $\mathrm{triv}: \mathcal{P}_1(Y) \to \mathrm{Gr}$ and
 a natural equivalence
\begin{equation*}
t:  \pi^{*}F \to i \circ \mathrm{triv}\text{.}
\end{equation*}
\end{definition}
The natural equivalence $t$ is also depicted by the diagram
\begin{equation*}
\alxydim{@C=1.5cm@R=1.5cm}{\mathcal{P}_1(Y) \ar[r]^{\pi_{*}} \ar[d]_{\mathrm{triv}} & \mathcal{P}_1(M)
\ar@{=>}[dl]|*+{t}
\ar[d]^{F} \\  \mathrm{Gr} \ar[r]_{i} & T\text{.}}
\end{equation*}
To set up the familiar terminology,
we call a functor \emph{locally $i$-trivializable},
if it admits a $\pi$-local $i$-trivialization
for some choice of $\pi$. We call
a functor \emph{$i$-trivial}, if it admits
an $\id_M$-local $i$-trivialization,
i.e. if it is naturally equivalent
to the functor $i \circ \mathrm{triv}$.
To abbreviate the notation, we will often write $\mathrm{triv}_i$ instead
of $i \circ \mathrm{triv}$. 

Note
that local trivializations can
be pulled back: if $\zeta: Z \to Y$  and $\pi:Y \to M$ are  surjective submersions,
and $(\mathrm{triv},t)$ is a $\pi$-local $i$-trivialization of a functor
$F$, we obtain a $(\pi \circ \zeta)$-local $i$-trivialization $(\zeta^{*}\mathrm{triv},\zeta^{*}t)$
of $F$. In terms of open covers, this corresponds to a refinement of the
cover. 

\begin{definition}
\label{def6}
Let $\mathrm{Gr}$ be a Lie groupoid and let $i: \mathrm{Gr} \to T$ be a functor. The category $\loctrivfunct{i}{\mathrm{Gr}}{1}$ of \emph{functors with $\pi$-local $i$-trivialization}
is defined as follows:
\begin{enumerate}
\item
its objects are triples $(F,\mathrm{triv},t)$ consisting of 
a functor $F: \mathcal{P}_1(M) \to T$
and a  $\pi$-local $i$-trivialization $(\mathrm{triv},t)$ of $F$.
\item
a morphism 
\begin{equation*}
\alxy{
(F,\mathrm{triv},t) \ar[r]^-\alpha &(F',\mathrm{triv}',t')
}
\end{equation*}
is a natural transformation $\alpha: F \to F'$.
Composition of morphisms is simply composition of these natural transformations.
\end{enumerate}
\end{definition}

\medskip

Motivated by transition functions
of fibre bundles, we extract a similar datum
from a functor $F$ with $\pi$-local $i$-trivialization $(\mathrm{triv},t)$; this datum is a natural equivalence
\begin{equation*}
g: \pi_1^{*}\mathrm{triv}_i \to
\pi_2^{*}\mathrm{triv}_i
\end{equation*}
between the two functors $\pi_1^{*}\mathrm{triv}_i$
and $\pi_2^{*}\mathrm{triv}_i$ from
$\mathcal{P}_1(Y^{[2]})$ to $T$,
where $\pi_1$ and $\pi_2$ are the
 projections from the two-fold
fibre product $Y^{[2]}:=Y \times_M
Y$ of $Y$ to the components. In the case that the
surjective submersion comes from
an open cover of $M$, $Y^{[2]}$
is the disjoint union of all two-fold
intersections of open subsets.
 The natural equivalence
$g$
is defined by
\begin{equation*}
  g := \pi_2^* t \circ \pi_1^* t^{-1}\text{;}
\end{equation*}
its component at a point $\alpha\in Y^{[2]}$
is the morphism $t(\pi_2(\alpha)) \circ t(\pi_1(\alpha))^{-1}$ in
$T$. The composition is well-defined because $\pi \circ \pi_1=\pi \circ \pi_2$. 

Transition functions of fibre bundles
satisfy a cocycle condition over
three-fold intersections. The
natural equivalence $g$ has a similar
property when pulled back to the
three-fold fibre product $Y^{[3]}:=Y
\times_M Y \times_M Y$. 

\begin{proposition}
The diagram
\begin{equation*}
  \alxydim{@R=1.5cm}{
    & \pi_2^*\mathrm{triv}_i
    \ar[dr]^{\pi_{23}^* g}
    \\
    \pi_1^* \mathrm{triv}_i 
    \ar[ur]^{\pi_{12}^* g}
    \ar[rr]_{\pi_{13}^* g}
    && 
    \pi_3^* \mathrm{triv}_i
  }
\end{equation*}
of natural equivalences between functors from $\mathcal{P}_1(Y^{[3]})$
to $T$ is commutative.
\end{proposition}

Now that we have defined the
data $(\mathrm{triv},g)$ associated to an object $(F,\mathrm{triv},t)$
in $\loctrivfunct{i}{\mathrm{Gr}}{1}$, we consider
a morphism 
\begin{equation*}
\alpha: (F,\mathrm{triv},t) \to (F',\mathrm{triv}',t')
\end{equation*}
between two functors with $\pi$-local $i$-trivializations, i.e. a natural
transformation $\alpha:F \to F'$. We define a natural transformation
\begin{equation*}
h: \mathrm{triv}_i \to \mathrm{triv}'_i
\end{equation*}
by $h := t' \circ \pi^{*}\alpha \circ t^{-1}$, whose component at $x \in Y$
is the morphism $t'(x) \circ \alpha(\pi(x)) \circ t(x)^{-1}$ in $T$. From
the definitions of $g$, $g'$ and $h$ one obtains the commutative diagram
\begin{equation}
\label{32}
 \alxydim{@C=1.5cm@R=1.5cm}{
   \pi_1^* \mathrm{triv}_i
   \ar[r]^g
   \ar[d]_{\pi_1^* h}
   &
   \pi_2^* \mathrm{triv}_i
   \ar[d]^{\pi_2^* h}
   \\
   \pi_1^* \mathrm{triv}'_i
   \ar[r]_{g'}
   &
   \pi_2^* \mathrm{triv}'_i\text{.}
 }
\end{equation}

The behaviour of the natural equivalences data $g$ and $h$ leads  to  the following definition of a category $\trans{i}{\mathrm{Gr}}{1}$ of descent data. This
terminology will be explained in the next section.

\begin{definition}
\label{def7}
The category $\trans{i}{\mathrm{Gr}}{1}$ of \emph{descent data} of $\pi$-locally $i$-trivialized
functors is defined as follows:
\begin{enumerate}
\item 
its objects are pairs $(\mathrm{triv},g)$ of a
functor $\mathrm{triv}: \mathcal{P}_1(Y) \to \mathrm{Gr}$ and a natural equivalence
\begin{equation*}
g: \pi_1^{*}\mathrm{triv}_{i} \to \pi_2^{*}\mathrm{triv}_{i}\text{,}
\end{equation*}
such that
the diagram
\begin{equation}
\label{200}
  \alxydim{@R=1.5cm}{
    & \pi_2^*\mathrm{triv}_i
    \ar[dr]^{\pi_{23}^* g}
    \\
    \pi_1^* \mathrm{triv}_i 
    \ar[ur]^{\pi_{12}^* g}
    \ar[rr]_{\pi_{13}^* g}
    && 
    \pi_3^* \mathrm{triv}_i
  }
\end{equation}
is commutative. 
\item
a morphism $(\mathrm{triv},g) \to (\mathrm{triv}',g')$ is a natural transformation
\begin{equation*}
h: \mathrm{triv}_{i} \to \mathrm{triv}'_{i}
\end{equation*}
such that
the diagram
\begin{equation}
\label{201}
 \alxydim{@C=1.5cm@R=1.5cm}{
   \pi_1^* \mathrm{triv}_i
   \ar[r]^g
   \ar[d]_{\pi_1^* h}
   &
   \pi_2^* \mathrm{triv}_i
   \ar[d]^{\pi_2^* h}
   \\
   \pi_1^* \mathrm{triv}'_i
   \ar[r]_{g'}
   &
   \pi_2^* \mathrm{triv}'_i\text{.}
 }
\end{equation}
is commutative. The composition is the composition of these natural transformations.
\end{enumerate}
\end{definition}
Summarizing, we have defined a functor
\begin{equation}
\label{4}
\ex{\pi}:\loctrivfunct{i}{\mathrm{Gr}}{1} \to \trans{i}{\mathrm{Gr}}{1}\text{,}
\end{equation}
that extracts descent data from functors
with local trivialization and of morphisms of those in the way described above.  

\subsection{Reconstructing a Functor from Descent Data}

\label{sec5}

In this section we show that extracting descent data from a functor $F$ preserves all  information about $F$. We also justify the terminology \emph{descent data}, see  Remark
\ref{re1} below.  

\begin{theorem}
\label{th3}
The functor
\begin{equation*}
\ex{\pi}: \loctrivfunct{i}{\mathrm{Gr}}{1} \to \trans{i}{\mathrm{Gr}}{1}
\end{equation*}
is an equivalence of categories. 
\end{theorem}
For the proof we define a weak inverse functor
\begin{equation}
\label{17}
\mathrm{Rec}_{\pi}: \trans{i}{\mathrm{Gr}}{1} \to \loctrivfunct{i}{\mathrm{Gr}}{1}
\end{equation}
that reconstructs a functor (and a $\pi$-local $i$-trivialization) from given descent data. The definition of  $\mathrm{Rec}_{\pi}$ is given in three steps: \begin{enumerate}
\item 
We construct a groupoid $\upp{M}{\pi}$ covering the path groupoid $\mathcal{P}_1(M)$
by means of a surjective functor $\uppp: \upp{M}{\pi}
\to \mathcal{P}_1(M)$, and show that any object $(\mathrm{triv},g)$ in $\trans{i}{\mathrm{Gr}}{1}$
gives rise to a functor
\begin{equation*}
R_{(\mathrm{triv},g)}: \upp{M}{\pi} \to T\text{.}
\end{equation*}
We enhance this to a functor
\begin{equation}
\label{24}
R:\trans{i}{\mathrm{Gr}}{1} \to \funct(\upp{M}{\pi},T)\text{,}
\end{equation}
where \label{p4} $\funct(\upp{M}{\pi},T)$ is the category of functors from $\upp{M}{\pi}$
to $T$ and natural transformations between those.

\item
We show that the  functor $\uppp: \upp{M}{\pi} \to \mathcal{P}_1(M)$
is an equivalence of categories and construct a weak inverse
\begin{equation*}
s: \mathcal{P}_1(M) \to \upp{M}{\pi}\text{.}
\end{equation*}
The pullback along $s$ is  the functor
\begin{equation}
\label{25}
s^{*}:\funct(\upp{M}{\pi},T) \to \funct(\mathcal{P}_1(M),T)
\end{equation}
obtained by pre-composition with $s$.
\item

By constructing canonical $\pi$-local $i$-trivializations of functors in the image of the composition $s^{*} \circ R$ 
of the functors (\ref{24}) and (\ref{25}), we extend this composition to a functor
\begin{equation*}
\mathrm{Rec}_{\pi}:=s^{*}\circ R:\trans{i}{\mathrm{Gr}}{1} \to \loctrivfunct{i}{\mathrm{Gr}}{1}\text{.}
\end{equation*}
Finally, we give in Appendix \ref{appA} the proof that $\mathrm{Rec}_{\pi}$
is a weak inverse of the functor $\ex{\pi}$ and thus show that $\ex{\pi}$
is an equivalence of categories.
\end{enumerate}

Before we perform the steps 1 to 3, let us make the following remark about
the nature of the category $\trans{i}{i}{1}$ and the functor $\mathrm{Rec}_{\pi}$.

\begin{remark}
\label{re1}
\normalfont
We consider the case $i := \id_\mathrm{Gr}$. Now, the forgetful functor $v: \loctrivfunct{i}{\mathrm{Gr}}{1} \to \funct(\mathcal{P}_{1}(M),\mathrm{Gr})$
has a canonical weak inverse, which associates to a functor $F:\mathcal{P}_1(M) \to
\mathrm{Gr}$ the $\pi$-local $i$-trivialization $(\pi^{*}F,\id_{\pi^{*}F})$.
Under this identification, $\trans{i}{\mathrm{Gr}}{1}$
is the descent category of the functor category $\funct(M,\mathrm{Gr})$
with respect to $\pi$ in the sense of a stack \cite{moerdijk1,street}. The functor 
\begin{equation*}
\mathrm{Rec}_{\pi}:\trans{i}{\mathrm{Gr}}{1} \to \funct(\mathcal{P}_1(M),\mathrm{Gr}) \end{equation*}
realizes the descent.
\end{remark}

\paragraph*{Step 1: The Groupoid $\upp{M}{\pi}$.}
The groupoid $\upp{M}{\pi}$ we introduce  
is 
the \emph{universal path pushout} 
associated to the surjective submersion $\pi: Y \to M$. Heuristically, $\upp{M}{\pi}$ is the path groupoid of the covering
$Y$ combined with \quot{jumps} in the fibres of $\pi$. We explain its universality
in Appendix \ref{app2}  for completeness and introduce here a concrete realization (see Lemma \ref{lem3}).

\begin{definition}
\label{def21}
The groupoid $\upp{M}{\pi}$ is defined as follows. Its objects 
are points $x\in Y$ and its morphisms are  formal (finite) compositions of   two types
of basic morphisms: thin homotopy classes $\overline{\gamma}:x
\to y$ of paths in $Y$, and 
points $\alpha\in Y^{[2]}$ regarded as  morphisms $\alpha:\pi_1(\alpha) \to \pi_2(\alpha)$. Among the morphisms, 
we impose three relations:
\begin{enumerate}
\item[(1)]
for any thin homotopy class $\overline{\Theta}:\alpha \to \beta$ of paths  in $Y^{[2]}$, we demand that the diagram
\begin{equation*}
 \alxydim{@C=1.5cm@R=1.5cm}{
   \pi_1(\alpha)
   \ar[r]^{\alpha}
   \ar[d]_{(\pi_1)_{*}(\overline{\Theta})}
   &
   \pi_2(\alpha)
   \ar[d]^{(\pi_2)_{*}(\overline{\Theta})}
   \\
   \pi_1(\beta)
   \ar[r]_\beta
   &
   \pi_2(\beta)\text{.}   
 }
\end{equation*}
of morphisms in $\upp{M}{\pi}$ is commutative.
\item[(2)] 
for any point $\Xi \in Y^{[3]}$, we demand that the diagram
\begin{equation*}
 \alxydim{@R=1.5cm}{
   & \pi_2(\Xi)
   \ar[dr]^{\pi_{23}(\Xi)}
   \\
   \pi_1(\Xi) 
   \ar[ur]^{\pi_{12}(\Xi)}
   \ar[rr]_{\pi_{13}(\Xi)}
   && 
   \pi_3(\Xi)
 }
\end{equation*}
of morphisms in $\upp{M}{\pi}$ is commutative.
\item[(3)]
we impose the equation
$\overline{\id_x}
= (x,x)\in Y^{[2]}$ for any $x\in Y$.
\end{enumerate}
\end{definition}

It is clear that this definition  indeed gives a groupoid. It is important for us because it
provides   the two following natural definitions.
\begin{definition}
\label{def13}
For an object $(\mathrm{triv},g)$ in $\trans{i}{\mathrm{Gr}}{1}$, we have a functor
\begin{equation*}
R_{(\mathrm{triv},g)}: \upp{M}{\pi} \to T
\end{equation*}
that sends  an object $x \in Y$ to $\mathrm{triv}_i(x)$, a basic morphism $\overline{\gamma}:x
\to y$ to $\mathrm{triv}_i(\overline{\gamma})$ and a basic morphism $\alpha$ to $g(\alpha)$.
\end{definition}

The definition is well-defined since it respects  the relations among the morphisms:
(1) is respected due to the commutative diagram for  the natural transformation $g$, (2) is the
cocycle condition (\ref{200}) for $g$ and (3)  follows from the latter since $g$ is
invertible. 

\begin{definition}
\label{def14}
For a morphism 
$h:(\mathrm{triv},g) \to (\mathrm{triv}',g')$
in $\trans{i}{\mathrm{Gr}}{1}$ we have a
natural transformation
\begin{equation*}
R_h : R_{(\mathrm{triv},g)} \to R_{(\mathrm{triv}',g')}
\end{equation*}
that sends an object $x\in Y$ to the morphism $h(x)$ in $T$. 
\end{definition}

The commutative diagram for the natural transformation $R_h$
for a basic morphism $\overline{\gamma}:x \to y$ follows from the one of $h$, and  for a basic morphism $\alpha\in Y^{[2]}$  from the condition (\ref{201})
on the morphisms of $\trans{i}{\mathrm{Gr}}{1}$. 

We explain in Appendix \ref{app2} that Definitions (\ref{def13}) and (\ref{def14}) 
are consequences
of the universal property of the groupoid $\upp{M}{\pi}$, as specified
in Definition \ref{def4} and calculated  in Lemma \ref{lem3}. Here we summarize the definitions above in the following
way:
\begin{lemma}
Definitions (\ref{def13}) and (\ref{def14})
yield a functor
\begin{equation}
\label{16}
R: \trans{i}{\mathrm{Gr}}{1} \to \funct(\upp{M}{\pi},T)\text{.}
\end{equation}
\end{lemma}

\medskip

\paragraph*{Step 2: Pullback to $M$.}
To continue the reconstruction of a functor from given
descent data
let us  introduce
the  projection functor
\begin{equation}
\label{51}
\uppp: \upp{M}{\pi} \to \mathcal{P}_1(M)
\end{equation}
sending an object $x \in Y$ to $\pi(x)$, a basic morphism $\overline{\gamma}:x \to y$ to $\pi_{*}(\overline{\gamma})$
and a basic morphism $\alpha\in Y^{[2]}$ to $\id_{\pi(\pi_1(\alpha))}$ ($=\id_{\pi(\pi_2(\alpha))}$).
In other words, it is just the functor $\pi_{*}$ and forgets the jumps in the
fibres of $\pi$. 
More precisely,
\begin{equation*}
\uppp \circ \iota =\pi_{*}\text{,}
\end{equation*}
where $\iota:\mathcal{P}_1(Y) \to \upp{M}{\pi}$ is the obvious
inclusion functor.

\begin{lemma}
\label{lem8}
The projection functor $\uppp:\upp{M}{\pi} \to \mathcal{P}_1(M)$ is a surjective equivalence of categories.
\end{lemma}

\begin{proof}
Since $\pi:Y \to M$ is surjective, it is clear that $\uppp$ is surjective on objects. It remains to show that  the map
\begin{equation}
\label{1}
(\uppp)_1: \mathrm{Mor}_{\upp{M}{\pi}}(x,y) \to \mathrm{Mor}_{\mathcal{P}_1(M)}(\pi(x),\pi(y))
\end{equation}
is bijective for all $x,y\in Y$. Let $\gamma: \pi(x) \to \pi(y)$ be any path in $M$. Let $\lbrace U_i \rbrace_{i\in I}$ an open cover of $M$ with sections $s_i:U_i \to Y$. Since the image of $\gamma:[0,1] \to M$ is compact, there exists a \textit{finite} subset $J \subset I$ such that $\lbrace U_i \rbrace_{i\in J}$ covers the image. Let $\gamma = \gamma_n \circ ... \circ \gamma_1$ be a decomposition of $\gamma$ such that $\gamma_i \in PU_{j(i)}$ for some assignment $j: \lbrace 1,...,n\rbrace \to J$. Let $\tilde\gamma_i := (s_{j(i)})_{*}\gamma_i \in PY$ be lifts of the pieces, $\tilde\gamma_i: a_i \to b_i$ with $a_i,b_i \in Y$. Now we consider the path
\begin{equation*}
\tilde\gamma := (b_n,y) \circ \tilde\gamma_n \circ (b_{n-1},a_n) \circ ... \circ \tilde\gamma_2 \circ  (b_1,a_2) \circ \tilde\gamma_1 \circ (x,a_1)\text{,}
\end{equation*}
whose thin homotopy class is evidently a preimage of the thin homotopy class of $\gamma$ under $(\uppp)_1$. The injectivity of (\ref{1}) follows from the identifications (1), (2) and (3) of morphisms in the groupoid $\upp{M}{\pi}$.  
\end{proof}

Since $\uppp$ is an equivalence of categories, there exists a (up to natural isomorphism) unique weak inverse functor $s: \mathcal{P}_1(M) \to \upp{M}{\pi}$ together with natural equivalences $\lambda : s \circ \uppp \to \id_{\upp{M}{\pi}}$ and $\rho: \uppp \circ s \to \id_{\mathcal{P}_1(M)}$. The inverse functor $s$ can be constructed explicitly: for a fixed choice of lifts $s(x) \in Y$ for every point $x\in M$, and a fixed choice of an open cover, each path can be lifted as described in the proof of Lemma \ref{lem8}. In this case we  have $\rho=\id$, and the component of $\lambda$ at $x\in Y$ is the morphism $(s(\pi(x)),x)$ in $\upp{M}{\pi}$. Now we have a canonical functor
\begin{equation*}
s^{*} \circ R: \trans{i}{\mathrm{Gr}}{1} \to \funct(\mathcal{P}_1(M),T)\text{.}
\end{equation*}
It reconstructs a functor $s^{*}R_{(\mathrm{triv},g)}$ from a given object
$(\mathrm{triv},g)$ in $\trans{i}{\mathrm{Gr}}{1}$ and a natural transformation $s^{*}R_h$
from a given morphism $h$ in $\trans{i}{\mathrm{Gr}}{1}$. 

\paragraph*{Step 3: Local Trivialization.} What remains   to enhance the functor $s^{*}
\circ R$ to a functor 
\begin{equation*}
\mathrm{Rec}_{\pi}:\trans{i}{\mathrm{Gr}}{1} \to \loctrivfunct{i}{\mathrm{Gr}}{1}
\end{equation*}
is finding a $\pi$-local
$i$-trivialization $(\mathrm{triv},t)$ of each reconstructed functor $s^{*}R_{(\mathrm{triv},g)}$.
Of course  the given functor $\mathrm{triv}:\mathcal{P}_1(Y)
\to \mathrm{Gr}$ serves as the first component of the trivialization, and
it remains to define the natural equivalence
\begin{equation}
\label{30}
t: \pi^{*}s^{*}R_{(\mathrm{triv},g)} \to \mathrm{triv}_i\text{.}
\end{equation}
We use the natural equivalence $\lambda : s \circ \uppp \to \id_{\upp{M}{\pi}}$
associated to the functor $s$ and obtain a natural equivalence
\begin{equation*}
\iota^{*}\lambda: s \circ \pi_{*} \to \iota
\end{equation*}
between functors from $\mathcal{P}_1(Y)$ to $\upp{M}{\pi}$. Its component at 
$x\in Y$ is the morphism $(s(\pi(x)),x)$ going from $s(\pi(x))$ to $x$. 
Using 
\begin{equation*}
\pi^{*}s^{*}R_{(\mathrm{triv},g)}=(s \circ \pi_{*})^{*}R_{(\mathrm{triv},g)}
\quad\text{ and }\quad\mathrm{triv}_i=\iota^{*}R_{(\mathrm{triv},g)}\text{,}
\end{equation*}
we define by
\begin{equation*}
t:= g \circ \iota^{*}\lambda
\end{equation*}
the natural equivalence (\ref{30}).
Indeed, its component at $x \in Y$ is the morphism $g((s(\pi(x)),x)): \mathrm{triv}_i(s(\pi(x)))
\to \mathrm{triv}_i(x)$, these are natural in $x$ and isomorphisms because $g$ is one. Diagrammatically, it is
\begin{equation*}
\alxy{\mathcal{P}_1(Y) \ar[dr]_{\iota} \ar[rr]^{\pi_{*}} \ar[dd]_{\mathrm{triv}} 
    && \mathcal{P}_1(M) \ar@{=>}[dl]|*+{\lambda} \ar@/^1pc/[dl]^{s}
\ar[dd]^{s^{*}R_{(\mathrm{triv},g)}} 
   \\ 
   & \upp{M}{\pi} \ar@/^1pc/[ur]^<<<<{\uppp} \ar[dr]_{R_{(\mathrm{triv},g)}}
& \\ \mathrm{Gr} \ar[rr]_{i} && T\text{.}}
\end{equation*}
This shows
\begin{lemma}
The pair $(\mathrm{triv},t)$ is a $\pi$-local $i$-trivialization of the functor
$s^{*}R_{(\mathrm{triv},g)}$.
\end{lemma}

This finishes the definition of the reconstruction functor $\mathrm{Rec}_{\pi}$.
The remaining proof that $\mathrm{Rec}_{\pi}$ is a weak inverse of $\ex{\pi}$ is postponed to Appendix \ref{appA}.

\section{Transport Functors}

\label{sec3}

Transport functors are locally trivializable functors whose descent data
is smooth. Wilson lines are restrictions of a  functor to
paths between two fixed points. We deduce a characterization of  transport functors
by the smoothness of their Wilson lines.

\subsection{Smooth Descent Data}

In this section we specify a subcategory $\transsmooth{i}{1}{\mathrm{Gr}}$ of the category
 $\trans{i}{i}{1}$ of descent data we have defined in the previous section. This subcategory  is
 supposed to contain \emph{smooth} descent data.  The main issue 
 is  to decide, when a functor $F: \mathcal{P}_1(X) \to \mathrm{Gr}$ is smooth:
 in contrast to the objects and the morphisms of the Lie groupoid $\mathrm{Gr}$, the set $P^1X$ of morphisms of $\mathcal{P}_1(X)$ is not a smooth manifold.
 
\begin{definition}
\label{def1}
Let $\mathrm{Gr}$ be a Lie groupoid and let $X$ be a smooth manifold. A functor $F: \mathcal{P}_1(X) \to \mathrm{Gr}$ is called \emph{smooth}, if the following two
conditions are satisfied:
\begin{enumerate}
\item 
On objects, $F: X \to \mathrm{Obj}(\mathrm{Gr})$ is a smooth map.
\item
For every $k \in \N_0$, every open subset $U\subset \R^k$  and every map $c: U \to PX$ such that the composite
\begin{equation}
\label{31}
\alxydim{@C=1.5cm}{U \times [0,1] \ar[r]^-{c \times \id} & PX \times [0,1] \ar[r]^-{\mathrm{ev}}
& X}
\end{equation}
is smooth, also
\begin{equation*}
\alxy{U \ar[r]^-{c} & PX \ar[r]^{\mathrm{pr}} & P^1 X \ar[r]^-{F} & \mathrm{Mor}(\mathrm{Gr})}
\end{equation*}
is smooth.
\end{enumerate}
\end{definition}

In (\ref{31}), $\mathrm{ev}$ is the evaluation map $\mathrm{ev}(\gamma,t):=\gamma(t)$.
Similar definitions of smooth maps defined on thin homotopy classes of paths have also been used in \cite{barret1,caetano,mackaay1}.
We explain in Appendix \ref{appB} how Definition \ref{def1} is motivated and how
it arises from the  general concept of  diffeological spaces  \cite{chen1},
a generalization of the concept of a smooth manifold, cf. Proposition
\ref{prop4} i). 

\begin{definition}
\label{def15}
A natural transformation 
$\eta: F \to G$ between smooth functors $F,G:\mathcal{P}_1(X) \to \mathrm{Gr}$ is called \emph{smooth}, if its components form a smooth map $X \to \mathrm{Mor}(\mathrm{Gr}): X \mapsto \eta(X)$.
\end{definition}

Because the composition in the Lie groupoid $\mathrm{Gr}$ is smooth, compositions
of smooth natural transformations are again smooth.  Hence,
smooth functors and smooth natural transformations 
 form a category $\funct^{\infty}(\mathcal{P}_1(X),\mathrm{Gr})$. Notice that if
 $f:M \to X$ is a smooth map, and $F:\mathcal{P}_1(X) \to \mathrm{Gr}$ is
 a smooth functor, the pullback $f^{*}F$ is also smooth. Similarly, pullbacks
 of smooth natural transformations are smooth.

\begin{definition}
\label{def2}
Let $\mathrm{Gr}$ be a Lie groupoid and let $i:\mathrm{Gr} \to T$ be a functor. An object $(\mathrm{triv},g)$ in $\trans{i}{\mathrm{Gr}}{1}$ is called \emph{smooth},
if the  following two conditions are satisfied:
\begin{enumerate}
\item 
The functor 
\begin{equation*}
\mathrm{triv}:\mathcal{P}_1(Y) \to \mathrm{Gr}
\end{equation*}
is smooth in
the sense of Definition \ref{def1}.
\item
The natural equivalence 
\begin{equation*}
g: \pi_1^{*}\mathrm{triv}_i \to \pi_2^{*}\mathrm{triv}_i
\end{equation*}
factors through $i$ by a  natural equivalence $\tilde g: \pi_1^{*}\mathrm{triv}
\to \pi_2^{*}\mathrm{triv}$ which is smooth in the sense of Definition \ref{def15}. For the components at a point $\alpha\in Y^{[2]}$, the factorization means
$g(\alpha)=i(\tilde g(\alpha))$.
\end{enumerate}
In the same sense,
a morphism 
\begin{equation*}
h:(\mathrm{triv},g) \to (\mathrm{triv}',g')
\end{equation*}
between smooth objects is called smooth, if it factors through $i$ by
a smooth natural
equivalence $\tilde h: \mathrm{triv} \to \mathrm{triv}'$.
\end{definition}

\begin{remark}
If $i$ is  faithful,  the natural equivalences $\tilde g$ and $\tilde h$ in Definition \ref{def2} are uniquely determined, provided that they exist. If $i$ is additionally
 full, also the existence of $g$ and $h$ is guaranteed.
\end{remark}

Smooth objects and morphisms in $\trans{i}{\mathrm{Gr}}{1}$ form the subcategory $\transsmooth{i}{1}{\mathrm{Gr}}$.
Using the equivalence $\ex{\pi}$ defined in Section \ref{sec6}, we obtain a subcategory 
$\loctrivfunctsmooth{i}{\mathrm{Gr}}{\pi}{1}$ of $\loctrivfunct{i}{\mathrm{Gr}}{1}$ consisting of those objects
$(F,\mathrm{triv},t)$ for which $\ex{\pi}(F,\mathrm{triv},t)$ is smooth and of
those morphisms $h$ for which $\ex{\pi}(h)$ is smooth. 

\begin{proposition}
 The functor $\mathrm{Rec}_{\pi}: \trans{i}{\mathrm{Gr}}{1} \to \loctrivfunct{i}{\mathrm{Gr}}{1}$
 restricts to an equivalence
of categories
\begin{equation*}
\mathrm{Rec}_{\pi}: \transsmooth{i}{1}{\mathrm{Gr}} \to \loctrivfunctsmooth{i}{\mathrm{Gr}}{\pi}{1}\text{.}
\end{equation*}
\end{proposition}

\begin{proof}
This follows from the fact that $\ex{\pi} \circ \mathrm{Rec}_{\pi}=\id_{\trans{i}{\mathrm{Gr}}{1}}$,
see the proof of Theorem \ref{th3} in Appendix \ref{appA}.
\end{proof}

Now we are ready to define transport functors. 

\begin{definition}
\label{def3}
Let $M$ be a smooth manifold, $T$ a category, $\mathrm{Gr}$ a Lie groupoid and $i:\mathrm{Gr}
\to T$ a functor.
\begin{enumerate}
\item 
A \emph{transport functor on $M$ in $T$ with $\mathrm{Gr}$-structure} is a functor
\begin{equation*}
\mathrm{tra}: \mathcal{P}_1(M) \to T
\end{equation*}
such that there exists a surjective submersion $\pi:Y \to M$ and a $\pi$-local $i$-trivialization $(\mathrm{triv},t)$, such that $\ex{\pi}(\mathrm{tra},\mathrm{triv},t)$
is smooth. 

\item
A \emph{morphism between transport functors on $M$ in $T$ with $\mathrm{Gr}$-structure}
is a natural equivalence $\eta:\mathrm{tra}
\to \mathrm{tra}'$ such that there exists a surjective submersion $\pi:Y
\to M$ together with $\pi$-local $i$-trivializations of $\mathrm{tra}$ and
$\mathrm{tra}'$, such that $\ex{\pi}(\eta)$ is smooth. 
\end{enumerate}
\end{definition}

It is clear that the identity natural
transformation of a transport functor
$\mathrm{tra}$ is a morphism in
the above sense. To show that the
composition of morphisms between
transport functors is possible,
note that if $\pi:Y \to M$ is a
surjective submersion for which
$\ex{\pi}(\eta)$ is smooth, and
$\zeta:Z \to Y$ is another surjective
submersion, then also $\ex{\pi
\circ \zeta}(\eta)$ is smooth.
If now
$\eta: \mathrm{tra} \to \mathrm{tra}'$
and $\eta': \mathrm{tra}' \to \mathrm{tra}''$
are morphisms of transport functors,
 and
 $\pi:Y \to M$ and $\pi':Y' \to
 M$ are surjective submersions
 for which $\ex{\pi}(\eta)$ and
 $\ex{\pi'}(\eta')$ are smooth,
 the fibre product $\tilde\pi: Y \times_M
 Y' \to M$ is a surjective submersion
 and factors through $\pi$ and
 $\pi'$ by surjective submersions.
 Hence, $\ex{\tilde\pi}(\eta' \circ
  \eta)$ is smooth. 
  
\begin{definition}
\label{def16}
The category of all transport functors on $M$ in $T$ with $\mathrm{Gr}$-structure 
and all morphisms between those  is denoted by $\transport{i}{1}{\mathrm{Gr}}{T}$.
\end{definition}

From the definition of a transport functor with $\mathrm{Gr}$-structure  it is not clear that, for a fixed
surjective submersion $\pi:Y \to M$, all choices of  $\pi$-local $i$-trivializations $(\mathrm{triv},t)$ with smooth descent data give rise to isomorphic objects
in $\transsmooth{i}{1}{\mathrm{Gr}}$. This is at least true for full functors
$i: \mathrm{Gr} \to T$ and contractible surjective
submersions: a surjective submersion $\pi:Y \to M$ is called \textit{contractible},
if there exists a smooth map $c: Y \times [0,1] \to Y$ such that $c(y,0)=y$
for all $y\in Y$ and $c(y,1)=y_k$ for some fixed choice of $y_k\in Y_k$ for
each connected component $Y_k$ of $Y$. We may assume without loss of generality,
that $c$ has a sitting instant with respect to the second parameter, so that
we can regard $c$ also as a map $c:Y \to PY$. For example, if $Y$
is the disjoint union of the open sets of a good open cover of $M$, $\pi:Y
\to M$ is contractible.

\begin{lemma}
\label{prop3}
Let $i:\mathrm{Gr} \to T$ be a full functor,  let $\pi:Y \to M$ be a contractible surjective
submersion and let $(\mathrm{triv},t)$ and $(\mathrm{triv}',t')$
be two $\pi$-local $i$-trivializations
of a transport functor $\mathrm{tra}:\mathcal{P}_1(M)
\to T$ with $\mathrm{Gr}$-structure. Then, the identity natural transformation
$\id_{\mathrm{tra}}:\mathrm{tra} \to \mathrm{tra}$ defines a morphism
\begin{equation*}
\id_{\mathrm{tra}}: (\mathrm{tra},\mathrm{triv},t) \to (\mathrm{tra},\mathrm{triv}',t')
\end{equation*} 
in $\loctrivfunctsmooth{i}{\mathrm{Gr}}{\pi}{1}$, in particular, $\ex{\pi}(\mathrm{tra},\mathrm{triv},t)$ and $\ex{\pi}(\mathrm{tra},\mathrm{triv}',t')$ are isomorphic objects in $\transsmooth{i}{1}{\mathrm{Gr}}$.
\end{lemma}

\begin{proof}
Let $c : Y \times [0,1] \to Y$ be a smooth
contraction, regarded as a map $c : Y \to PY$. For each $y \in Y_k$ we have a path $c(y):y \to y_k$,
and the commutative diagram for the natural transformation $t$ gives 
\begin{equation*}
t(y)  = \mathrm{triv}_i(\overline{c(y)})^{-1}  \circ t(y_k)\circ \mathrm{tra}(\pi_{*}(\overline{c(y)}))\text{,}
\end{equation*}
and analogously for $t'$. The descent datum of the natural equivalence $\id_{\mathrm{tra}}$ is the natural equivalence
\begin{equation*}
h:=\ex{\pi}(\id)=t' \circ t^{-1} : \mathrm{triv}_i \to \mathrm{triv}'_i\text{.}
\end{equation*}
Its component at $y \in Y_k$ is the morphism
\begin{equation}
h(y)=\mathrm{triv}_i'(\overline{c(y)})^{-1} \circ t'(y_k) \circ  t(y_k)^{-1}\circ \mathrm{triv}_i(\overline{c(y)}): \mathrm{triv}_i(y) \to \mathrm{triv}'_i(y)
\end{equation}
in $T$. Since $i$ is full, $t'(y_k) \circ t(y_k)^{-1} = i(\kappa_k)$
for some morphism $\kappa_k: \mathrm{triv}(y_k) \to \mathrm{triv}'(y_k)$, so that
 $h$ factors through $i$ by
\begin{equation*}
\tilde h(y):=\mathrm{triv}'(\overline{c(y)})^{-1} \circ \kappa_k\circ \mathrm{triv}(\overline{c(y)})\in
\mathrm{Mor}(\mathrm{Gr})\text{.}
\end{equation*}
Since $\mathrm{triv}$ and $\mathrm{triv}'$ are smooth functors, $\mathrm{triv}\circ \mathrm{pr} \circ c$ and $\mathrm{triv}'\circ \mathrm{pr} \circ c$ are smooth maps, so that the components of $\tilde h$ form a smooth map $Y \to \mathrm{Mor}(\mathrm{Gr})$.
Hence, $h$ is a morphism
in $\transsmooth{i}{1}{\mathrm{Gr}}$. 
\end{proof}

To keep track of all the categories we have defined, consider the following
diagram of functors which is strictly commutative:
\begin{equation}
\label{33}
\alxydim{@C=1.5cm}{\transsmooth{i}{1}{\mathrm{Gr}} \ar@{^(->}[d] \ar[r]^{\mathrm{\mathrm{Rec}_{\pi}}} & \loctrivfunctsmooth{i}{\mathrm{Gr}}{\pi}{1} \ar@{^(->}[d] \ar[r]^-{v^{\infty}} & \transport{i}{1}{\mathrm{Gr}}{T} \ar@{^(->}[d] \\ \trans{i}{\mathrm{Gr}}{1} \ar[r]_{\mathrm{Rec}_{\pi}}
& \loctrivfunct{i}{\mathrm{Gr}}{1} \ar[r]_{v} & \funct(M,T)}
\end{equation}
The vertical arrows are the inclusion functors, and $v^{\infty}$ and $v$ are forgetful functors. In the next subsection we show that the functor $v^{\infty}$ is an equivalence
of categories.

\subsection{Wilson Lines of Transport Functors}

\label{sec11}

We restrict functors to paths
between two fixed points and study the smoothness of these
restrictions. For this purpose we assume that the functor $i:\mathrm{Gr}
\to T$ is an equivalence of categories; this is the case in all   examples
of transport functors we give in Section \ref{sec4}.

\begin{definition}
\label{def22}
Let $F:\mathcal{P}_1(M) \to T$
be a functor, let $\mathrm{Gr}$
be a Lie groupoid and let $i:\mathrm{Gr}
\to T$ be an equivalence
of categories. Consider two points $x_1,x_2\in M$
together with a choice of objects
$G_k$ in $\mathrm{Gr}$ and isomorphisms
$t_k: F(x_k) \to i(G_k)$ in $T$ for $k=1,2$. Then, the map
\begin{equation*}
\mathcal{W}_{x_1,x_2}^{F,i}: \mathrm{Mor}_{\mathcal{P}_1(M)}(x,y)
\to \mathrm{Mor}_{\mathrm{Gr}}(G_1,G_2):
\overline{\gamma} \mapsto i^{-1}(t_2
\circ F(\overline{\gamma})
\circ t_1^{-1})
\end{equation*}
is called the \emph{Wilson line of $F$
from $x_1$ to $x_2$}. 
\end{definition}
Note that because $i$ is essentially
surjective, the choices
of objects $G_k$ and morphisms $t_k: F(x_k) \to G_k$ exist for all points
$x_k \in M$. Because $i$ is full
and faithful, the morphism $t_2
\circ F(\overline{\gamma})
\circ t_1^{-1}: i(G_1) \to i(G_2)$ has a unique
preimage under $i$, which is the
Wilson line. For a different choice $t_k':F(x_k) \to i(G_k')$ of objects
in $\mathrm{Gr}$ and isomorphisms in $T$ the Wilson line changes like
\begin{equation*}
\mathcal{W}^{F,i}_{x_1,x_2} \mapsto \tau_2^{-1} \circ \mathcal{W}^{F,i}_{x_1,x_2} \circ \tau_1
\end{equation*}
for $\tau_k:G_k' \to G_k$ defined by $i(\tau_k)=t_k \circ t_k'^{-1}$. 
\begin{definition}
\label{def19}
A Wilson line
$\mathcal{W}^{F,i}_{x_1,x_2}$ is
called \emph{smooth}, if for every $k \in \N_0$, every open
subset $U \subset \R^k$ and every
map $c:U \to PM$ such  that $c(u)(t)\in M$ is smooth on
$U \times [0,1]$, $c(u,0)=x_1$ and $c(u,1)=x_2$ for all $u\in U$, also the map
\begin{equation*}
\mathcal{W}_{x_1,x_2}^{F,i} \circ
\mathrm{pr} \circ c: U \to \mathrm{Mor}_{\mathrm{Gr}}(G_1,G_2)
\end{equation*}
is smooth.
\end{definition}

This definition of smoothness arises
again from the context
of diffeological spaces, see Proposition \ref{prop5} i)
in Appendix \ref{appB}. Notice that if a Wilson line is smooth for some choice
of objects $G_k$ and isomorphisms $t_k$, it is smooth for any other choice.
For this reason we have not labelled Wilson lines with additional
indices $G_1,G_2,t_1,t_2$. 

\begin{lemma}
\label{lem12}
Let $i:\mathrm{Gr} \to T$ be an equivalence of categories, let 
\begin{equation*}
F:\mathcal{P}_1(M) \to T
\end{equation*}
be a functor whose Wilson lines $\mathcal{W}^{F,i}_{x_1,x_2}$
are smooth for all points $x_1,x_2\in M$, and let $\pi:Y \to M$ be a contractible
surjective submersion. Then, $F$ admits a $\pi$-local $i$-trivialization
$(\mathrm{triv},t)$
whose descent data  $\ex{\pi}(\mathrm{triv},t)$ is smooth. 
\end{lemma}

\begin{proof}
We choose a smooth contraction $r:Y \to PY$ and make, for every connected component $Y_k$ of $Y$,   a choice of objects
$G_k$ in $\mathrm{Gr}$ and isomorphisms $t_k:F(\pi(y_k)) \to i(G_k)$.
First we set $\mathrm{triv}(y) := G_k$ for all $y \in Y_k$, and define 
morphisms
\begin{equation*}
t(y):= t_k \circ F(\pi_{*}(\overline{r(y)})): F(\pi(y)) \to i(G_k)
\end{equation*}
in $T$. For a path $\gamma:y \to y'$, we define
the morphism
\begin{equation*}
\mathrm{triv}(\overline{\gamma}) := i^{-1}(t(y') \circ  F(\pi_{*}(\overline{\gamma}))
\circ t(y)^{-1}):G_k \to G_k
\end{equation*}
in $\mathrm{Gr}$. By construction, the morphisms $t(y)$ are the components of a natural equivalence
$t:\pi^{*}F \to \mathrm{triv}_i$, so that we have defined a $\pi$-local
$i$-trivialization $(\mathrm{triv},t)$ of $F$. 
Since $\mathrm{triv}$ is locally
constant on objects, it satisfies condition 1 of Definition \ref{def1}.
To check condition 2, notice that, for any path $\gamma:y \to y'$,
\begin{equation}
\label{47}
\mathrm{triv}(\overline{\gamma}) = \mathcal{W}^{F,i}_{y_k,y_k}(\pi_{*}(\overline{r(y')} \circ \overline{\gamma} \circ
\overline{r(y)}^{-1}))\text{.}
\end{equation} 
More generally, if $c:U \to PY$ is a map, we have, for every $u\in U$, a path
\begin{equation*}
\tilde c(u):= \pi_{*} (r(c(u)(1)) \circ c(u) \circ r(c(u)(0))^{-1})
\end{equation*}
in $M$. Then, equation (\ref{47}) becomes
\begin{equation*}
\mathrm{triv} \circ \mathrm{pr} \circ c = \mathcal{W}^{F,i}_{y_k,y_k} \circ
\mathrm{pr} \circ \tilde c\text{.} 
\end{equation*}
Since the right hand side is by assumption a smooth function;    $\mathrm{triv}$ is a smooth functor.
The component of the natural equivalence $g:=\pi_2^{*}t \circ \pi_1^{*}t^{-1}$ at a point $\alpha=(y,y')\in
Y^{[2]}$ with $y\in Y_k$ and $y'\in Y_l$ is the morphism
\begin{equation*}
g(\alpha) =t_l \circ F(\pi(c(y'))) \circ F(\pi(c(y)))^{-1} \circ t_k^{-1}: i(G_k) \to
i(G_l)\text{,} 
\end{equation*}
and hence of the form $g(\alpha)=i( \tilde g(\alpha))$. Now consider a chart
$\varphi: V \to Y^{[2]}$ with an open subset $V\in\R^n$, and the path $c(u):=r(\pi_2(\varphi(u))) \circ r(\pi_1(\varphi(u)))^{-1}$
in $Y$. We find
\begin{equation*}
\tilde g \circ \varphi  = \mathcal{W}^{F,i}_{y_k,y_l} \circ \mathrm{pr} \circ c 
\end{equation*}
as functions from $U$ to $\mathrm{Mor}(G_k,G_l)$.
Because the right hand side is by assumption a smooth function,  $\tilde g$
is smooth on every chart, and hence also a smooth function. 
\end{proof}

\begin{theorem}
\label{th4}
Let $i:\mathrm{Gr} \to T$ be an equivalence of categories. A functor 
\begin{equation*}
F:\mathcal{P}_1(M) \to T
\end{equation*}
is a transport functor with $\mathrm{Gr}$-structure if and only if for every pair
$(x_1,x_2)$ of points in $M$ the
Wilson line $\mathcal{W}^{F,i}_{x_1,x_2}$
is smooth. 
\end{theorem}

\begin{proof}
One implication is shown by Lemma \ref{lem12}, using the fact that contractible
surjective submersions always exist. To prove the other implication we express the Wilson line of the transport functor  locally in terms of the functor $R_{(\mathrm{triv},g)}:
\upp{M}{\pi} \to T$ from Section \ref{sec5}. We postpone this construction to Appendix \ref{sec12}.
\end{proof}

Theorem \ref{th4} makes it possible to check explicitly, whether a given
functor $F$ is a transport functor or not. Furthermore, because every transport
functor has smooth Wilson lines, we can apply Lemma \ref{lem12} and have

\begin{corollary}
\label{co2}
Every transport functor $\mathrm{tra}:\mathcal{P}_1(M) \to T$ with $\mathrm{Gr}$-structure (with $i:\mathrm{Gr} \to T$ an equivalence of categories)
admits a $\pi$-local $i$-trivialization with smooth descent data for any contractible
surjective submersion $\pi$.
\end{corollary}

This corollary can be understood analogously to the fact, that every fibre
bundle over $M$ is trivializable over every good open cover of $M$. 
\begin{proposition}
\label{prop7}
For an equivalence of categories $i:\mathrm{Gr} \to T$ and a contractible surjective submersion $\pi:Y \to M$,
the forgetful functor 
\begin{equation*}
v^{\infty}: \loctrivfunctsmooth{i}{\mathrm{Gr}}{\pi}{1} \to \transport{i}{1}{\mathrm{Gr}}{T}
\end{equation*}
is a surjective equivalence of categories.
\end{proposition}

\begin{proof}
By Corollary \ref{co2} $v^{\infty}$ is  surjective. Since it is certainly faithful, it remains
to prove that it is full. Let $\eta$ be a morphism of transport functors
with $\pi$-local $i$-trivialization, i.e. there exists a surjective submersion
$\pi':Y' \to M$ such that $\ex{\pi'}(\eta)$ is smooth. Going to a contractible
surjective submersion $Z \to Y \times_M Y'$ shows that also $\ex{\pi}(\eta)$
is smooth. 
\end{proof}

Summarizing, we have for $i$ an equivalence of categories and $\pi$ a contractible
surjective submersion,  the following equivalences of categories: 
\begin{equation*}
\alxydim{@C=1.5cm}{\transsmooth{i}{1}{\mathrm{Gr}}  \ar@/^2pc/[r]^{\mathrm{\mathrm{Rec}_{\pi}}} & \loctrivfunctsmooth{i}{\mathrm{Gr}}{\pi}{1} \ar@/^2pc/[l]^{\ex{\pi}}  \ar[r]^-{v^{\infty}} & \transport{i}{1}{\mathrm{Gr}}{T}\text{.}  }
\end{equation*}

\section{Differential Forms and smooth Functors}

\label{sec7}

We establish a relation between smooth
descent data we have defined in the previous section  and more familiar geometric objects like differential forms, motivated by 
\cite{baez} and \cite{baez3}. The
relation we find can be expressed as a path ordered exponential, understood as the solution of an  initial value problem.

\begin{lemma}
\label{lem2}
Let $G$ be a Lie group with Lie algebra $\mathfrak{g}$. There is a canonical
 bijection between the set $\Omega^{1}(\R,\mathfrak{g})$ of $\mathfrak{g}$-valued 1-forms on $\R$ and the set of smooth
maps 
\begin{equation*}
f: \R \times \R \to G
\end{equation*}
satisfying the cocycle condition
\begin{equation}
\label{6}
f(y,z) \cdot f(x,y) = f(x,z)\text{.}
\end{equation}
\end{lemma}

\begin{proof}
The idea behind this bijection is that  $f$ is the path-ordered exponential of a 1-form $A$,
\begin{equation*}
f(x,y)= \mathcal{P}\exp \left (
  \int_{x}^y A \right )\text{.}
\end{equation*}
Let us explain in detail what that means. 
Given the 1-form $A$, we pose the initial value problem
\begin{equation}
\label{5}
\frac{\partial}{\partial t}u(t) =- \mathrm{d}r_{u(t)}|_1(A_t\left(\frac{\partial}{\partial
t}\right))
\quad\text{ and }\quad
u(t_0)=1
\end{equation}
for a smooth function $u: \R \to G$ and a number $t_0\in \R$. Here, $r_{u(t)}$ is the right multiplication
in $G$ and $\mathrm{d}r_{u(t)}|_1: \mathfrak{g} \to T_{u(t)}G$ is its differential
evaluated at $1\in G$. The sign in (\ref{5}) is a convention  well-adapted to the examples in Section \ref{sec4}. Differential equations of this type have a unique solution $u(t)$ defined on all of $\R$, such that $f(t_0,t) := u(t)$ depends smoothly
on both parameters. To see that $f$ satisfies the cocycle condition
(\ref{6}), define for fixed $x,y\in \R$ the function $\Psi(t):=f(y,t) \cdot f(x,y)$. Its derivative is
\begin{eqnarray*}
\frac{\partial}{\partial t}\Psi(t) &=&\mathrm{d}r_{f(x,y)}|_1 \left (
\frac{\partial}{\partial t}f(y,t)
\right )\nonumber
 \\&=&-  \mathrm{d}r_{f(x,y)}|_1(\mathrm{d}r_{f(y,t)}|_1(A_t\left(\frac{\partial}{\partial
 t}\right)))\\&=&-\mathrm{d}r_{\Psi(t)}(A_t \left (\frac{\partial}{\partial t}\right))\nonumber
\end{eqnarray*}
and furthermore $\Psi(y)=f(x,y)$. So, by uniqueness
\begin{equation*}
f(y,t) \cdot f(x,y) = \Psi(t) = f(x,t)\text{.}
\end{equation*}
Conversely, for a smooth function $f:\R \times \R \to G$, let 
$u(t):=f(t_0,t)$ for some $t_0\in \R$, and define
\begin{equation}
\label{27}
A_t\left(\frac{\partial}{\partial t}\right) := - \mathrm{d}r_{u(t)}|^{-1}_1\frac{\partial}{\partial t}u(t)\text{,}
\end{equation}
which yields a 1-form on $\R$.
If $f$ satisfies the cocycle condition, this 1-form is independent of the choice
of $t_0$. The definition of the 1-form $A$ is obviously inverse to (\ref{5}) and thus establishes the claimed bijection.
\end{proof}

We also need a relation between the  functions $f_A$ and $f_{A'}$ corresponding
to 1-forms $A$ and $A'$, when $A$
and $A'$ are related by a gauge transformation. In the following we denote  the left and
right invariant Maurer-Cartan forms on $G$ forms by $\theta$ and $\bar\theta$ respectively.

\begin{lemma}
\label{lem14}
Let $A\in\Omega^1(\R,\mathfrak{g})$ be a $\mathfrak{g}$-valued 1-form on
$\R$, let $g:\R \to G$ be a smooth function and let $A':=\mathrm{Ad}_{g}(A)
- g^{*}\bar\theta$. If $f_A$ and $f_{A'}$ are the smooth functions corresponding
to $A$ and $A'$ by Lemma \ref{lem2}, we have 
\begin{equation*}
 g(y) \cdot f_A(x,y)  =  f_{A'}(x,y) \cdot g(x)\text{.}
\end{equation*}
\end{lemma}

\begin{proof}
By direct verification, the function $g(y)\cdot f_{A}(x,y) \cdot g(x)^{-1}$
solves the initial value problem (\ref{5}) for the 1-form $A'$. Uniqueness
gives the claimed equality. 
\end{proof}

In the following we use the two lemmata above for 1-forms on $\R$ to obtain
a similar correspondence between 1-forms on an arbitrary smooth manifold $X$ and certain smooth
functors defined on the path groupoid $\mathcal{P}_1(X)$. For a given 1-form $A \in \Omega^1(X,\mathfrak{g})$, we first define a map
\begin{equation*}
k_A:
PX \to G
\end{equation*}
in the following way: a path $\gamma:x \to y$ in $X$ can be continued  to a smooth function $\gamma:\R \to X$ with $\gamma(t)=x$ for $t < 0$ and $\gamma(t)=y$
for $t>1$, due to its sitting instants. Then,
the pullback $\gamma^{*}A\in\Omega^1(\R,\mathfrak{g})$ corresponds by Lemma \ref{lem2} to a smooth function
$f_{\gamma^{*}A}:\R \times \R \to G$. Now we define 
\begin{equation*}
k_A(\gamma):=f_{\gamma^{*}A}(0,1)\text{.}
\end{equation*}
The map $k_A$ defined like this comes with the following properties:
\begin{itemize}
\item[a)]
For the constant path $\id_x$ we obtain the constant
function $f_{\id_x^{*}A}(x,y)=1$ and thus 
\begin{equation}
\label{21}
k_A(\id_x)=1\text{.}
\end{equation}
\item[b)]
For two paths
$\gamma_1:x \to y$ and $\gamma_2:y \to z$, we have
\begin{equation*}
f_{(\gamma_2 \circ\gamma_1)^{*}A}(0,1) = f_{(\gamma_2 \circ\gamma_1)^{*}A}(\frac{1}{2},1) \cdot f_{(\gamma_2 \circ\gamma_1)^{*}A}(0,\frac{1}{2})= f_{\gamma_1^{*}A}(0,1) \cdot f_{\gamma_2^{*}A}(0,1)
\end{equation*}
and thus
\begin{equation}
\label{22}
k_A(\gamma_2 \circ \gamma_1)=k_A(\gamma_2) \cdot k_A(\gamma_1)\text{.}
\end{equation}
\item[c)]
If $g:X \to G$ is a smooth function and $A':=\mathrm{Ad}_{g}(A) - g^{*}\bar\theta$,
\begin{equation}
\label{23}
g(y) \cdot k_A(\gamma)  = k_{A'}(\gamma) \cdot g(x)
\end{equation}
for any path $\gamma:x
\to y$.
\end{itemize}

The
next proposition shows that the definition of $k_A(\gamma)$ depends only on the thin
 homotopy  class of $\gamma$.

\begin{proposition}
\label{lem6}
The map $k_A: PX \to G$ factors in a unique way through the set $P^1X$ of thin homotopy classes
of paths, i.e. there is a unique map 
\begin{equation*}
F_A: P^1X \to G
\end{equation*}
such that $k_A=F_A \circ \mathrm{pr}$ with $\mathrm{pr}:PX \to P^1X$ the
projection.
\end{proposition}

\begin{proof}
If $k_A$ factors through the surjective map $\mathrm{pr}:PX \to P^1X$,
the map $F_A$ is determined uniquely. So we only have to show that two thin homotopy
equivalent paths $\gamma_0:x \to y$ and $\gamma_1:x \to y$ are mapped to
the same group element, $k_A(\gamma_0)=k_A(\gamma_1)$.
We have moved this issue to Appendix \ref{appD}.
\end{proof}

In fact, the map $F_A: P^1X \to G$ is not just a map. To understand it correctly, we need the following category:

\begin{definition}
\label{def17}
Let $G$ be a Lie group. We denote by $\mathcal{B}G$ the following Lie groupoid: it
has only one object, and $G$ is its set of morphisms. The unit element $1\in
G$ is the identity morphism, and  group multiplication is the composition,
i.e. $g_2 \circ g_1 := g_2 \cdot g_1$.
\end{definition}

To understand the notation, notice that the geometric realization of the nerve of $\mathcal{B}G$ yields the classifying space of the group $G$, i.e. $|N(\mathcal{B}G)|=BG$.
We claim that the map $F_A$ defined by Proposition \ref{lem6} defines a functor
\begin{equation*}
F_A: \mathcal{P}_1(X) \to \mathcal{B}G\text{.}
\end{equation*}
Indeed, since $\mathcal{B}G$ has only one object one only has to check that $F_A$ respects the composition (which is
shown by (\ref{21})) and  the identity morphisms (shown
in (\ref{22})).

\begin{lemma}
\label{lem7}
The functor $F_A$ is smooth in the sense of Definition \ref{def1}.
\end{lemma}

\begin{proof}
Let $U \subset \R^k$ be an open subset of some $\R^k$ and let $c: U \to PX$ be a map such that $c(u)(t)$ is smooth on $U \times [0,1]$.  We denote the path associated to a point $x\in U$ and extended smoothly to $\R$ by $\gamma_x := c(x):\R \to X$. This means that $U
\to \Omega^1(\R,\mathfrak{g}):x \mapsto \gamma_x^{*}A$
is a smooth family of $\mathfrak{g}$-valued 1-forms on $\R$. We recall that
\begin{equation*}
(k_A \circ c)(x)=k_A(\gamma_x)=f_{\gamma_x^{*}A}(0,1)
\end{equation*}
is defined to be the solution of a differential equation, which now  depends smoothly on $x$. Hence,  $  k_A \circ c = F_A \circ\mathrm{pr} \circ c: U
\to G$ is a smooth function. 
\end{proof}

Let us summarize the correspondence between 1-forms on $X$ and smooth functors
developed in the Lemmata above in  terms of  an equivalence between categories. One category is a category  \label{p3}$\funct^{\infty}(\mathcal{P}_1(X),\mathcal{B}
G)$ of smooth functors and  smooth natural transformations. The second category is the  category of differential $G$-cocycles on $X$: 

\begin{definition}
\label{def11}
Let $X$ be a smooth manifold and $G$ be a Lie group with Lie algebra $\mathfrak{g}$.
We consider the following category $\gconn{X}$: objects are all $\mathfrak{g}$-valued 1-forms $A$ on $X$, and a morphism $A \to A'$ is a smooth function $g:X \to G$ such that
\begin{equation*}
A' = \mathrm{Ad}_{g}(A) - g^{*}\bar\theta\text{.}
\end{equation*}
The composition is the multiplication of functions, $g_2 \circ g_1=g_2g_1$. 
\end{definition}
We claim that the Lemmata above provide the structure of a functor
\begin{equation*}
\fu: \gconn{X} \to \funct^{\infty}(\mathcal{P}_1(X),\mathcal{B}G)\text{.}
\end{equation*}
It sends a $\mathfrak{g}$-valued 1-form $A$ on $X$ to the functor $F_A$ defined
uniquely in Proposition \ref{lem6} and which is shown by Lemma \ref{lem7}.
It sends a function $g:X \to G$ regarded as a morphism $A \to A'$ to the smooth
natural
transformation $F_A \to F_{A'}$ whose component at a point $x$\ is $g(x)$. This is natural in $x$ due to (\ref{23}).

\begin{proposition}
\label{th2}
The functor 
\begin{equation*}
\fu: \gconn{X} \to \funct^{\infty}(X,\mathcal{B}G)\text{.}
\end{equation*}
is an isomorphism of categories, which reduces on the level of objects to a bijection
\begin{equation*}
\Omega^1(X,\mathfrak{g}) \cong \lbrace \text{Smooth functors }F:\mathcal{P}_1(X) \to \mathcal{B}G \rbrace\text{.}
\end{equation*}
\end{proposition}

\begin{proof} 
If $A$ and $A'$ are two $\mathfrak{g}$-valued 1-forms on $X$, the set of
morphisms between them is the set of smooth functions $g:X \to G$
satisfying the condition $A'=\mathrm{Ad}_g(A)-g^{*}\bar\theta$. The set of morphisms between
the functors
$F_A$ and $F_{A'}$ are smooth natural transformations, i.e. smooth maps $g:X
\to G$, whose naturality square is equivalent to the same condition. So,
the functor $\fu$ is manifestly full and faithful. It remains to show that
it is  a bijection on the level of objects.  This is done in Appendix \ref{appE} by an explicit construction
of a 1-form $A$ to a given smooth functor $F$.
\end{proof}

One can also enhance the category $\gconn{X}$ in such a way that it becomes
the familiar category of local data of principal $G$-bundles with connection.

\begin{definition}
\label{def12}
  The category $\diffco{G}{1}{\pi}$ of 
  \emph{differential $G$-cocycles} of the surjective submersion $\pi$ is
  the category whose objects are pairs $(g,A)$  consisting of a
  1-form $A \in \Omega^1(Y,\mathfrak{g})$
  and a smooth function $g:Y^{[2]} \to G$ such that
\begin{equation*}
\pi_{13}^* g = \pi_{23}^* g \cdot \pi_{12}^* g
\quad\text{ and }\quad
\pi_2^* A = \mathrm{Ad}_g(\pi_1^*A) - g^* \bar\theta\text{.}
\end{equation*}
A morphism
\begin{equation*}
h: (g,A) \to (g',A')
\end{equation*}
is a smooth function $h :Y \to G$
such that
\begin{equation*}
A' = \mathrm{Ad}_h(A) - h^* \bar\theta
\quad\text{ and }\quad
 \pi_2^* h \cdot g =g' \cdot \pi_1^*h\text{.}
\end{equation*}
Composition of morphisms is given by the product of these functions, $h_2 \circ h_1=h_2h_1$.
\end{definition}

To explain the notation, notice that for $\pi=\id_X$ we obtain $\gconn{X}=\diffco{G}{1}{\pi}$.
As an example, we consider the group $G=U(1)$ and a surjective submersion $\pi:Y
\to M$ coming from a good open cover $\mathfrak{U}$ of $M$. Then, the group of isomorphism classes of $\diffco{U(1)}{1}{\pi}$ is the Deligne hypercohomology
group $H^1(\mathfrak{U},\mathcal{D}(1))$, where $\mathcal{D}(1)$ is the Deligne
sheaf complex $0 \to \sheaf{U(1)} \to \sheaf{\Omega}^1$.

\begin{corollary}
\label{co1}
The functor $\fu$ extends to an equivalence of categories
\begin{equation*}
\diffco{G}{1}{\pi} \cong \transsmooth{i_G}{1}{\mathrm{Gr}} \text{,}
\end{equation*}
where $i_{G}: \mathcal{B}G \to G\text{-}\mathrm{Tor}$  sends the object of $\mathcal{B}
G$ to
the group $G$ regarded as a $G$-space, and a morphism $g\in G$ to the equivariant
smooth map  which   multiplies with $g$ from the left.
 
\end{corollary}
This corollary is an important step towards our main theorem, to which we
come in the next section.

\section{Examples}

\label{sec4}

Various structures in the theory of bundles with connection are special
cases of transport functors
with $\mathrm{Gr}$-structure for particular choices of the structure groupoid $\mathrm{Gr}$. In
this section we  spell out some prominent examples. 

\subsection{Principal Bundles with Connection}

\label{sec5_1}

In this section, we fix a Lie group
$G$. Associated to this Lie group,
we have the Lie groupoid $\mathcal{B}G$ from Definition \ref{def17},
 the category
 $G\text{-}\mathrm{Tor}$  of smooth manifolds with right $G$-action
and $G$-equivariant smooth maps between those, and the functor $i_G: \mathcal{B}G \to G\text{-}\mathrm{Tor}$ that sends the object of $\mathcal{B}G$ to the $G$-space $G$ and a morphism $g\in G$ to the $G$-equivariant
diffeomorphism that multiplies from the left by $g$. The functor $i_G$
is an equivalence of categories.

As we have outlined in the introduction, a principal $G$-bundle $P$ with
connection over $M$ defines a functor 
\begin{equation*}
\mathrm{tra}_{P}: \mathcal{P}_1(M) \to G\text{-}\mathrm{Tor}\text{.}
\end{equation*}
Before we show that $\mathrm{tra}_P$
is a transport functor with 
$\mathcal{B}G$-structure, let us recall
its definition in detail. To an object $x\in M$ it assigns the fibre $P_{x}$ of the bundle $P$ over
the point $x$. To a path $\gamma:x \to y$, it assigns the parallel transport
map $\tau_{\gamma}: P_x \to P_y$.

For preparation, we recall the
basic definitions concerning local
trivializations of principal bundles with connections.
In the spirit of  this article, we
use surjective submersions instead
of coverings by open sets. In
this language, a local trivialization
of the principal bundle $P$ is
a surjective submersion $\pi:Y
\to M$ together with a $G$-equivariant
diffeomorphism
\begin{equation*}
\phi: \pi^{*}P \to Y \times G
\end{equation*}
that covers the identity on $Y$. Here,
the fibre product
$\pi^{*}P=Y \times_M P$ comes with
the projection $p: \pi^{*}P \to
P$ on the second factor. It induces
a section 
\begin{equation*}
s:Y \to P: y \mapsto p(\phi^{-1}(y,1))\text{.}
\end{equation*}
The transition function $\tilde g_{\phi}:Y^{[2]} \to G$ associated to the local trivialization $\phi$ is defined by 
\begin{equation}
\label{29}
s(\pi_1(\alpha)) = s(\pi_2(\alpha)) \cdot \tilde g_{\phi}(\alpha)
\end{equation}
for every point $\alpha\in Y^{[2]}$. A connection on $P$
is a $\mathfrak{g}$-valued 1-form
$\omega\in \Omega^1(P,\mathfrak{g})$ that obeys
\begin{equation}
\label{3}
\omega_{\rho g}\left ( \frac{\mathrm{d}}{\mathrm{d}t}(\rho
g) \right ) = \mathrm{Ad}^{-1}_{g}
\left (\omega_\rho \left( \frac{\mathrm{d}\rho}{\mathrm{d}t} \right) \right) + \theta_{g}\left
( \frac{\mathrm{d}g}{\mathrm{d}t}
\right )
\end{equation}
for smooth maps $\rho:[0,1]\to
P$ and $g:[0,1] \to G$. In this
setup, a tangent vector $v \in
T_{p}P$ is called horizontal, if it
is in the kernel of $\omega$. 

Notice that all our conventions are chosen such that the transition function $\tilde g_{\phi}:Y^{[2]} \to G$ and the local connection 1-form $\tilde A_\phi:=s^{*}\omega\in \Omega^1(Y,\mathfrak{g})$ define an object in the category $\diffco{G}{1}{\pi}$ from Definition \ref{def12}. 

To define the parallel transport
map $\tau_{\gamma}$ associated
to a path $\gamma: x \to y$ in
$M$, we assume
first that $\gamma$ has a lift $\tilde \gamma:
\tilde x \to \tilde y$ in $Y$, that is,
$\pi_{*} \tilde\gamma= \gamma$.
Consider then the path $s_{*}\tilde\gamma$
in $P$, which can be modified by
the pointwise action of a path $g$ in $G$
from the right, $(s_{*}\tilde\gamma
)g$. This modification has  now to be chosen
such that every tangent vector
to $(s_{*}\tilde \gamma)g$  is
horizontal, i.e.
\begin{equation*}
0=\omega_{(s_{*}\tilde\gamma)g}
\left (\frac{\mathrm{d}}{\mathrm{d}t}((s_{*}\tilde\gamma)g)
\right)
\stackrel{(\ref{3})}{=} \mathrm{Ad}^{-1}_{g}\left(\omega_{s_{*}\tilde\gamma}\left(\frac{\mathrm{d}(s_{*}\tilde\gamma)}{\mathrm{d}t}
\right) \right)
+ \theta_{g} \left (\frac{\mathrm{d}g}{\mathrm{d}t}\right)
\end{equation*}
This is a
linear differential equation for
$g$, which
has together with the initial condition
$g(0)=1$ a unique
solution $g=g(\tilde\gamma)$. Then,
 for any $p\in P_x$, 
\begin{equation}
\label{28}
\tau_{\gamma}(p) := s(y)(g(1)\cdot
h)\text{,}
\end{equation}
where $h$ is the unique  group element
with $s(x)h=p$.  It is evidently smooth in
$p$ and $G$-equivariant.
Paths $\gamma$ in $M$ which do not have a lift to $Y$ have to be split up in pieces
which admit lifts; $\tau_{\gamma}$ is then the composition of the 
parallel transport maps of those.

\begin{lemma}
\label{lem10}
Let $P$ be a principal $G$-bundle
over $M$ with connection $\omega\in\Omega^1(P,\mathfrak{g})$.
For a surjective submersion
$\pi:Y \to M$  and a trivialization $\phi$ with associated section $s:Y \to P$, we consider the 
 smooth functor 
\begin{equation*}
F_{\omega}
:= \fu(s^{*}\omega): \mathcal{P}_1(Y)
\to \mathcal{B}G
\end{equation*}
associated to the 1-form $s^{*}\omega\in\Omega^1(Y,\mathfrak{g})$
by Proposition \ref{th2}. Then,
\begin{equation}
\label{36}
i_G(F_{\omega}(\overline{\gamma})) = \phi_y \circ \tau_{\pi_{*}\gamma} \circ
\phi_x^{-1}
\end{equation}
for any path $\gamma:x \to y$ in
$PY$.  
\end{lemma}

\begin{proof}
Recall the definition of the functor
$F_{\omega}$: for a path $\gamma:
x \to  y$, we have to consider
the 1-form $\gamma^{*}s^{*}\omega\in\Omega^1(\R,\mathfrak{g})$,
which defines a smooth function
$f_{\omega}:\R \times \R \to G$.
Then, $F_{\omega}(\overline{\gamma}):=f_{\omega}(0,1)$.
We claim the equation
\begin{equation}
\label{18}
f_{\omega}(0,t)=g(t)\text{.}
\end{equation}
This comes from the fact that both
functions are solutions of the
same differential equation, with
the same initial value for $t=0$.
Using (\ref{18}),
\begin{equation*}
i_{G}(F_{\omega}(\overline{\gamma}))(h)=F_{\omega}(\overline{\gamma})\cdot
h = g(1) \cdot h
\end{equation*}
for some $h\in G$. On the other hand,
\begin{equation*}
\phi_{y} (\tau_{\pi_{*}\gamma}
(\phi_x^{-1}(h))) =\phi_{y} (\tau_{\pi_{*}\gamma}
(s(x)h)) \stackrel{(\ref{28})}{=} \phi_y(s(y)(g(1)\cdot
h))=g(1)\cdot h\text{.}
\end{equation*}
This proves   equation (\ref{36}). 
\end{proof}

Now we are ready to formulate the basic relation between principal $G$-bundles
with connection and transport functors with $\mathcal{B}G$-structure.

\begin{proposition}
\label{prop1}
The functor 
\begin{equation*}
\mathrm{tra}_P: \mathcal{P}_1(M)
\to G\text{-}\mathrm{Tor}
\end{equation*}
obtained from parallel transport in a principal
$G$-bundle $P$, is a transport functor with $\mathcal{B}G$-structure 
in the sense of Definition \ref{def3}.
\end{proposition}

\begin{proof}
The essential ingredient is, that $P$ is locally trivializable: we choose
a surjective submersion $\pi:Y \to M$ and a trivialization $\phi$. 
The construction of a  functor $\mathrm{triv}_{\phi}: \mathcal{P}_1(Y)
\to \mathcal{B}G$ and a natural equivalence
\begin{equation*}
\alxydim{@=1.2cm}{\mathcal{P}_1(Y) \ar[r]^{\pi_{*}} \ar[d]_{\mathrm{triv}_{\phi}} & \mathcal{P}_1(M)
\ar@{=>}[dl]|{t_{\phi}}
\ar[d]^{\mathrm{tra}_P} \\ \mathcal{B}G \ar[r]_-{i_G} & G\text{-}\mathrm{Tor}}
\end{equation*}
is as follows. We let $\mathrm{triv}_{\phi}
:= \fu(s^{*}\omega)$ be the smooth
functor associated to the 1-form
$s^{*}\omega$ by Proposition \ref{th2}. To define the natural equivalence $t_{\phi}$, consider a point
$x \in Y$. We find $\pi^{*}\mathrm{tra}_P(x)=P_{\pi(x)}$ and $(i_G \circ \mathrm{triv}_{\phi})(x)=G$.
So we define the component of $t_{\phi}$ at $x$ by
\begin{equation*}
t_{\phi}(x):=\phi_x: P_{\pi(x)}
\to G\text{.}
\end{equation*}
This is natural in $x$ since the diagram
\begin{equation*}
\alxydim{@=1.25cm}{P_{\pi(x)} \ar[d]_{\tau_{\pi_{*}\overline{\gamma}}} \ar[r]^-{\phi_x}
& G \ar[d]^{i_G(\mathrm{triv}_{\phi}(\overline{\gamma}))} \\ P_{\pi(y)} \ar[r]_-{\phi_y} & G}
\end{equation*}
is commutative by Lemma \ref{lem10}. Notice that the natural
equivalence 
\begin{equation}
\label{38}
g_{\phi} := \pi_2^{*}t_{\phi} \circ \pi_1^{*}t_{\phi}
\end{equation} 
factors through the smooth transition
function $\tilde g_{\phi}$ from
(\ref{29}), i.e.  $g_{\phi}=i_G(\tilde
g_{\phi})$. Hence, the pair $(\mathrm{triv}_{\phi},g_{\phi})$ is a smooth object in $\transsmooth{i}{1}{i_G}$.
\end{proof}

Now we consider the morphisms. Let $\varphi: P \to P'$ be a morphism of principal $G$-bundles over $M$
(covering the identity on $M$)
which respects the connections,
i.e. $\omega = \varphi^{*}\omega'$. For any point $p\in M$, its restriction $\varphi_x:P_x
\to P_x'$ is a smooth $G$-equivariant map.
For any path $\gamma:x \to y$, the parallel transport map satisfies
\begin{equation*}
\varphi_y \circ \tau_{\gamma} = \tau'_{\gamma} \circ \varphi_x\text{.}
\end{equation*}
This is
nothing but the commutative diagram for the components $\eta_{\varphi}(x):=\varphi_x$ natural transformation $\eta_{\varphi}:\mathrm{tra}_P
\to \mathrm{tra}_{P'}$.

\begin{proposition}
\label{prop2}
The natural transformation 
\begin{equation*}
\eta_{\varphi}:
\mathrm{tra}_P \to \mathrm{tra}_{P'}
\end{equation*}
obtained from a morphism $\varphi:P
\to P'$ of principal
$G$-bundles, is a morphism of transport functors in the sense of Definition \ref{def3}.
\end{proposition}

\begin{proof}
Consider a surjective submersion $\pi:Y \to M$ such that $\pi^{*}P$ and $\pi^{*}P'$
are trivializable, and choose  trivializations $\phi$ and $\phi'$. The descent
datum of $\eta_{\varphi}$ is the natural equivalence $h:=t_{\phi}' \circ \pi^{*}\eta_{\varphi} \circ t_{\phi}^{-1}$. Now define the 
 map
\begin{equation*}
\tilde h: Y \to G: x \mapsto p_G(\phi'(x,\varphi(s(x))))
\end{equation*}
where $p_G$ is the projection to $G$. This map  is smooth and satisfies $h=i_G(\tilde h)$. Thus, $\eta_{\varphi}$ is a morphism of transport functors.
\end{proof}

Taking the  Propositions \ref{prop1}
and \ref{prop2}  together, we have defined a functor
\begin{equation}
\label{48}
\gbun \to \transport{i}{1}{\mathrm{Gr}}{G\text{-}\mathrm{Tor}}
\end{equation}
from the category of principal $G$-bundles over $M$ with connection to the
category of transport functors on $M$ in $G\text{-}\mathrm{Tor}$ with $\mathcal{B}G$-structure. In particular, this functor provides us with
lots of examples of transport functors.  
\begin{theorem}
\label{th1}
The functor
\begin{equation}
\label{19}
\gbun \to \transport{i}{1}{\mathcal{B}G}{G\text{-}\mathrm{Tor}} 
\end{equation}
is an  equivalence of categories.
\end{theorem}

We give two proofs of this Theorem: the first is short and the second is
explicit.
\medskip

\begin{proofblank}{First Proof} Let $\pi:Y \to M$ be a contractible
surjective submersion, over which every principal $G$-bundle
is trivializable. Extracting a connection 1-form $\tilde A_{\phi}\in\Omega^1(Y,\mathfrak{g})$
and the transition function (\ref{29}) yields a functor 
\begin{equation*}
\gbun \to \diffco{G}{1}{\pi}
\end{equation*}
to the category of differential $G$-cocycles for $\pi$, which is in fact an
equivalence of categories. We claim that the composition of
this equivalence with the sequence
\begin{equation}
\label{56}
\alxydim{@=1.01cm}{\diffco{G}{1}{\pi} \ar[r]^-{\fu}
& \transsmooth{i}{1}{i_{G}} \ar[r]^-{\mathrm{Rec}_{\pi}} & \loctrivfunctsmooth{i}{\mathrm{Gr}}{\pi}{1} \ar[r]^-{v^{\infty}} & \transport{i}{1}{\mathrm{Gr}}{G\text{-}\mathrm{Tor}}}
\end{equation}
of functors is naturally equivalent to the functor (\ref{19}). By Corollary \ref{co1}, Theorem \ref{th3} and Proposition \ref{prop7} all functors in
(\ref{56}) are equivalences of categories, and so is (\ref{19}). To show
the claim recall that in the proof of Proposition \ref{prop1} we have defined a local trivialization of $\mathrm{tra}_P$, whose descent data $(\mathrm{triv}_{\phi},g_{\phi})$
is the image of the local data $(\tilde A_{\phi},\tilde g_{\phi})$ of the
principal $G$-bundle under
the functor $\fu$. This reproduces exactly the steps in the sequence
(\ref{56}).  
\end{proofblank}

\begin{proofblank}{\secondproof} We show that the functor (\ref{19}) is faithful, full and essentially
surjective. In fact, this proof shows that it is even  surjective.  So let $P$ and $P'$ two principal $G$-bundles with connection
over $M$, and let $\mathrm{tra}_P$ and $\mathrm{tra}_{P'}$ be the associated
transport functors. 

Faithfulness follows directly from the definition, so
assume now that  $\eta: \mathrm{tra}_P \to \mathrm{tra}_{P'}$ is a morphism of  transport
functors. We define a morphism $\varphi:P \to P$ pointwise as $\varphi(x):=\eta(p(x))(x)$
for any $x\in P$, where $p:P \to M$ is the projection of the bundle $P$.
This is clearly a preimage of $\eta$ under the functor (\ref{19}), so that we only
have to show that $\varphi$ is a smooth map. We choose a surjective
submersion such that   $P$ and $P'$ are trivializable and such that   $h:=\ex{\pi}(\eta)=t_{\phi'} \circ \pi^{*}\eta
\circ t_{\phi}^{-1}$
is a smooth morphism in $\transsmooth{i}{1}{i_G}$. Hence it factors through
a smooth map $\tilde h: Y \to G$, and from the definitions of $t_{\phi}$
and $t_{\phi'}$ it follows that $\pi^{*}\varphi$ is the function
\begin{equation*}
\pi^{*}\varphi:\pi^{*}P \to \pi^{*}P':(y,p) \mapsto
\phi'^{-1}(\phi(y,p)\tilde h(y))\text{,}
\end{equation*}
and thus  smooth. Finally, since $\pi$ is a surjective submersion,
$\varphi$ is smooth. 

It  remains to prove that the functor (\ref{19}) is essentially surjective.
First we construct, for a given transport functor $\mathrm{tra}:\mathcal{P}_1(M) \to G\text{-}\mathrm{Tor}$
a principal $G$-bundle $P$ with connection over $M$, performing exactly the
inverse steps of (\ref{56}). We choose a surjective submersion $\pi:Y \to
M$ and a $\pi$-local $i$-trivialization $(\mathrm{triv},t)$ of the transport functor $\mathrm{tra}$. By construction, its descent data $(\mathrm{triv},g):=\ex{\pi}(\mathrm{triv},g)$
is an object in $\transsmooth{i}{1}{i_G}$. 
By Corollary \ref{co1}, there exists a 1-form $A\in\Omega^1(Y,\mathfrak{g})$, and a smooth function $\tilde g: Y^{[2]} \to G$, forming an object $(A,\tilde
g)$
in the category $\diffco{G}{1}{\pi}$ of differential cocycles such that
\begin{equation}
\label{34}
\fu(A,\tilde g) = (\mathrm{triv},g)
\end{equation}
in $\transsmooth{i}{1}{i_G}$. In particular $g = i_G(\tilde g)$.
The pair $(A,\tilde g)$ is local data for a principal $G$-bundle
$P$ with connection $\omega$. The reconstructed bundle comes with a canonical trivialization $\phi: \pi^{*}P \to Y
\times G$, for which the associated section $s:Y \to P$ is such that 
$A=s^{*}\omega$, and whose transition function is $\tilde g_{\phi}=\tilde g$.

Let us extract descent data of the transport functor $\mathrm{tra}_P$
of $P$:  as described in the proof of Proposition \ref{th1}, the
trivialization $\phi$ of the bundle $P$ gives rise to a $\pi$-local $i_G$-trivialization
$(\mathrm{triv}_{\phi}, t_{\phi})$  of the transport functor $\mathrm{tra}_P$,
namely
\begin{equation}
\label{40}
\mathrm{triv_{\phi}}:=F_{\omega}:=\fu(s^{*}\omega)=\fu(A)
\end{equation}
and $t_{\phi}(x) := \phi_x$. Its natural equivalence $g_{\phi}$ from (\ref{38})
is just $g_{\phi}=i_G(\tilde g_{\phi})$. 

Finally we construct an isomorphism $\eta: \mathrm{tra}_P \to \mathrm{tra}$
of transport functors.  Consider
the natural equivalence
\begin{equation*}
\zeta := t^{-1} \circ  t_{\phi}: \pi^{*}\mathrm{tra}_{P} \to
\pi^{*}\mathrm{tra}\text{.}
\end{equation*}
From condition (\ref{201}) it follows that $\zeta(\pi_1(\alpha))=\zeta(\pi_2(\alpha))$
for every point $\alpha\in Y^{[2]}$. So $\zeta$ descends to a natural equivalence
\begin{equation*}
\eta(x):=\zeta(\tilde x)
\end{equation*}
for $x\in M$ and any $\tilde x\in Y$ with $\pi(\tilde x)=x$. An easy computation
shows that
$\ex{\pi}(\eta)=t \circ \zeta \circ  t_{\phi}^{-1}=\id$,
which is in particular smooth and thus proves that $\eta$ is an isomorphism in $\transsmooth{i}{1}{i_G}$.
\end{proofblank}

\subsection{Holonomy Maps}

\label{sec9}

In this section, we show that important results of \cite{barret1,caetano} on holonomy
maps of principal $G$-bundles with connection can
be reproduced as particular cases. 

\begin{definition}[\cite{caetano}]
\label{def20}
A holonomy  map on a smooth manifold $M$ at a point $x\in M$ is a group homomorphism \begin{equation*}
\mathcal{H}_x:\pi_1^1(M,x)
\to G\text{,}
\end{equation*}
which is smooth in the following sense: for every open subset $U \subset
\R^{k}$ and every map $c: U \to L_xM$ such that $\Gamma(u,t):=c(u)(t)$ is smooth
on $U \times [0,1]$, also
\begin{equation*}
\alxy{U \ar[r]^-{c} & L_xM \ar[r]^-{\mathrm{pr}} & \pi_1^1(M,x) \ar[r]^-{\mathcal{H}} & G}
\end{equation*}
is smooth.
\end{definition}

Here, $L_xM \subset PM$ is the set of paths $\gamma:x \to x$, whose image
under the projection $\mathrm{pr}:PM \to P^1M$ is, by definition, the thin
homotopy group $\pi_1^1(M,x)$ of $M$ at $x$. Also notice, that
\begin{itemize}
\item 
in the context of diffeological spaces  reviewed in Appendix \ref{appB},
the definition of smoothness given here just means that $\mathcal{H}$ is
a morphism between diffeological spaces,
cf. Proposition \ref{prop5} ii). 
\item
the notion of \textit{intimate paths} from \cite{caetano} and the notion
of thin homotopy from \cite{mackaay1} coincides with
our notion of thin homotopy, while the notion of thin homotopy used in \cite{barret1}
is different from ours.  
\end{itemize}
In \cite{caetano} it has been shown that parallel transport in a principal $G$-bundle over $M$
around based loops defines
a holonomy map.  For connected
manifolds $M$ it was also shown how to
reconstruct a principal $G$-bundle with connection from a given holonomy
map $\mathcal{H}$ at $x$, such that the holonomy of this bundle around loops based at $x$
equals $\mathcal{H}$. This establishes a bijection between holonomy maps
and principal $G$-bundles with connection over connected manifolds. The same result has been proven
(with the before mentioned different notion of thin homotopy)
in \cite{barret1}. 

\medskip

To relate these results to Theorem \ref{th1}, we consider again transport
functors $\mathrm{tra}:\mathcal{P}_1(M) \to G\text{-}\mathrm{Tor}$ with $\mathcal{B}G$-structure. Recall from Section \ref{sec10} that for   any point
$x \in M$ and any identification $F(x) \cong G$ the functor $\mathrm{tra}$ produces a group homomorphism $F_{x,x}: \pi_1^1(M,x) \to G$.

\begin{proposition}
\label{prop6}
Let $\mathrm{tra}: \mathcal{P}_1(M) \to G\text{-}\mathrm{Tor}$ be a transport functor on $M$ with
$\mathcal{B}G$-structure. Then, for any point $x \in M$  and any identification
 $F(x) \cong G$, the group
homomorphism
\begin{equation*}
\mathrm{tra}_{x,x}:\pi_1^1(M,x) \to G
\end{equation*}
is a holonomy map.
\end{proposition}

\begin{proof}
The group homomorphism $\mathrm{tra}_{x,x}$ is a Wilson line of the transport
functor $\mathrm{tra}$, and hence smooth by Theorem \ref{th4}. 
\end{proof}

For illustration, let us combine Theorem \ref{th1} and Proposition \ref{prop6} to the following
diagram, which is evidently commutative:
\begin{equation*}
\alxydim{@C=1cm@R=2.1cm}{\gbun \ar[dr]_{\txt{\small
Theorem
\ref{th1} \;\;}} \ar[rr]^{\txt{\small\cite{caetano}}} && \left \lbrace \txt{Holonomy
maps\\ on $M$ at $x$} \right \rbrace   \\ &\transport{i}{1}{\mathrm{Gr}}{G\text{-}\mathrm{Tor}}\text{.} \ar[ur]_{\txt{\;\;\;\;\;\small
Proposition
\ref{prop6}}}&} 
\end{equation*}

\subsection{Associated Bundles and Vector Bundles with Connection}

Recall that a principal $G$-bundle $P$ together with a faithful representation
$\rho:G \to \mathrm{Gl}(V)$ of
the Lie group $G$ on a vector space $V$ defines a vector bundle
$P \times_{\rho} V$
with structure group $G$, called the vector bundle associated to $P$ by the
representation $\rho$. One can regard a (say, complex) representation of a group $G$ conveniently
as a functor $\rho: \mathcal{B}G \to \mathrm{Vect}(\C)$ from the one-point-category
$\mathcal{B}G$ into the category of complex vector spaces: the object of $\mathcal{B}
G$ is sent to the vector space $V$ of the representation, and a group element
$g \in G$  is sent to an isomorphism $g:V \to V$ of this vector space. The
axioms of a functor are precisely the axioms one demands for a representation.
Furthermore, the representation is faithful, if and only if the functor is faithful. 

\begin{definition}
Let
\begin{equation*}
    \rho : \mathcal{B}G \to \mathrm{Vect}(\C)
\end{equation*}
be any representation of the Lie group $G$.
A transport functor 
\begin{equation*}
    \mathrm{tra} : \mathcal{P}_1(M) \to \mathrm{Vect}(\C)
\end{equation*}
with $\mathcal{B}G$-structure is called \emph{associated transport functor}.
\end{definition}

As an example, we consider the
defining representation of the
Lie group $U(n)$ on the vector
space $\C^n$, considered as a functor
\begin{equation}
\label{57}
\rho_n: \mathcal{B} U(n) \to \mathrm{Vect}(\C_h^n)
\end{equation}
to the category of $n$-dimensional
hermitian vector spaces and isometries between those. Because we only include
isometries in $\mathrm{Vect}(\C^n_h)$, the functor $\rho_n$ is an equivalence
of categories. 

 Similarly to Theorem
\ref{th1}, we find a geometric interpretation for associated transport functors on $M$ with
 $\mathcal{B} U(n)$-structure, namely hermitian vector bundles of rank $n$
with (unitary) connection over $M$. We denote the category of those vector
bundles
by \label{p5} $\mathrm{VB}(\C^n_h)^{\nabla}_M$. 
 Let us just outline
the very basics: given such a vector bundle $E$, we   associate a
functor
\begin{equation*}
\mathrm{tra}_E: \mathcal{P}_1(M)
\to \mathrm{Vect}(\C_h^n)\text{,}
\end{equation*}
which sends a point $x\in M$ to the vector
space $E_x$, the fibre of $E$ over
$x$, and   a path $\gamma:x
\to y$ to the parallel transport
map $\tau:E_x \to E_y$, which is
 linear and an isometry. 

\begin{theorem}
\label{th7}
The functor $\mathrm{tra}_E$ obtained from a hermitian vector bundle $E$
with connection over $M$ is a transport functor on $M$ with $\mathcal{B} U(n)$-structure; furthermore, the assignment $E \mapsto \mathrm{tra}_E$ yields a
functor
\begin{equation}
\label{58}
\mathrm{VB}(\C^n_h)^{\nabla}_M \to \transport{i}{1}{\mathcal{B} U(n)}{\mathrm{Vect}(\C_h^n)}\text{,}
\end{equation}
which is an equivalence of categories.
\end{theorem}

\begin{proof}
We proceed like in the first proof of Theorem \ref{th1}. Here we use 
the correspondence between hermitian vector bundles with connection and their
local data in $\diffco{U(n)}{1}{\pi}$, for contractible surjective submersions
$\pi$. Under this correspondence  the functor (\ref{58})  becomes naturally equivalent
to the composite
\begin{equation}
\alxydim{@=1.0cm}{\diffco{U(n)}{1}{\pi} \ar[r]^-{\Xi}
& \transsmooth{i}{1}{\rho_n} \ar[r]^-{\mathrm{Rec}_{\pi}} & \loctrivfunctsmooth{i}{\mathrm{Gr}}{\pi}{1} \ar[r]^-{v^{\infty}} & \transport{i}{1}{\mathcal{B} U(n)}{\mathrm{Vect}(\C_h^n)}}
\nonumber
\end{equation}
which is, by Corollary \ref{co1}, Theorem \ref{th3} and Proposition \ref{prop7},
an equivalence of categories.  
\end{proof}

Let us also consider the Lie groupoid $\mathrm{Gr}_{U} := \bigsqcup_{n\in \N} \mathcal{B} U(n)$, whose set of objects is $\N$ (with the discrete smooth structure) and whose morphisms are
\begin{equation*}
\mathrm{Mor}_{\mathrm{Gr}_U}(n,m) = \begin{cases} U(n) & \text{if } n=m \\
\emptyset & \text{if }n \neq m \\
\end{cases}
\end{equation*}
so that $\mathrm{Mor}(\mathrm{Gr}_U)$ is a disjoint union of Lie groups.
The functors $\rho_n$ from (\ref{57}) induce a functor
\begin{equation*}
\rho_{U}: \mathrm{Gr}_U \to \mathrm{Vect}(\C_h)
\end{equation*} 
to the category of hermitian vector spaces (without a fixed dimension) and
isometries between those. 

The category $\mathrm{Vect}(\C_h)$ in fact a monoidal category, and its monoidal
structure induces monoidal structures on the category $\mathrm{VB}(\C_h)^{\nabla}_M$
of hermitian vector bundles with  connection over $M$ as well as on
the category of transport functors $\transport{i}{1}{\mathrm{Gr}_U}{\mathrm{Vect}(\C_h)}$,  as outlined
in Section \ref{sec10}. Since parallel transport in vector bundles is compatible
with tensor products, we have
\begin{corollary}
The functor 
\begin{equation*}
 \mathrm{VB}(\C_h)^{\nabla}_M
\to \transport{i}{1}{\mathrm{Gr}_U}{\mathrm{Vect}(\C_h)}
\end{equation*}
is a monoidal equivalence of monoidal categories. 
\end{corollary} 
 
In particular, we have the unit transport functor $\trivlin_{\C}$  which sends every point to the complex numbers $\C$, and every path to the identity $\id_{\C}$. The following fact is easy to verify:

\begin{lemma}
\label{lem9}
Let $\mathrm{tra}: \mathcal{P}_1(M) \to \mathrm{Vect}(\C_h)$ be a transport functor with $\mathrm{Gr}_U$-structure, corresponding to a hermitian vector bundle $E$ with connection over $M$. Then, there is a canonical bijection between morphisms
\begin{equation*}
\eta: \trivlin_{\C} \to \mathrm{tra}
\end{equation*}
of transport functors with $\mathrm{Gr}_{U}$-structure and smooth flat section of $E$.
\end{lemma}

\subsection{Generalized Connections}

In this section we  consider 
 functors 
\begin{equation*}
F : \mathcal{P}_1(M) \to \mathcal{B}G\text{.}
\end{equation*}
By now, we can arrange such functors in three types: 
\begin{enumerate}
\item 
We demand nothing of $F$:  such functors are  addressed as \textit{generalized connections}
\cite{ashtekar1}. 

\item
We demand that $F$ is a transport functor with $\mathcal{B}G$-structure: it corresponds to an ordinary principal $G$-bundle with connection. 

\item
We demand that $F$ is smooth in the sense of  Definition \ref{def1}:  by Proposition \ref{th2},
one can replace such functors by 1-forms $A \in \Omega^1(M,\mathfrak{g})$,
so that we can speak of a trivial $G$-bundle. 
\end{enumerate}
Note that for a functor $F:\mathcal{P}_1(M) \to \mathcal{B}G$ and the identity
functor $\id_{\mathcal{B}G}$ on $\mathcal{B}G$ the Wilson line 
\begin{equation*}
\mathcal{W}^{F,\id_{\mathcal{B}G}}_{x_1,x_2}: \mathrm{Mor}_{\mathcal{P}_1(M)}(x_1,x_2)
\to G
\end{equation*}
does not depend on choices of objects $G_1$, $G_2$ and morphisms $t_k: i(G_k)
\to F(x_k)$ as in the general setup described in Section \ref{sec11}, since
$\mathcal{B}G$ has only one object and one can canonically choose $t_k=\id$.
So, generalized connections have a particularly good Wilson lines. Theorem \ref{th4} provides a precise criterion to decide
when a generalized
connection is regular: if and only if all its Wilson lines are smooth.

\section{Groupoid Bundles with
Connection}

\label{sec15}

In all examples we have discussed so far the Lie groupoid $\mathrm{Gr}$ is of the
form $\mathcal{B}G$, or a union of those. In this section we discuss  transport
functors with $\mathrm{Gr}$-structure for a general Lie groupoid $\mathrm{Gr}$. We start with the local aspects of such transport functors, and then discuss two examples of target categories. Our main example is related to the notion of principal groupoid bundles \cite{moerdijk}. In contrast to the examples in Section \ref{sec4}, transport functors with $\mathrm{Gr}$-structure do not only reproduce the existing definition of a principal groupoid bundle, but also reveal  precisely what a connection on such a bundle must be.

We start with the local aspects of transport functors with $\mathrm{Gr}$-structure by considering smooth functors
\begin{equation}
\label{grpdfunctor}
F: \mathcal{P}_1(X) \to \mathrm{Gr}\text{.}
\end{equation} 
Our aim is to obtain a correspondence between such functors and certain 1-forms, generalizing the one derived in Section \ref{sec7}.
If we denote the objects of $\mathrm{Gr}$  by $\mathrm{Gr}_0$ and the morphisms by $\mathrm{Gr}_1$,  $F$ defines in the first place a smooth map $f: X \to \mathrm{Gr}_0$. Using the technique introduced in Section \ref{sec7}, we obtain further a  1-form $A$ on $X$ with values in the vector bundle $f^{*}\id^{*}T\mathrm{Gr}_1$ over $X$. Only the fact that $F$ respects targets and sources imposes two new  conditions: 
\begin{equation*}
f^{*}\mathrm{d}s \circ A = 0
\quad\text{ and }\quad
f^{*}\mathrm{d}t \circ A + \mathrm{d}f = 0\text{.}
\end{equation*}
Here we regard $\mathrm{d}f$ as a 1-form on $X$ with values in $f^{*}T\mathrm{Gr}_0$, and $\mathrm{d}s$ and $\mathrm{d}t$ are the differentials of the source and target maps.

Now we recall that the \emph{Lie algebroid} $E$ of $\mathrm{Gr}$  is the vector bundle
\begin{equation*}
E := \id^{*}\mathrm{ker}(\mathrm{d}s)
\end{equation*}
 over $\mathrm{Gr}_0$
where $\id:\mathrm{Gr}_0 \to \mathrm{Gr}_1$ is the identity embedding. The \emph{anchor} is the morphism $a := \mathrm{d}t: E \to T\mathrm{Gr}_0$ of vector bundles over $\mathrm{Gr}_0$. Using this terminology, we see that
the smooth functor \erf{grpdfunctor} defines a smooth map $f: X \to \mathrm{Gr}_0$ plus a 1-form $A \in \Omega^1(X,f^{*}E)$ such that $f^{*}a \circ A + \mathrm{d}f=0$.

In order to deal with smooth natural transformations, we introduce the following notation. We denote by
\begin{equation*}
c: \mathrm{Gr}_1 \lli s \times_t \mathrm{Gr}_1 \to \mathrm{Gr}_1: (h,g) \mapsto h \circ g
\end{equation*}
the composition in the Lie groupoid $\mathrm{Gr}$, and for $g:x \to y$ a morphism by
\begin{equation*}
r_g: s^{-1}(y) \to s^{-1}(x): h \mapsto h \circ g
\end{equation*}
the composition by $g$ from the right. Notice that $c$ and $r_g$ are smooth maps.  
It is straightforward to check that one has a well-defined map
\begin{equation*}
\mathrm{AD}_g: T_g \Gamma_1 \lli{\mathrm{d}s} \times_{a} E_{s(g)} \to E_{t(g)}
\end{equation*}
which is defined by
\begin{equation}
\label{defAD}
\mathrm{AD}_g(X,Y) := \mathrm{d}r_{g^{-1}}|_g (\mathrm{d}c|_{g,\id_{s(g)}}(X,Y))\text{.}
\end{equation}
For example, if $\mathrm{Gr}=\mathcal{B}G$ for a Lie group $G$, the Lie algebroid is the trivial bundle $E = \Gamma_0 \times \mathfrak{g}$, the composition $c$ is the multiplication of $G$, and \erf{defAD} reduces to
\begin{equation*}
\mathrm{AD}_g(X,Y) = \bar\theta_g(X) +  \mathrm{Ad}_g(Y) \in \mathfrak{g}\text{.}
\end{equation*}

Suppose now that
\begin{equation*}
\eta: F \Rightarrow F'
\end{equation*}
is a smooth natural transformation between smooth functors $F$  and $F'$ which correspond to pairs $(f,A)$ and $(f',A')$, respectively. It defines a smooth  map $g: X \to \mathrm{Gr}_1$ such that 
\begin{equation}
\label{nattransts}
s \circ g = f
\quad\text{ and }\quad
t \circ g = f'\text{.}
\end{equation}
Generalizing Lemma \ref{lem14}, the naturality of $\eta$ implies additionally
\begin{equation}
\label{12}
A' + \mathrm{AD}_g(\mathrm{d}g,-A)=0\text{.}
\end{equation}

The structure obtained like this forms a category $Z^1_X(\mathrm{Gr})$ of \emph{$\mathrm{Gr}$-connections}: its objects are pairs $(f,A)$ of smooth functions $f: X \to \mathrm{Gr}_0$ and 1-forms $A\in \Omega^1(X,f^{*}E)$ satisfying $f^{*}\mathrm{d}t \circ A+\mathrm{d}f=0$, and its morphisms are smooth maps $g:X \to \mathrm{Gr}_1$ satisfying \erf{nattransts} and  \erf{12}. The category $Z^1_X(\mathrm{Gr})$ generalizes the category of $G$-connections from Definition \ref{def11} in the sense that  $Z_X^1(\mathcal{B}G)=Z_X^1(G)$ for $G$ a Lie group. We obtain the following generalization of Proposition \ref{th2}.
\begin{proposition}
\label{prop8}
There is a canonical isomorphism of categories
\begin{equation*}
\mathrm{Funct}^{\infty}(X,\mathrm{Gr}) \cong Z_X^1(\mathrm{Gr})\text{.}
\end{equation*}
\end{proposition} 

We remark that examples of smooth of functors with values in a Lie groupoid naturally appear in the discussion of transgression to loop spaces, see Section 4 of the forthcoming paper \cite{schreiber5}.

\medskip

Now we come to the global aspects of transport functors with $\mathrm{Gr}$-structure. We introduce the category of $\mathrm{Gr}$-torsors as an interesting target category of such transport functors. 
A \textit{smooth  $\mathrm{Gr}$-manifold} \cite{moerdijk} is a triple $(P,\lambda,\rho)$ consisting a smooth manifold $P$, a surjective submersion $\lambda: P \to \mathrm{Gr}_0$ and
a smooth map
\begin{equation*}
\rho: P\; {}_\lambda\!\! \times_t \mathrm{Gr}_1
\to P
\end{equation*}
such that 
\begin{enumerate}
\item 
$\rho$ respects $\lambda$ in the sense that
$\lambda(\rho(p,\varphi))=s(\varphi)$ for all $p\in P$ and
$\varphi\in\mathrm{Gr}_1$ with $\lambda(p)=t(\varphi)$, 

\item 
$\rho$ respects the composition $\circ$ of morphisms of $\mathrm{Gr}$.
\end{enumerate}
A \emph{morphism between $\mathrm{Gr}$-manifolds} is a smooth map $f: P \to P'$ which
respects $\lambda$, $\lambda'$ and $\rho$, $\rho'$.
A \emph{$\mathrm{Gr}$-torsor} is a $\mathrm{Gr}$-manifold for which $\rho$ acts in a free and transitive way. $\mathrm{Gr}$-torsors form a category denoted $\mathrm{Gr}\text{-}\mathrm{Tor}$. For a fixed object
$X \in \mathrm{Obj}(\mathrm{Gr})$, $P_{X}:=t^{-1}(X)$ is a $\mathrm{Gr}$-torsor
with $\lambda=s$ and $\rho=\circ$. Furthermore, a morphism $\varphi:X \to
Y$ in $\mathrm{Gr}$ defines a morphism $P_X \to P_Y$ of $\mathrm{Gr}$-torsors.
Together, this defines a functor
\begin{equation}
\label{60}
  i_{\mathrm{Gr}} : \mathrm{Gr} \to \mathrm{Gr}\text{-}\mathrm{Tor}\text{.}
\end{equation}

The functor \erf{60} allows us to study transport functors
\begin{equation*}
\mathrm{tra}: \mathcal{P}_1(M) \to \mathrm{Gr}\text{-}\mathrm{Tor}
\end{equation*}
with $\mathrm{Gr}$-structure. By a straightforward adaption of the Second Proof of Theorem \ref{th1} one can construct the total space $P$ of a fibre bundle over $M$ from the transition function $\tilde g: Y^{[2]} \to \mathrm{Gr}_1$ of $\mathrm{tra}$, in such a way that $P$ is fibrewise a $\mathrm{Gr}$-torsor. More precisely, we reproduce the following definition.

\begin{definition}[\cite{moerdijk}]
A \emph{principal $\mathrm{Gr}$-bundle} over $M$ is a $\mathrm{Gr}$-manifold $(P,\lambda,\rho)$ together with a smooth map $p:P \to M$ which is preserved by the action, such that there exists a surjective submersion $\pi:Y \to M$ with a smooth map $f:Y \to \mathrm{Gr}_1$ and a morphism
\begin{equation*}
\phi: P \times_M Y \to Y \!\!\lli{f}\times_{t}\mathrm{Gr}_1
\end{equation*}
of $\mathrm{Gr}$-manifolds that preserves the projections to $Y$. 
\end{definition}

Here we have used surjective submersions instead of open covers, like we already did for principal bundles (see Section \ref{sec5_1}).
Principal $\mathrm{Gr}$-bundles over $M$ form a category denoted $\mathrm{Gr}\text{-}\mathfrak{Bun}^{\nabla}(M)$, whose morphisms are morphisms of $\mathrm{Gr}$-manifolds that preserve the projections to $M$.

The descent data of the transport functor $\mathrm{tra}$ not only consists of the transition function $\tilde g$ but also of a smooth functor $\mathrm{triv}: \mathcal{P}_1(Y) \to \mathrm{Gr}$. Now, Proposition \ref{prop8} \emph{predicts} the notion of a connection 1-form on a principal $\mathrm{Gr}$-bundle:

\begin{definition}
\label{grpdbcon}
Let $\mathrm{Gr}$ be a Lie groupoid and $E$ be its Lie algebroid. 
A \emph{connection on a principal $\mathrm{Gr}$-bundle} $P$ is a 1-form $\omega\in\Omega^1(P,\lambda^{*}E)$ such that 
\begin{equation*}
\lambda^{*}\mathrm{d}t \circ \omega + \mathrm{d}\lambda=0
\quad\text{ and }\quad
p_1^{*}\omega + \mathrm{AD}_g(\mathrm{d}g,-\rho^{*}\omega)=0\text{,}
\end{equation*}
where $\lambda: P \to \mathrm{Gr}_0$ and $\rho: P \lli{\lambda}\times_t \mathrm{Gr}_1 \to P$ are the structure of the $\mathrm{Gr}$-manifold $P$, $p_{1}$ and $g$ are the projections to $P$ and $\mathrm{Gr}_1$, respectively.  
\end{definition}

By construction, we have
\begin{theorem}
\label{th5}
There is a canonical equivalence of categories
\begin{equation*}
\mathrm{Gr}\text{-}\mathfrak{Bun}^{\nabla}(M) \cong \mathrm{Trans}_{\mathrm{Gr}}(M,\mathrm{Gr}\text{-}\mathrm{Tor})\text{.}
\end{equation*}
\end{theorem}

Indeed, choosing a local trivialization $(Y,f,\phi)$ of a principal $\mathrm{Gr}$-bundle $P$, one obtains a section $s:Y \to P:y \mapsto p(\phi^{-1}(y,\id_{f(y)}))$. This section satisfies $\lambda \circ s = f$, so that the pullback of a connection 1-form $\omega \in \Omega^1(P,\lambda^{*}E)$ along $s$ is a 1-form $A:=s^{*}\omega\in\Omega^1(Y,f^{*}E)$. The first  condition in Definition \ref{grpdbcon} implies that $(A,f)$ is an object in $Z_Y^1(\mathrm{Gr})$, and thus by Proposition \ref{prop8} a smooth functor $\mathrm{triv}: \mathcal{P}_1(Y) \to \mathrm{Gr}$. The second condition implies that the transition function $\tilde g:Y^{[2]} \to \mathrm{Gr}$ defined by $s(\pi_1(\alpha))=\rho(s(\pi_2(\alpha)),\tilde g (\alpha))$ is a morphism in $Z_{Y^{[2]}}^1(\mathrm{Gr})$ from $\pi_1^{*}F$ to $\pi_2^{*}F$. All together, this is descent data for a transport functor on $M$ with $\mathrm{Gr}$-structure.

\medskip

We remark that this automatically induces a notion of parallel transport for  a connection $A$ on a principal $\mathrm{Gr}$-bundle $P$:  let $\mathrm{tra}_{P,A}: \mathcal{P}_1(M) \to \mathrm{Gr}\text{-}\mathrm{Tor}$ be the transport corresponding to $(P,A)$ under the equivalence of Theorem \ref{th5}. Then, the parallel transport of $A$ along a path $\gamma:x \to y$ is the $\mathrm{Gr}$-torsor morphism
\begin{equation*}
\mathrm{tra}_{P,A}(\gamma): P_x \to P_y\text{.}
\end{equation*}

\medskip

\def\act{\!\!/\!\!/\!}

In the remainder of this section we discuss a class of groupoid bundles with connection related to action groupoids. We recall that for $V$ a complex vector space with an action of a Lie group $G$, the action groupoid $V\act G$ has $V$ as its objects and $G \times V$ as its morphisms. The source map is the projection to $V$, and the target map is the action $\rho: G \times V \to V$. Every action groupoid $V \act G$ comes with a canonical functor
\begin{equation*}
i_{V \act G} : V \act G \to \mathrm{Vect}_{*}(\C)
\end{equation*}
to the category of pointed complex vector spaces, which sends an object $v\in V$ to the pointed vector space $(V,v)$ and a morphism $(v,g)$ to the linear map $\rho(g,-)$, which respects the base points.

\begin{proposition}
A transport functor $\mathrm{tra}: \mathcal{P}_1(M) \to \mathrm{Vect}_{*}(\C)$ with $V\act G$-structure is  a complex vector bundle over $M$ with structure group $G$ and a smooth flat section.  
\end{proposition}

\begin{proof}
We consider the strictly commutative diagram
\begin{equation*}
\alxydim{@C=1.5cm}{V \act G \ar[r]^-{i_{V\act G}}  \ar[d]_{\mathrm{pr}} & \mathrm{Vect}_{*}(\C) \ar[d]^{f} \\ \mathcal{B}G \ar[r]_-{\rho} & \mathrm{Vect}(\C)}
\end{equation*}
of functors, in which $f$ is the functor that forgets the base point, $\mathrm{pr}$ is the functor which sends a morphism $(g,v)$ in $V \act G$ to $g$, and $\rho$ is the given representation. The diagram shows that the composition
\begin{equation*}
f \circ \mathrm{tra}: \mathcal{P}_1(M) \to \mathrm{Vect}(\C)
\end{equation*} 
is a transport functor with $\mathcal{B}G$-structure, and hence the claimed vector bundle $E$ by (a slight generalization of) Theorem \ref{th7}. Remembering the forgotten base point defines a natural transformation
\begin{equation*}
\eta: \trivlin_{\C} \to f \circ \mathrm{tra}\text{.}
\end{equation*}
If we regard the identity transport functor $\trivlin_{\C}$ as a transport functor with $\mathcal{B}G$-structure, the natural transformation $\eta$ becomes a morphism of transport functors with $\mathcal{B}G$-structure, and thus defines by Lemma \ref{lem9} a smooth flat section in $E$.
\end{proof}

\section{Generalizations and further Topics}

\label{sec1}

The concept of transport functors has generalizations
in many aspects, some of which we want to outline in this section. 

\subsection{Transport $n$-Functors}

The motivation to write this article was to find a formulation of parallel
transport along curves, which can be generalized to higher dimensional parallel
transport. Transport functors have a natural generalization to transport
$n$-functors. In particular the case $n=2$ promises
relations between transport 2-functors and gerbes with connective structure
\cite{baez2},
similar to the relation between transport 1-functors and bundles with connections
presented in Section \ref{sec4}. 
 We address these issues in a  further publication \cite{schreiber2}.

\medskip

Let us briefly describe the generalization of the concept of transport functors
to transport $n$-functors. The first generalization is that of the path groupoid $\mathcal{P}_1(M)$
to a path $n$-groupoid $\mathcal{P}_n(M)$. Here, $n$-groupoid means that
every $k$-morphism is an equivalence,
i.e. invertible up to  $(k+1)$-isomorphisms.
 The set of objects is again
the manifold $M$,  the $k$-morphisms are  smooth
maps $[0,1]^{k} \to M$ with sitting instants on each boundary of the $k$-cube,
and the top-level morphisms $k=n$ are additionally
taken up to thin homotopy in the appropriate sense. 

We then consider $n$-functors 
\begin{equation}
\label{26}
F: \mathcal{P}_n(M) \to T
\end{equation}
from the path $n$-groupoid $\mathcal{P}_n(M)$ to some target $n$-category $T$. Local
trivializations of such $n$-functors are considered with respect to an $n$-functor
$i: \mathrm{Gr} \to T$, where $\mathrm{Gr}$ is a Lie $n$-groupoid, and to
a surjective
submersions $\pi:Y \to M$. A $\pi$-local
$i$-trivialization then consists of an $n$-functor $\mathrm{triv}:\mathcal{P}_n(Y)
\to \mathrm{Gr}$ and an equivalence
\begin{equation}
\label{13}
\alxydim{@C=1.5cm@R=1.5cm}{\mathcal{P}_n(Y) \ar[r]^{\pi_{*}} \ar[d]_{\mathrm{triv}} & \mathcal{P}_n(M)
\ar@{=>}[dl]|*+{t}
\ar[d]^{F} \\  \mathrm{Gr} \ar[r]_{i} & T}
\end{equation}
of $n$-functors. Local trivializations   lead to an $n$-category
$\trans{i}{i}{n}$ of descent data, which are descent $n$-categories
in the sense of \cite{street}, similar to Remark \ref{re1} for $n=1$. 

The category $\trans{i}{i}{n}$ has  a natural notion of smooth
objects and smooth $k$-morphisms. Then, $n$-functors (\ref{26}) which allow
 local trivializations with smooth descent data will be called transport
$n$-functors, and form an $n$-category $\transport{i}{n}{\mathrm{Gr}}{T}$. 

\medskip
 
In the case $n=1$, the procedure described above reproduces the framework
of transport functors described in this article. The case $n=2$ will be considered in detail
in two forthcoming papers. First we settle the local aspects: we derive a correspondence between smooth 2-functors and differential 2-forms (Theorem 2.20 in \cite{schreiber5}). Then we continue with the global aspects in  \cite{schreiber2}. 

As a further example, we now describe  the case $n=0$.
Note that a $0$-category is
a set, a Lie $0$-groupoid is a smooth manifold, and a $0$-functor is a map. To start with, we have the set $\mathcal{P}_0(M)=M$, a set $T$, a smooth manifold $G$ and an injective map $i: G \to T$.  Now we consider
maps $F: M \to T$. Following the general concept, such a map is $\pi$-locally $i$-trivializable,
if there exists a map $\mathrm{triv}:Y \to G$ such that the diagram
\begin{equation*}
\alxydim{@=1.3cm}{Y \ar[r]^{\pi} \ar[d]_{\mathrm{triv}} & M
\ar[d]^{F} \\  G \ar[r]_{i} & T}
\end{equation*}
is commutative. Maps $F$ together with $\pi$-local $i$-trivializations form
the set $\loctrivfunct{i}{\mathrm{Gr}}{0}$. The set $\trans{i}{i}{0}$ of descent data is just
  the set of  maps $\mathrm{triv}:Y \to G$ satisfying the equation
\begin{equation}
\label{10}
\pi_1^{*}\mathrm{triv}_{i}
= \pi_2^{*}\mathrm{triv}_{i}\text{,}
\end{equation}
where we have used the notation $\pi_k^{*}\mathrm{triv}_i=i \circ \mathrm{triv}
\circ \pi_k$ from Section \ref{sec2}. 
It is easy to see that every $\pi$-local $i$-trivialization $\mathrm{triv}$
of a map $F$ satisfies
this condition. This defines the map
 \begin{equation*}
\ex{\pi}:\loctrivfunct{i}{\mathrm{Gr}}{0} \to \trans{i}{i}{0}\text{.}
\end{equation*}
Similar to Theorem \ref{th3} in the case $n=1$, this  is indeed a bijection: every function $\mathrm{triv}:Y \to G$ satisfying
(\ref{10}) with $i$ injective factors through $\pi$. Now  it is easy to say when an element in $\trans{i}{i}{0}$ is called
smooth: if and only if the map $\mathrm{triv}:Y \to G$ is smooth.  Such maps
form the set $\transsmooth{i}{0}{i}$, which in turn defines the set  $\transport{i}{0}{G}{T}$ of transport $0$-functors with $G$-structure. Due to (\ref{10}),  there  is a canonical bijection $\transsmooth{i}{0}{i} \cong
C^{\infty}(M,G)$.
So, we have
\begin{equation*}
\transport{i}{0}{G}{T} \cong C^{\infty}(M,G)\text{,}
\end{equation*}
in other words: transport 0-functors  on $M$ with $G$-structure are
smooth functions from $M$ to $G$. 

\medskip

Let us revisit Definition \ref{def2} of the category $\transsmooth{i}{1}{i}$ of smooth descent data, which  now can equivalently be reformulated as follows:

\begin{itemize}
\item[]\it
Let $\mathrm{Gr}$ be a Lie groupoid and let $i:\mathrm{Gr} \to T$ be a functor. An object $(\mathrm{triv},g)$ in $\trans{i}{i}{1}$ is called \emph{smooth},
if the functor $\mathrm{triv}:\mathcal{P}_1(Y) \to \mathrm{Gr}$ is smooth
the sense of Definition \ref{def1}, and if the natural equivalence $g:Y^{[2]}
\to \mathrm{Mor}(T)$ is a transport $0$-functor with $\mathrm{Mor}(T)$-structure.
A morphism 
\begin{equation*}
h:(\mathrm{triv},g) \to (\mathrm{triv}',g')
\end{equation*}
between smooth objects is called \emph{smooth}, if $h: Y \to \mathrm{Mor}(T)$ is
a transport 0-functor with $\mathrm{Mor}(T)$-structure.
\end{itemize}
This gives an outlook  how the definition of the $n$-category $\transsmooth{i}{n}{i}$ of smooth descent data will
be for higher $n$: it will  recursively use  transport $(n-1)$-functors. 

\subsection{Curvature of Transport
Functors}

\label{sec6_2}

When we describe parallel transport in terms of functors, it is a natural question
how related notions like curvature can be seen in this formulation.
Interestingly, it turns out that the curvature of a  transport functor is
a transport 2-functor. More generally, the curvature
of a transport $n$-functor is a transport $(n+1)$-functor. This becomes
evident with a view to Section \ref{sec7}, where we have related smooth
functors
and differential 1-forms. In a similar way, 2-functors can be related to  2-forms. A comprehensive discussion of the curvature of transport
functors is therefore beyond the scope of this article, and has to be postponed
until
after the discussion of transport 2-functors \cite{schreiber2}. 

We shall briefly indicate the basic ideas. We recall from Section \ref{sec10}
when a functor $F:\mathcal{P}_1(M) \to T$ is flat: if it factors through
the fundamental groupoid $\Pi_1(M)$, whose morphisms are smooth homotopy classes
of paths in $M$. In general, one can associate to a transport functor $\mathrm{tra}$ a  2-functor
\begin{equation*}
  \mathrm{curv}(\mathrm{tra}) : \mathcal{P}_{2}(M) \to \mathrm{Grpd}
\end{equation*}
into the 2-category of groupoids. This 2-functor is particularly trivial if  $\mathrm{tra}$ is flat. Furthermore, the 2-functor $\mathrm{curv}(\mathrm{tra})$
is itself flat in the sense that it factors through the fundamental 2-groupoid
of $M$: this is nothing but the Bianchi identity.

For smooth functors $F:\mathcal{P}_1(M) \to \mathcal{B}G$, which corresponding by Proposition
\ref{th2} to 1-forms $A\in\Omega^1(M,\mathfrak{g})$, it turns out that the
2-functor $\mathrm{curv}(F)$ corresponds to a 2-form $K\in\Omega^2(M,\mathfrak{g})$
which is related to $A$ by the usual equality $K=\mathrm{d}A + A \wedge A$.

\subsection{Alternatives to smooth Functors}

The definition of transport functors concentrates on the smooth aspects of
parallel transport. As we have outlined in Appendix \ref{appB}, our definition
of smooth descent data $\transsmooth{i}{1}{i}$ can be regarded as the internalization
of functors and natural transformations in the category $\smsp$ of diffeological
spaces and diffeological maps.

Simply by choosing another ambient category $C$, we obtain possibly weaker notions
of parallel transport. Of particular interest is the situation where the
ambient category is the category Top of topological spaces and continuous
maps. Indicated by results of \cite{stasheff}, one would expect
that reconstruction theorems as discussed in Section \ref{sec5} should also
exist for Top, and also for transport $n$-functors for $n>1$.
Besides, parallel transport along topological paths of bounded
variation can be defined, and is of interest for its own right, see, for example, \cite{baudoin}.

\subsection{Anafunctors}

The notion of smoothly locally trivializable functors  is
closely related to the concept of anafunctors. Following \cite{makkai1}, an anafunctor
$F: A \to B$ between categories $A$ and $B$ is a category $|F|$ together
with a functor $\tilde F: |F| \to
B$ and a surjective equivalence
$p: |F| \to A$, denoted as a diagram
\begin{equation}
\label{41}
    \alxydim{@=1.3cm}{
      |F| \ar[r]^{\tilde F} 
      \ar[d]_{p}
      & B
      \\
      A
    }
\end{equation}
called a span.
 It has been shown in  \cite{Bartels}  how
to formulate the concept of an anafunctor internally to any category $C$.  

Note that an anafunctor in $C$ gives rise to an ordinary functor $A \to B$
in $C$, if the
epimorphism $p$ has a section. 
In the category of sets, $C=\set$, every epimorphism has a section, if one
assumes the axiom of choice (this is what we do). The original motivation
for introducing anafunctors was, however, to deal with situations where one
does not assume the axiom of choice \cite{makkai1}. In the category $C=\sm$
of smooth manifolds,  surjective submersions are particular epimorphisms, as they arise for example
as   projections
of 
smooth fibre bundles. Since not every bundle has a global smooth section,
an anafunctor in $\sm$ does not produce a functor. The same applies
to the category $C=\smsp$ of diffeological spaces described in Appendix \ref{appB}.

Let us indicate how anafunctors arise
from smoothly locally trivialized functors. Let $\mathrm{tra}:\mathcal{P}_1(M) \to T$ be a transport functor with $\mathrm{Gr}$-structure. We choose a $\pi$-local $i$-trivialization
$(\mathrm{triv},t)$, whose descent data $(\mathrm{triv},g)$ is smooth.
Consider the functor
\begin{equation*}
  R_{(\mathrm{triv},g)} : \upp{M}{\pi} \to T
\end{equation*}
that we have defined in Section \ref{sec5} from this descent data. By Definition \ref{def2} of smooth descent data, the functor $\mathrm{triv}:\mathcal{P}_1(Y) \to \mathrm{Gr}$ is smooth and the natural equivalence $g$ factors through a smooth natural equivalence $\tilde
g: Y \to \mathrm{Mor}(\mathrm{Gr})$. So, the functor $R_{(\mathrm{triv},g)}$
factors through $\mathrm{Gr}$,
\begin{equation*}
R_{(\mathrm{triv},g)} =  i \circ A
\end{equation*}
for a   functor $A: \upp{M}{\pi} \to \mathrm{Gr}$. In fact, the category $\upp{M}{\pi}$
can be considered as a category internal to $\smsp$, so that the
functor $A$ is internal to $\smsp$ as described in Appendix
\ref{appB}, Proposition \ref{prop4}
ii).  Hence the reconstructed functor yields a span
\begin{equation*}
 \alxydim{@=1.3cm}{
   \upp{M}{\pi}
   \ar[d]_{\uppp}
   \ar[r]^-{A}
   &
   \mathrm{Gr}
   \\
   \mathcal{P}_1(M)\text{,}
 }
\end{equation*}
internal to $\smsp$,
i.e. an anafunctor $\mathcal{P}_1(M) \to \mathrm{Gr}$. Because the epimorphism
$\uppp$ is not invertible in $\smsp$,  we do not get
an ordinary functor $\mathcal{P}_1(M) \to \mathrm{Gr}$ internal to $\smsp$: the weak inverse functor $s: \mathcal{P}_1(M)
\to \upp{M}{\pi}$ we have constructed in Section \ref{sec5} is not internal to $\smsp$.

\paragraph*{Acknowledgements}
We thank Bruce Bartlett, Uwe Semmelmann, Jim Stasheff and Danny Stevenson for helpful correspondences
and
Christoph Schweigert for helpful discussions. U.S.  thanks John Baez for
many valuable suggestions and  discussions.  We acknowledge support from the Sonderforschungsbereich
\quot{Particles, Strings and the Early Universe - the Structure of Matter and 
Space-Time}. 

\begin{appendix}

\section{More Background}

\subsection{The universal Path Pushout}

\label{app2}

Here we motivate Definition \ref{def21} of the groupoid $\upp{M}{\pi}$.
Let $\pi:Y \to M$ be a surjective submersion. 
A \emph{path pushout  of  $\pi$} is a  triple $(A,b,\nu)$ consisting of a groupoid $A$, a functor $b: \mathcal{P}_1(Y)
\to A$ and a natural equivalence $\nu: \pi_1^{*}b \to \pi_{2}^{*}b$
with
\begin{equation*}
\pi_{13}^{*}\nu = \pi_{23}^{*}\nu \circ \pi_{12}^{*}\nu\text{.}
\end{equation*}
A morphism 
\begin{equation*}
(R,\mu): (A,b,\nu) \to (A',b',\nu')
\end{equation*}
between path pushouts is a functor $R: A \to A'$ and a natural equivalence
$\mu: R \circ b \to b'$ such that
\begin{equation}
\label{14}
\alxy{\mathcal{P}_1(Y^{[2]}) \ar[r]^{\pi_1} \ar[d]_{\pi_2} & \mathcal{P}_1(Y)
\ar@{=>}[dl]|{\nu}
\ar@/^1.3pc/[ddr]^{b'}="1" \ar@{=>}"1";[d]|{\mu^{-1}}
\ar[d]^{b} & \\ \mathcal{P}_1(Y) \ar[r]_{b} \ar@/_1.3pc/[drr]_{b'}="2" & A \ar@{=>}"2"|{\mu} \ar[rd]^{R}
& \\ && A'}
=
\alxy{\mathcal{P}_1(Y^{[2]}) \ar[r]^{\pi_1} \ar[d]_{\pi_2} & \mathcal{P}_1(Y)
\ar@{=>}[dl]|{\nu'}
 \ar[d]^{b'}  \\ \mathcal{P}_1(Y) \ar[r]_{b'}  & A'\text{.}}
\end{equation}

Among all path pushouts of $\pi$ we  distinguish some having a universal
property.

\begin{definition}
\label{def4}
A path pushout  $(A,b,\nu)$ is \emph{universal}, if, 
given any other path pushout  $(T,F,g)$, there exists a morphism $(R,\mu):
(A,b,\nu) \to (T,F,g)$ such that, given any other such morphism $(R',\mu')$, there
is a unique natural equivalence  $r: R \to R'$ with 
\begin{equation}
\label{15}
\alxydim{@C=2cm}{\mathcal{P}_1(Y) \ar@/_1.5pc/[dr]_<<<<{F}="3" \ar[r]^{b} & A \ar@{=>}"3"|{\mu'} \ar@/^1.2pc/[d]^{R}="1" \ar@/_1.2pc/[d]_{R'}="2" \ar@{=>}"1";"2"|{r} \\ & T}
=
\alxydim{@C=2cm}{\mathcal{P}_1(Y) \ar@/_1.5pc/[dr]_{F}="3" \ar[r]^{b} & A \ar@{=>}"3"|{\mu} \ar[d]^{R}="1"  \\ & T\text{.}}
\end{equation}
\end{definition}

Now we show how two path pushouts having both the universal property, are
related.

\begin{lemma}
\label{lem1}
Given two universal path pushouts $(A,b,\nu)$ and $(A',b',\nu')$ of the
same surjective submersion $\pi:Y \to M$, there is an  equivalence of categories $a: A \to A'$.
\end{lemma}

\begin{proof}
We use the universal properties of both triples applied to each other. We
obtain  two choices of morphisms $(R,\mu)$ and $(\tilde R, \tilde \mu)$, namely 
\begin{equation*}
\alxydim{@C=0.8cm}{\mathcal{P}_1(Y^{[2]}) \ar[r]^{\pi_1} \ar[d]_{\pi_2} & \mathcal{P}_1(Y)
\ar@{=>}[dl]|{\nu}
\ar[d]^{b} \ar@/^1.2pc/[ddr]^{b'}="1" \ar@/^2.5pc/[dddrr]^{b}="2" && \\ \mathcal{P}_1(Y) \ar[r]_{b}
\ar@/_1.2pc/[drr]_{b'}="3" \ar@/_2.5pc/[ddrrr]_{b}="4" & A \ar@{=>}"3"|{\mu} \ar@{<=}"1"|{\mu^{-1}} \ar[dr]_{a} && \\ && A' \ar@{=>}"4"|{\mu'} \ar@{<=}"2"|{\mu'^{-1}} \ar[dr]_{a'} & \\
&&& A }
\text{ and }
\alxydim{@C=0.7cm}{\mathcal{P}_1(Y^{[2]}) \ar[r]^{\pi_1} \ar[d]_{\pi_2} & \mathcal{P}_1(Y)
\ar@{=>}[dl]|{\nu}
\ar@/^1.3pc/[ddr]^{b}="1" 
\ar[d]^{b} & \\ \mathcal{P}_1(Y) \ar[r]_{b} \ar@/_1.3pc/[drr]_{b}="2" & A  \ar[rd]^{\id_A}
& \\ && A}
\end{equation*} 
The unique natural transformation we get from the universal property is here $r: a' \circ a \to \id_A$.  Doing the same thing in the other order, we obtain
a unique natural transformation $r': a \circ a' \to \id_{A'}$. Hence $a: A \to
A'$ is an equivalence of categories.
\end{proof}

We also need

\begin{lemma}
\label{lem4}
Let $(A,b,\nu)$ be a universal path pushout of $\pi$, let $(T,F,g)$ and $(T,F',g')$
two other path pushouts and let $h: F \to F'$ be a natural transformation
with $\pi_2^{*}h \circ g = \pi_1^{*}h
\circ g'$. For any choice of morphisms 
\begin{equation*}
(R,\mu):(A,b,\nu) \to (T,F,g)
\quad\text{ and }\quad
(R',\mu'):(A,b,\nu) \to (T,F',g')
\end{equation*}
there is a unique natural transformation
$r: R \to R'$
with $\mu \bullet (\id_b \circ r) = \mu'$.
\end{lemma}

\begin{proof}
Note that the natural equivalence $h$ defines a morphism 
\begin{equation*}
(\id_T,h): (T,F,g)
\to (T,F',g')
\end{equation*}
of path pushouts. The composition $(\id_T,h) \circ (R,\mu)$
gives a morphism 
\begin{equation*}
(R,h \circ \mu):(A,b,\nu) \to (T',F',g')\text{.}
\end{equation*}
Since $(R',\mu')$
is universal, we obtain a unique natural transformation $r: R \to R'$. 
\end{proof}

Now consider the groupoid $\upp{M}{\pi}$ from Definition \ref{def21}, together
with the inclusion functor $\iota: \mathcal{P}_1(Y) \to \upp{M}{\pi}$ and the
identity $\id_{Y^{[2]}}:\pi_1^{*}\iota \to \pi_2^{*}\iota$ whose component at a point
$\alpha\in Y^{[2]}$ is the morphism $\alpha$ in $\upp{M}{\pi}$. Its commutative diagram
follows from relations (1) and (2), depending on the type of morphism you
apply it to. Its cocycle condition follows from (3). So, the triple $(\upp{M}{\pi},\iota,\id_{Y^{[2]}})$
is a path pushout.

\begin{lemma}
\label{lem3}
The triple $(\upp{M}{\pi},\iota,\id_{Y^{[2]}})$ is universal.
\end{lemma}

\begin{proof}
Let $(T,F,g)$ any path pushout. We construct the morphism
\begin{equation*}
(R,\mu): (\upp{M}{\pi},\iota,\id_{Y^{[2]}})\to (T,F,g)
\end{equation*}
as follows. The functor 
\begin{equation*}
R:
\upp{M}{\pi} \to T
\end{equation*}
 sends an object $x \in Y$ to $F(x)$, a morphism $\gamma:x
\to y$ to $F(\gamma)$ and a morphism $\alpha$ to $g(\alpha)$.
This definition is well-defined under the relations among the morphisms:
(1) is the commutative diagram for the natural transformation $g$, (2) is the
cocycle condition for $g$ and (3)  follows from the latter since $g$ is
invertible. The natural equivalence $\mu: R \circ \iota \to F$ is the identity. By definition equation (\ref{14}) is satisfied, so that $(R,\mu)$ is
a morphism of path pushouts. Now we assume that there is another morphism
$(R',\mu')$. The component of the natural equivalence $r: R \to R'$ at a point  $x\in Y$ is 
$\mu'^{-1}(x)$, its naturality with respect to a morphisms $\gamma:x
\to y$ is then just the one of $\mu'$, and with respect to morphisms
$\alpha\in Y^{[2]}$ comes from condition (\ref{14}) on morphisms of path
pushouts. It also satisfies the equality (\ref{15}). Since this equation
already determines $r$, it is unique.    
\end{proof}

Notice that the construction of the functor $R$ reproduces Definition \ref{def13}.
Let us finally apply Lemma \ref{lem4} to the universal path pushout $(\upp{M}{\pi},\iota,\id_{Y^{[2]}})$.
Given the two functors $F,F': \mathcal{P}_1(Y) \to T$, the natural transformation
$h: F \to F'$, and the universal morphisms
$(R,\mu)$ and $(R',\mu')$ as constructed in the proof of Lemma \ref{lem3},
the natural transformation $r: R \to R'$ has the component $h(x)$ at $x$. This reproduces Definition \ref{def14}. 

\subsection{Diffeological Spaces and smooth Functors}

\label{appB}

This section puts Definition \ref{def1} of a smooth
functor into the wider perspective of functors internal to some category
$C$, here the category $\smsp$ of diffeological spaces \cite{chen1,souriau1}.
Diffeological spaces  generalize  the concept of a smooth manifold. 
While the set $\sm(X,Y)$ of smooth maps between smooth manifolds
$X$ and $Y$ does not
form, in general, a smooth manifold itself, the set $\smsp(X,Y)$
of diffeological maps between diffeological spaces is again  a diffeological space in a canonical way. 
In other words, the category $\smsp$ of diffeological spaces  is
closed. 
\begin{definition}
\label{def18}
A \emph{diffeological space} is a set $X$ together with a collection of  plots: maps
\begin{equation*}
c : U \to X\text{,}
\end{equation*}
each of them defined on an open subset $U \subset \R^{k}$ for some $k\in
\N_0$, such that
\begin{itemize}
\item[a)]
for any plot $c: U \to X$ and any smooth function $f:V \to U$ also 
\begin{equation*}
c \circ f: V \to X
\end{equation*}
is a plot. 
\item[b)]
every constant map $c:U \to X$ is a plot.
\item[c)]
if $f:U \to X$ is a map defined on $U \subset \R^k$ and  $\lbrace U_i \rbrace_{i\in I}$ is an open cover of $U$ for which all restrictions $f|_{U_{i}}$ are plots of $X$, then also $f$ is a plot.

\end{itemize}
A \emph{diffeological map} between diffeological spaces  $X$ and $Y$ is a map $f:X \to Y$ such that
for every plot $c: U \to X$ of $X$ the map $f \circ c:U \to Y$ is a plot of $Y$. The set of all diffeological maps is denoted by $\smsp(X,Y)$. 
\end{definition}

\medskip

In fact Chen  originally used convex subsets $U \subset \R^k$ instead of open ones, but this will not be of any importance for this article. For a comparison of various concepts see \cite{kriegl1}. 
The following examples of diffeological spaces are important for us. 
\begin{itemize}
\item[(1)]
First of all, every smooth manifold is a diffeological space, the plots being all
smooth maps defined on all open subset of all $\R^{n}$. A map between two manifolds is smooth if and only if it is diffeological. 
\item[(2)] 
For diffeological spaces $X$ and $Y$ the space $\smsp(X,Y)$ of all
diffeological maps
from $X$ to $Y$ is a diffeological space in the following way: a map
\begin{equation*}
c: U \to \smsp(X,Y)
\end{equation*}
is a plot if and only if for any plot $c':V \to X$ of $X$ the composite 
\begin{equation*}
\alxydim{@C=1.5cm}{U \times V \ar[r]^-{c \times c'} & \smsp(X,Y) \times X \ar[r]^-{\mathrm{ev}}
& Y}
\end{equation*}
is a plot of $Y$. Here, $\mathrm{ev}$ denotes the evaluation map $\mathrm{ev}(f,x):=f(x)$.
\item[(3)]
Every subset $Y$ of a diffeological space $X$ is a diffeological space: its plots are those plots of $X$ whose image is contained in $Y$.

\item[(4)]
For a diffeological space $X$, a set $Y$ and a map  $p :X \to Y$, the set $Y$ is
also a diffeological space: a map $c:U \to Y$ is a plot if and only if there exists a cover of $U$ by open sets $U_{\alpha}$ together with plots $c_{\alpha}:U_{\alpha} \to X$ of $X$ such that $c|_{U_{\alpha}}=p\circ c_{\alpha}$.

\item[(5)]
Combining (1) and (2) we obtain the following important example: for smooth manifolds $X$ and $Y$ the space $C^{\infty}(X,Y)$ of smooth maps
from $X$ to $Y$ is a diffeological space in the following way: a map
\begin{equation*}
c: U \to C^{\infty}(X,Y)
\end{equation*}
is a plot if and only if 
\begin{equation*}
\alxydim{@C=1.5cm}{U \times X \ar[r]^-{c \times \id_X} & C^{\infty}(X,Y) \times X \ar[r]^-{\mathrm{ev}}
& Y}
\end{equation*}
is a smooth map. This applies for instance to the free loop space $LM=C^{\infty}(S^1,M)$.

\item[(6)]
Combining (3) and (5), the based loop space $L_xM$ and the path space $PM$ of a smooth manifold are diffeological spaces. 

\item[(7)]
Combining (4) and (6) applied to the projection $\mathrm{pr}:PM \to P^1M$ to thin
homotopy classes of paths, $P^1M$ is a diffeological space. In the same way, the thin homotopy group $\pi_1^1(M,x)$ is a diffeological space. 
\end{itemize}

From Example (7) we see that diffeological
spaces arise naturally in the setup
of transport functors introduced
in this article.

\begin{proposition}
\label{prop5}
During this article, we encountered two examples of diffeological maps:
\begin{enumerate}
\item[i)]
A Wilson line 
\begin{equation*}
\mathcal{W}^{F,i}_{x_1,x_2}: \mathrm{Mor}_{\mathcal{P}_1(M)}(x_1,x_2)
\to \mathrm{Mor}_{\mathrm{Gr}}(G_1,G_2)
\end{equation*}
is smooth in the sense of Definition
\ref{def19} if and only if it is
diffeological.
\item[ii)] 
A group homomorphism $\mathcal{H}:\pi_1^1(M,x)
\to G$ is a holonomy map in the
sense of Definition \ref{def20},
if and only if it is diffeological.
\end{enumerate}
\end{proposition}

Diffeological spaces and diffeological maps form a category
$\smsp$ in which we can
 internalize categories and functors. Examples of such categories are:

\begin{itemize}
\item 
the path groupoid 
$\mathcal{P}_1(M)$: its set of objects is the smooth manifold $M$,
which is by example (1) a diffeological space. Its set of morphisms $P^1X$  is
 a diffeological space by example (7).
\item
the universal path pushout $\upp{M}{\pi}$
of a surjective submersion $\pi:Y
\to M$: its set of objects 
is the smooth manifold $Y$, and
hence a diffeological space. A map
\begin{equation*}
    \phi : U \to \mathrm{Mor}(\upp{M}{\pi})
\end{equation*}
  is a plot if and only if there is a collection of plots $f_i : U \to P^1 Y$ and a collection of
smooth maps $g_i: U \to Y^{[2]}$
such that
\begin{equation*}
    g_N(x) \circ f_N(x) \cdots g_2(x) \circ  f_2(x) \circ g_1(x) \circ f_1(x) = \phi(x)\text{.}
\end{equation*}
\end{itemize}
We also
have examples of functors internal
to $\smsp$:

\begin{proposition}
\label{prop4}
During this article, we encountered two examples of functors
internal to $\smsp$:
\begin{enumerate}
\item[i)]
A functor $F: \mathcal{P}_1(M) \to \mathrm{Gr}$ is internal to $\smsp$
if and only if it is smooth in the sense of Definition \ref{def1}.

\item[ii)]
For a smooth object $(\mathrm{triv},g)$
in $\trans{i}{i}{1}$, the functor $R_{(\mathrm{triv},g)}$
factors smoothly through $i:\mathrm{Gr}\to
T$, i.e. there is a  functor
$A: \upp{M}{\pi} \to \mathrm{Gr}$ internal to $\smsp$
such that $i \circ A = R_{(\mathrm{triv},g)}$.
\end{enumerate}
\end{proposition}

\section{Postponed Proofs}

\subsection{Proof of Theorem \ref{th3}}

\label{appA}

Here we prove that the functor
\begin{equation*}
\ex{\pi}: \loctrivfunct{i}{\mathrm{Gr}}{1} \to \trans{i}{i}{1}
\end{equation*}
is an equivalence of categories. In Section \ref{sec5} we have defined a reconstruction functor $\mathrm{Rec}_{\pi}$
going in the opposite direction. Now we show that $\mathrm{Rec}_{\pi}$ is a weak inverse of $\ex{\pi}$ and thus prove that both are
 equivalences of categories. For this purpose, we show (a) the equation $\ex{\pi} \circ \mathrm{Rec}_{\pi} = \id_{\trans{i}{i}{1}}$ and (b) that there exists a natural equivalence
\begin{equation*}
\zeta : \id_{\loctrivfunct{i}{\mathrm{Gr}}{1}} \to \mathrm{Rec}_{\pi} \circ \ex{\pi}\text{.}
\end{equation*}

\medskip

To see (a), let $(\mathrm{triv},g)$ be an object in $\trans{i}{i}{1}$, and let
$\mathrm{Rec}_{\pi}(\mathrm{triv},g)=s^{*}R_{(\mathrm{triv},g)}$ be the reconstructed
functor, coming with the $\pi$-local $i$-trivialization $(\mathrm{triv},t)$ with $t:=g \circ
\iota^{*}\lambda$. Extracting descent data as described in Section \ref{sec6}, we find 
\begin{equation*}
(\pi_2^{*}t \circ \pi_1^{*}t^{-1})(\alpha) = g((\pi_2^{*}\lambda \circ \pi_1^{*}\lambda^{-1})(\alpha))
=  g(\alpha)
\end{equation*}
so that $\ex{\pi}(\mathrm{Rec}_{\pi}(\mathrm{triv},g)) = (\mathrm{triv},g)$. Similar, if $h:(\mathrm{triv},g) \to (\mathrm{triv}',g')$ is a morphism in
$\trans{i}{i}{1}$, the reconstructed natural equivalence is $\mathrm{Rec}_{\pi}(h)
:= s^{*}R_h$. Extracting descent data, we obtain for the component at a point $x\in Y$
\begin{eqnarray*}
(t' \circ \pi^{*}s^{*}R_h\circ t^{-1})(x)&=&g'(\lambda(x)) \circ R_h(s(\pi(x)))\circ g^{-1}(\lambda(x))\nonumber\\&=& g'(x,s(\pi(x))) \circ h(s(\pi(x))) \circ g^{-1}(x,s(\pi(x))) \nonumber\\ &=&h(x)
\end{eqnarray*}
 where we have used Definition \ref{def14} and the commutativity of diagram  (\ref{32}). This shows that  $\ex{\pi}(\mathrm{Rec}_{\pi}(h))=h$.

\medskip

To see (b), let $F: \mathcal{P}_1(M) \to T$ be a functor with  $\pi$-local $i$-trivialization $(\mathrm{triv},t)$. Let us first describe the functor
\begin{equation*}
\mathrm{Rec}_{\pi}(\ex{\pi}(F)): \mathcal{P}_1(M) \to T\text{.}
\end{equation*}
We extract descent data $(\mathrm{triv},g)$ in $\trans{i}{i}{1}$ as described
in Section $\ref{sec6}$
by setting
\begin{equation}
\label{7}
g:=\pi_2^{*}t \circ \pi_1^{*}t^{-1}\text{.}
\end{equation}
Then, we have
\begin{equation*}
\mathrm{Rec}_{\pi}(\ex{\pi}(F))(x) = \mathrm{triv}_{i}(s(x))
\end{equation*}
for any point $x\in M$. A morphism $\overline{\gamma}: x \to y$ is mapped by the functor $s$ to some finite composition
\begin{equation*}
s(\bar\gamma) = \alpha_{n} \circ \gamma_n \circ  \alpha_{n-1} \circ ... \circ \gamma_2 \circ \alpha_{1} \circ \gamma_1 \circ \alpha_{0}
\end{equation*} 
of basic morphisms $\gamma_i: x_i \to y_i$ and $\alpha_i\in Y^{[2]}$, so that  we have
\begin{equation}
\label{11}
\mathrm{Rec}_{\pi}(\ex{\pi}(F))(\overline{\gamma}) = g(\alpha_{n}) \circ \mathrm{triv}_{i}(\overline{\gamma_n}) \circ g(\alpha_{n-1}) \circ ... \circ \mathrm{triv}_i(\overline{\gamma_1}) \circ g(\alpha_0)\text{.}
\end{equation}
Now we  are ready define the component of the natural equivalence $\zeta$ at a functor $F$. This component is a morphism in $\loctrivfunct{i}{\mathrm{Gr}}{1}$ and thus 
itself a natural equivalence
\begin{equation*}
\zeta(F):  F \to \mathrm{Rec}_{\pi}(\ex{\pi}(F))\text{.}
\end{equation*}
We define the component
of $\zeta(F)$ at a point $x\in M$ by
\begin{equation}
\label{8}
\zeta(F)(x) := t(s(x)): F(x) \to \mathrm{triv}_{i}(s(x))\text{.}
\end{equation}
Now we check that this is natural in $x$: let $\overline{\gamma}:x \to y$ be a morphism like the above one. The diagram whose commutativity we have to show splits along the decomposition (\ref{11}) into diagrams of two types:
\begin{equation*}
\alxydim{}{F(\pi(x_{i})) \ar[r]^-{t(x_i)} \ar[d]_{\pi_{*}\gamma_i} & \mathrm{triv}_i(x_i)  \ar[d]^{\mathrm{triv}_i(\gamma_i)} \\ F(\pi(y_i)) \ar[r]_-{t(y_i)} & \mathrm{triv}_i(y_i)}
\quad\text{ and }\quad
\alxydim{}{F(\pi(\pi_1(\alpha))) \ar[r]^-{t(\pi_1(\alpha))} \ar[d]_{\id} & \mathrm{triv}_i(\pi_1(\alpha))  \ar[d]^{g(\gamma_i)} \\ F(\pi(\pi_2(\alpha))) \ar[r]_-{t(\pi_2(\alpha))} & \mathrm{triv}_i(\pi_2(\alpha))\text{.}}
\end{equation*}
Both diagrams are indeed commutative, the one on the left because $t$ is natural in $y \in Y$ and the one on the right because of (\ref{7}).

 It remains to show that $\zeta$ is natural in $F$, i.e. we
have to prove the commutativity of the naturality diagram
\begin{equation}
\label{9}
\alxydim{@C=2cm@R=1.3cm}{F \ar[r]^-{\zeta(F)} \ar[d]_{\alpha}
& \mathrm{Rec}_{\pi}(\ex{\pi}(F)) \ar[d]^{\mathrm{Rec}_{\pi}(\ex{\pi}(\alpha))} \\ F' \ar[r]_-{\zeta(F')}
& \mathrm{Rec}_{\pi}(\ex{\pi}(F'))}
\end{equation}
for any natural transformation $\alpha:F \to F'$. Recall that $\ex{\pi}(\alpha)$
is the natural transformation
\begin{equation*}
h := t \circ \pi^{*}\alpha \circ t^{-1}: \mathrm{triv}_{i} \to \mathrm{triv}'_{i}
\end{equation*}
and that $\mathrm{Rec}_{\pi}(\ex{\pi}(\alpha))$ is the natural transformation
whose component at a point $x\in M$ is the morphism
\begin{equation*}
h(s(x)):\mathrm{triv}_{i}(s(x)) \to \mathrm{triv}'_{i}(s(x))
\end{equation*}
in $T$. Then, with definition (\ref{8}), the commutativity of the naturality
square (\ref{9}) becomes obvious.

\subsection{Proof of Theorem \ref{th4}}

\label{sec12}

We show that a Wilson line $\mathcal{W}^{\mathrm{tra},i}_{x_1,x_2}$
of a transport functor $\mathrm{tra}$ with $\mathrm{Gr}$-structure  is smooth. 
Let $c:U \to PM$ be a map such that $\Gamma(u,t):=c(u)(t)$ is smooth, let
$\pi:Y \to M$ be a surjective submersion, and let $(\mathrm{triv},t)$ be
a $\pi$-local $i$-trivialization of the transport functor $\mathrm{tra}$,
for which $\ex{\pi}(\mathrm{triv},t)$ is smooth. 
Consider the pullback diagram
\begin{equation*}
\alxydim{@=1.3cm}{\Gamma^{-1}Y \ar[r]^-{a} \ar[d]_{p} & Y \ar[d]^{\pi} \\
U \times [0,1] \ar[r]_-{\Gamma}& M}
\end{equation*}
with the surjective submersion $p:\Gamma^{-1}Y \to U \times [0,1]$. We have
to show that
\begin{equation}
\label{49}
\mathcal{W}_{x_1,x_2}^{\mathrm{tra},i} \circ \mathrm{pr} \circ c:U \to G
\end{equation}
is a smooth map. This can be checked locally in a neighbourhood of  a point $u\in U$. Let $t_j\in
I$ for $j=0,...,n$ be numbers with $t_{j-1} < t_{j}$ for $j=1,...,n$, and
$V_j$ open neighbourhoods of $u$ chosen
small enough to admit smooth local sections
\begin{equation*}
s_j: V_j \times [t_{j-1},t_{j}] \to \Gamma^{-1}Y\text{.}
\end{equation*}
Then, we restrict all these sections to the  intersection $V$ of all
the $V_j$. Let $\beta_j: t_{j-1} \to t_{j}$ be paths through $I$ defining smooth
maps
\begin{equation}
\label{46}
\tilde \Gamma_j: V \times I \to Y : (v,t) \mapsto a(s_j(V,\beta_j(t)))\text{,}
\end{equation}
 which can be considered as maps $\tilde c_j: V \to PY$. Additionally,
we define the smooth maps 
\begin{equation*}
\tilde\alpha_j:V \to Y^{[2]}: v \mapsto (\tilde\Gamma_{j-1}(v,1),\tilde\Gamma_{j}(v,0))\text{.}
\end{equation*}
Note that for any $v\in V$, both $\mathrm{pr}(\tilde c_j(v))$ and $\tilde\alpha_j(v)$
are morphisms in the universal path pushout $\upp{M}{\pi}$, namely
\begin{equation*}
\mathrm{pr}(\tilde c_j(v)): \tilde \Gamma_j(v,0) \to \tilde\Gamma_j(v,1)
\quad\text{ and }\quad
\tilde\alpha_j(v): \tilde\Gamma_{j-1}(v,1) \to \tilde\Gamma_{j}(v,0)\text{.} \end{equation*}
Taking their composition, we obtain a map
\begin{equation*}
\phi: V \to \mathrm{Mor}(\upp{M}{\pi}): v \to \tilde c_{n}(v) \circ \tilde\alpha_j(v)
\circ ... \circ \tilde\alpha_1(v) \circ \tilde c_0(v)\text{.}
\end{equation*}
Now we claim two assertions for  the composite
\begin{equation}
\label{45}
i^{-1} \circ (R_{(\mathrm{triv},g)})_1 \circ \phi: V \to \mathrm{Mor}(\mathrm{Gr})
\end{equation}
of $\phi$ with
the functor $R_{(\mathrm{triv},g)}: \upp{M}{\pi} \to \mathrm{Mor}(T)$ we have defined
in Section \ref{sec5}: first, it is smooth, and second, it coincides with
the restriction of $\mathcal{W}_{x_1,x_2}^{\mathrm{tra},i} \circ \mathrm{pr} \circ c$ to $V$, both assertions together
prove the smoothness of (\ref{49}). To show the first assertion, note that (\ref{45}) is the following assignment:
\begin{equation}
\label{50}
v \mapsto \mathrm{triv}(\tilde c_j(v)) \cdot \tilde g(\tilde \alpha_j(v)) \cdot
... \cdot \tilde g(\tilde\alpha_1(v)) \circ \mathrm{triv}(\tilde c_0(v))\text{.}
\end{equation}
By definition, the descent data $(\mathrm{triv},g)$ is smooth. Because
$\mathrm{triv}$ is a smooth functor, and the maps $\tilde c_j$ satisfy the
relevant condition (\ref{46}), every factor $\mathrm{triv} \circ \tilde c_j :V \to G$ is smooth. Furthermore, the maps  $\tilde g:Y^{[2]} \to G$ are
smooth, so that
also the remaining factors are smooth in $v$. To show the second assertion,
consider a point $v\in V$. If we choose in the definition of the Wilson line
$\mathcal{W}_{x_1,x_2}^{\mathrm{tra},i}$ the objects $G_k := \mathrm{triv}(\tilde
x_k)$ and the isomorphisms $t_k := t(\tilde x_k)$ for some lifts $\pi(\tilde
x_k)=x_k$, where $t$ is the trivialization of $\mathrm{tra}$ from the beginning
of this section, we find
\begin{equation*}
(\mathcal{W}_{x_1,x_2}^{\mathrm{tra},i} \circ \mathrm{pr} \circ c)(v)=\mathrm{tra}(\overline{c(v)})\text{.}
\end{equation*}
The right hand side coincides with the right hand side of (\ref{50}). 

\subsection{Proof of Proposition \ref{lem6}}

\label{appD}

We are going  to prove that the map $k_A:PX \to G$, defined by 
\begin{equation*}
k_A(\gamma):=f_{\gamma^{*}A}(0,1)
\end{equation*}
for a path $\gamma:[0,1] \to X$ depends only on the thin homotopy class of
$\gamma$. Due to the multiplicative property
(\ref{22})
of $k_A$, it is enough to show $k_A(\gamma_0^{-1} \circ\gamma_1)=1$ for every thin homotopy equivalent paths $\gamma_0$ and $\gamma_1$. For this purpose we derive a relation to the pullback of the curvature $K:=\mathrm{d}A + [A \wedge A]$ of the 1-form $A$ along a homotopy between $\gamma_0$ and $\gamma_1$. If this homotopy is thin, the pullback vanishes. 

Let us fix the following notation. $Q:=[0,1] \times [0,1]$ is the unit square, $\gamma_{(a,b,c,d)} :(a,b) \to (c,d)$ is the straight
path in $Q$, and
\begin{figure}
\begin{center}
\includegraphics{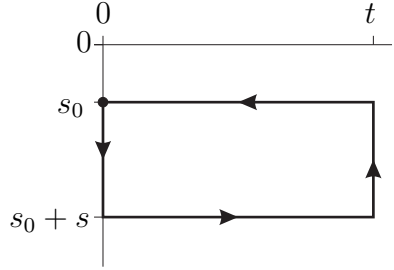}\setlength{\unitlength}{1pt}\begin{picture}(0,0)(204,496)\put(72.59414,554.30079){$s_0$}\put(54.59414,511.10079){$s_0+s$}\put(79.79414,577.04882){$0$}\put(86.99414,587.84882){$0$}\put(187.79414,587.84882){$t$}\end{picture}
\begin{minipage}[c]{0.8\textwidth}
\caption{The path $\tau_{s_0}(s,t)$.\label{fig2}
}
\end{minipage}
\end{center}
\end{figure}
\begin{equation*}
\tau_{s_0}: Q \to PQ
\end{equation*}
assigns for fixed $s_0\in [0,1]$ to a point $(s,t) \in Q$ the closed path
\begin{equation*}
\tau_{s_0}{(s,t)}:= \gamma_{(s_0,t,s_0,0)} \circ \gamma_{(s_0+s,t,s_0,t)} \circ \gamma_{(s_0+s,0,s_0+s,t)} \circ \gamma_{(s_0,0,s_0+s,0)}\text{,}
\end{equation*}
which goes counter-clockwise around the rectangle spanned by the points $(s_0,0)$ and $(s_0+s,t)$, see Figure \ref{fig2}. Now consider two paths $\gamma_0,\gamma_1:x \to y$ in $X$. Without loss of generality we can assume that the paths $\gamma_{(a,b,c,d)}$ used above have sitting instants, such that $\tau_{s_0}$ is smooth and  
\begin{equation}
\label{560}
\gamma_0(\gamma_{(0,1,0,0)}(t)) = \gamma_{0}^{-1}(t) 
\quad\text{ and }\quad
\gamma_1(\gamma_{(0,1,1,1)}(t)) = \gamma_1(t)\text{.}
\end{equation}

\begin{lemma}
\label{lem13}
Let $h:Q \to X$ be a smooth homotopy between the paths $\gamma_0,\gamma_1:x \to y$ with  $h(0,t)=\gamma_0(t)$
and $h(1,t)=\gamma_1(t)$. Then, the map
\begin{equation*}
u_{A,s_0}:= k_A \circ h_{*} \circ \tau_{s_0}: Q \to G
\end{equation*}
is smooth and has the following properties
\begin{enumerate}
\item[(a)]
$u_{A,0}(1,1)=k_A(\gamma_{0}^{-1}
\circ \gamma_{1})$

\item[(b)]
$u_{A,s_0}(s,1)= u_{A,s_0}(s',1) \cdot u_{A,s_0+s'}(s-s',1)$

\item[(c)]
with $\gamma_{s,t}$ the path defined by  $\gamma_{s,t}(\tau) := h(s,\tau t)$ and
 $K := \mathrm{d}A +[A \wedge A]$ the curvature of $A$ we have:
\begin{equation}
\label{52}
\left .    \frac{\partial}{\partial s}
    \frac{\partial}{\partial t}
    u_{A,s_0} \right |_{(0,t)}
    = -
    \mathrm{Ad}_{k_A(\gamma_{s_0,t})}^{-1} 
     \left ( h^*K \right )_{(s_{0},t)} \left (\frac{\partial}{\partial s},\frac{\partial}{\partial t} \right)
\end{equation}
\end{enumerate}
\end{lemma}

\begin{proof}
Since $h$ is constant for $t=0$ and $t=1$, (a) follows from \erf{560}. The multiplicative property
(\ref{22})
of $k_A$ implies (b). To prove (c), we define a further path $\gamma_{s_0,s,t}(\tau):= h(s_0+s\tau,t)$ and write
\begin{equation}
\label{20}
u_{A,s_0}(s,t) = f_{\gamma_{s_0,t}^{*}A}(0,1)^{-1} \cdot f_{\gamma^{*}_{s_0,s,t}A}(0,1)^{-1} \cdot f_{\gamma_{s_0+s,t}^{*}A}(0,1)
\end{equation}
where $f_{\varphi}: \R \times \R \to G$ are the smooth functions that correspond to the a 1-form $\varphi\in\Omega^1(\R,\mathfrak{g})$ by Lemma \ref{lem2} as the solution of initial value problems. By a uniqueness argument one can show that
$f_{\gamma_{s,t}^{*}A}(0,1) = f_{\gamma_{s,1}^{*}A}(0,t)$. 
Then, we calculate with (\ref{20}) and, for simplicity, in a faithful matrix representation of $G$,
\begin{multline*}
\frac{\partial}{\partial t}u_{A,s_0}(s,t) =f^{-1}_{\gamma_{s_0,t}^{*}A}(0,1) \cdot \left ( (h^{*}A)_{(s_0,t)}\left ( \frac{\partial }{\partial t}  \right ) \right . \cdot f_{\gamma_{s_0,s,t}^{*}A}(0,1)^{-1}   
\\ +\frac{\partial}{\partial t} f_{\gamma_{s_0,s,t}^{*}A}(0,1)^{-1}   \left . - f_{\gamma^{*}_{s_0,s,t}A}(0,1)^{-1} \cdot   (h^{*}A)_{(s_0+s,t)}\left (  \frac{\partial }{\partial t}  \right ) \right ) \cdot f_{\gamma_{s_0+s,t}^{*}A}(0,1)\text{.}
\end{multline*}
To take the derivatives along $s$, we use  $f_{\gamma_{s_0,s,t}^{*}A}(0,1)= f_{\gamma_{s_0,1,t}}(0,s)$ and $f_{\gamma_{s_0,0,t}}(0,1)=1$, both together show 
\begin{equation*}
\left . \frac{\partial}{\partial s}\right|_{0} f_{\gamma_{s_0,s,t}^{*}A}(0,1)^{-1}= (h^{*}A)_{(s_0,t)} \left (\frac{\partial}{\partial s} \right )\text{.}
\end{equation*}
Finally,
\begin{eqnarray*}
\left . \frac{\partial}{\partial s}\frac{\partial }{\partial t} u_{A,s_0} \right |_{s=0}\!\!\! &=& f^{-1}_{\gamma_{s_0,t}^{*}A}(0,1) \cdot \left ( \left ( (h^{*}A)_{(s_0,t)}\left ( \frac{\partial }{\partial t}  \right ) \right . \cdot(h^{*}A)_{(s_0,t)} \left (\frac{\partial}{\partial s} \right ) \right. \\&& \!\! + \frac{\partial}{\partial t}(h^{*}A)_{(s_0,t)} \left (\frac{\partial}{\partial s} \right ) - (h^{*}A)_{(s_0,t)} \left (\frac{\partial}{\partial s} \right ) \cdot  (h^{*}A)_{(s_0,t)} \left (\frac{\partial}{\partial t} \right ) \\&&\!\! \left .- \left .\frac{\partial}{\partial s}\right|_0   (h^{*}A)_{(s_0+s,t)}\left (  \frac{\partial }{\partial t}  \right ) \right ) \cdot f_{\gamma_{s_0,t}^{*}A}(0,1)\text{.} \end{eqnarray*}
This yields the claimed equality.
\end{proof}

Notice that if $h$ is a thin homotopy, $h^{*}K=0$, so that the right hand side
in (c) vanishes.  
Then we calculate at $(0,1)$
\begin{equation*}
\left .\frac{\partial}{\partial s} u_{A,s_0} \right |_{(0,1)} = \int_0^1  \left .\frac{\partial}{\partial s}\frac{\partial}{\partial t}
    u_{A,s_0} \right|_{(0,t)} \mathrm{d}t = 0\text{.}
\end{equation*}
Using (b) we obtain the same result
for all points $(s_0,1)$,
\begin{equation*}
\left .\frac{\partial}{\partial s} u_{A,0} \right |_{(s_0,1)}= u_{A,0}(s_0,1) \cdot  \left . \frac{\partial}{\partial s}  u_{A,s_0}  \right|_{(0,1)}= u_{A,0}(s_0,1)
\cdot 0 =0.  
\end{equation*}
This means that the function $u_{A,0}(s,1)$ is constant and thus determined
by its value at $s=0$, namely
\begin{equation*}
1=u_{A,0}(0,1)=u_{A,0}(1,1) \stackrel{\mathrm{(a)}}{=}k_A(\gamma_1^{-1} \circ \gamma_0)=k_A(\gamma_1)^{-1} \cdot k_A(\gamma_0)\text{.}
\end{equation*} 
This finishes the proof.

\subsection{Proof of Proposition \ref{th2}}

\label{appE}

We have to show that the functor $\fu: \diffco{G}{1}{X} \to \funct^{\infty}(\mathcal{P}_1(X),\mathcal{B}G)$ is bijective on objects. For this purpose, we define an inverse map $\fo$ that assigns to any smooth functor $F:\mathcal{P}_1(X) \to \mathcal{B}G$  a $\mathfrak{g}$-valued 1-form
$\fo(F)$, such
that $\fu(\fo(F))=F$ and such that $\fo(\fu(A))=A$ for any 1-form $A\in\Omega^1(X,\mathfrak{g})$. 

Let $F:\mathcal{P}_1(X) \to \mathcal{B}G$ be a smooth functor.
We define the 1-form $A:=\fo(F)$ at a point $p\in X$ and for a tangent vector $v\in T_pX$ in the following way. Let $\gamma: \R \to X$ be a smooth curve such that $\gamma(0)=p$ and $\dot \gamma(0)=v$.
We consider the map 
\begin{equation}
\label{43}
f_{\gamma}:= F_1 \circ (\gamma_{*})_1: \R \times \R \to G\text{.}
\end{equation}
The evaluation $\mathrm{ev}\circ ((\gamma_{*})_1 \times \id): U \times [0,1]
\to X$ with $U=\R \times \R$ is a smooth map
because $\gamma$ is smooth. Hence,
by Definition \ref{def1}, $f_{\gamma}$ is smooth.
The properties of the functor $F$ further imply the cocycle condition
\begin{equation}
\label{35}
f_{\gamma}(y,z) \cdot f_{\gamma}(x,y)=f_{\gamma}(x,z)\text{.}
\end{equation} 
By Lemma \ref{lem2}, the smooth map $f_{\gamma}: \R \times \R \to G$ corresponds to a $\mathfrak{g}$-valued 1-form $A_{\gamma}$ on $\R$. We define
\begin{equation}
\label{defA}
\alpha_{F,\gamma}(p,v) := A_{\gamma}|_0 \left (\frac{\partial}{\partial t} \right) \in \mathfrak{g}\text{.}
\end{equation}
With a view to the definition (\ref{27}) of $A_{\gamma}$, this is
\begin{equation}
\label{42}
\alpha_{F,\gamma}(p,v)=- \left . \frac{\mathrm{d}}{\mathrm{d}t}f_{\gamma}(0,t)
\right|_{t=0}\text{.}
\end{equation}

\begin{lemma}
\label{lem11}
$\alpha_{F,\gamma}(p,v)$ is independent of the choice of the smooth curve
$\gamma$.
\end{lemma}

\begin{proof}
We prove the following: if $\gamma_0,\gamma_1:\R \to X$ are smooth curves such that $\gamma_0(0) = \gamma_1(0)$ and $\left . \frac{\mathrm{d}\gamma_1(t)}{\mathrm{d}t} \right |_0=\left . \frac{\mathrm{d}\gamma_2(t)}{\mathrm{d}t} \right |_0$ then it follows
\begin{equation*}
 \left . \frac{\mathrm{d}}{\mathrm{d}t} \right |_{0} f_{\gamma_0}(0,t)= \left . \frac{\mathrm{d}}{\mathrm{d}t} \right |_{0} f_{\gamma_1}(0,t)\text{.}
\end{equation*}
Let $U \subset X$ be a convex open chart neighborhood of $x$. Let $\epsilon>0$ be such that  $(\gamma_{0\ast})_1(0,t)(x),(\gamma_{1\ast})_1(0,t)(x)\in U$ for all $t\in [0,\varepsilon)$ and $x\in[0,1]$. In the chart, we can form the difference vector
\begin{equation*}
d:[0,\varepsilon) \times [0,1] \to \R^n: (t,x) \mapsto  (\gamma_{1\ast})_1(0,t)(x)-(\gamma_{0\ast})_1(0,t)(x)\text{.}
\end{equation*}
We obtain a smooth map
\begin{equation*}
h: [0,\varepsilon) \times [0,1] \times [0,1] \to X: (t,\alpha,x) \mapsto (\gamma_{0\ast})_1(0,t)(x) + \alpha d(t,x)
\end{equation*}
and define by $H(t,\alpha)(x) := h(t,\alpha,x)$ a smooth homotopy
\begin{equation*}
H: [0,\varepsilon) \times [0,1] \to PX
\end{equation*}
between $(\gamma_{0\ast})_1(0,t)$ and $(\gamma_{1\ast})_1(0,t)$.

The difference map satisfies $d(0,x)=0$ and $\left  . \frac{\mathrm{d}}{\mathrm{d}t} \right |_0 d(t,x) =0$. By Hadamard's Lemma there exists a smooth map $e: [0,\varepsilon) \times [0,1] \to \R^n$ such that $d(t,x)=t^2 e(t,x)$.  We consider $Z := [0,\varepsilon) \times [0,\varepsilon^2)$ and the smooth maps
\begin{equation*}
f: [0,\varepsilon) \times [0,1] \to Z
\quad\text{ and }\quad
p:Z \to PX
\end{equation*} 
defined by $f(t,\alpha) := (t,t^2\alpha)$ and $p(t,\alpha)(x) := \gamma_0(t)(x) + \alpha e(t,x)$. We have  $H=p \circ f$. 
Now we compute via the chain rule
\begin{equation*}
\left . \frac{\mathrm{d}}{\mathrm{d}t} \right |_{0} F_1(H(t,\alpha)) = J(F_1 \circ p)|_{f(0,\alpha)}\cdot \left . \frac{\mathrm{d}}{\mathrm{d}t} \right |_{0}f(t,\alpha)= J(F_1 \circ p)|_{(0,0)}\cdot (1,0)\text{.}
\end{equation*}
The left hand side is $\left . \frac{\mathrm{d}}{\mathrm{d}t} \right |_{0} f_{\gamma_\alpha}(0,t)$, and the right hand side is independent of $\alpha$. This shows the claim.  
\end{proof} 

According to the result of Lemma \ref{lem11}, we drop the index $\gamma$,
and remain with a map $\alpha_F:TX \to \mathfrak{g}$ defined canonically
by the functor $F$. 
We show next that $\alpha_F$ is linear in $v$. For a multiple
$s v$ of $v$ we can choose the curve $\gamma_{s}$ with $\gamma_{s}(t):=
\gamma(s t)$. It is easy to see that then $f_{\gamma_{s}}(x,y)=f_{\gamma}(s
x,s y)$. Again by the chain rule 
\begin{equation*}
\alpha(p,sv)=-\frac{\mathrm{d}}{\mathrm{d}t}f_{\gamma_{s}}(0,t)|_{t=0}
=-  \frac{\mathrm{d}}{\mathrm{d}t}f_{\gamma}(0,s t)|_{t=0} = s \alpha_F(p,v)\text{.}
\end{equation*} 
In the same way one can show 
that $\alpha(p,v+w) = \alpha(p,v)+
\alpha(p,w)$.  

\begin{lemma}
The pointwise linear map $\alpha_F:TX \to \mathfrak{g}$ is smooth, and thus
defines a 1-form $A\in \Omega^1(X,\mathfrak{g})$ by $A|_p(v):= \alpha_F(p,v)$.

\end{lemma}

\begin{proof}
If $X$ is $n$-dimensional and $\phi : U \to X$ is a coordinate chart with
an open subset $U \subset \R^n$, the 
standard chart for the
tangent bundle $TX$ is
\begin{equation*}
\phi_{TX}:U\times \mathbb{R}^n \to TX:(u,v) \mapsto \mathrm{d}\phi|_u(v)\text{.} \end{equation*}
We prove the smoothness of $\alpha_F$ in the chart $\phi_{TX}$, i.e. we show that 
\begin{equation*}
A \circ \phi_{TX} : U \times \mathbb{R}^n \to \mathfrak{g}
\end{equation*}
is smooth. For this purpose, we define the map
\begin{equation*}
c : U \times \mathbb{R}^n  \times \R \to PX: (u,v,\tau) \mapsto (t \mapsto \phi(u+\beta (t \tau) v))
\end{equation*}
where $\beta$ is some smoothing function, i.e. an orientation-preserving diffeomorphism of $[0,1]$ with sitting instants. Now, $\mathrm{ev} \circ
(c \times \id)$ is evidently smooth in all parameters, and since  $F$ is a smooth functor,
\begin{equation*}
f_{c} := F_1 \circ \mathrm{pr} \circ c : U \times \mathbb{R}^n \times \mathbb{R}
  \to G
\end{equation*}
is a smooth function. Note that $\gamma_{u,v}(t) := c(u,v,t)(1)$ defines a smooth curve in $X$
with the properties
\begin{equation}
\label{54}
\gamma_{u,v}(0)=\phi(u)
\quad\text{ and }\quad
\dot\gamma_{u,v}=\mathrm{d}\phi|_u(v)\text{,}
\end{equation}
and which is in turn related to $c$ by
\begin{equation}
\label{55}
(\gamma_{u,v})_{*}(0,t)=c(u,v,t)\text{.}
\end{equation}
Using the path $\gamma_{u,v}$ in the definition of the 1-form $A$, we find
\begin{eqnarray*}
 (A \circ \phi_{TX})(u,v) &=&\alpha_F(\phi(u),\mathrm{d}\phi|_u(v))\nonumber
  \\ &\stackrel{\mathrm{(\ref{54})}}{=}& -\left .\frac{\mathrm{d}}{\mathrm{d}t}(F_1 \circ (\gamma_{u,v})_{*})(0,t)\right
  |_{t=0} \nonumber \\ &\stackrel{\mathrm{(\ref{55})}}{=}& -\left .
  \frac{\mathrm{d}}{\mathrm{d}t} f_c(u,v,t) \right |_{t=0}\text{.}
\end{eqnarray*}
The last expression  is, in particular, smooth in $u$ and $v$.
\end{proof}

Summarizing, we started with a given smooth functor $F:\mathcal{P}_1(X) \to
\mathcal{B}G$ and have derived a 1-form $\fo(F):=A\in\Omega^1(X,\mathfrak{g})$.
Next we show that this 1-form is the preimage of $F$ under the functor
\begin{equation*}
\fu: \gconn{X} \to \funct^{\infty}(X,\mathcal{B}G)
\end{equation*}
from Proposition \ref{th2}, i.e. we show 
\begin{equation*}
\fu(A)(\overline{\gamma})=F(\overline{\gamma})
\end{equation*}
for any path $\gamma\in PX$. We recall that the functor $\fu(A)$ was defined
by $\fu(A)(\overline{\gamma}):=f_{\gamma^{*}A}(0,1)$, where $f_{\gamma^{*}A}:\R
\times \R \to G$ solves the differential equation
\begin{equation}
\label{37}
\frac{\mathrm{d}}{\mathrm{d}t}f_{\gamma^{*}A}|_{(0,t)}=-\mathrm{d}r_{f_{\gamma^{*}A}(0,t)}|_1 \left ((\gamma^{*}A)_t\left
( \frac{\partial}{\partial t} \right )\right)
\end{equation}    
with the initial value $f_{\gamma^{*}A}(0,0)=1$. 
Now we use the construction of the 1-form $A$ from the given functor $F$. For the smooth function $f_{\gamma}:\R \times \R \to G$ from (\ref{43}) we obtain using  $\gamma_t(\tau) := \gamma(t+\tau)$ with $p:=\gamma_t(0)$ and $v:=\dot\gamma_t(0)$ 
\begin{multline}
\label{39}
\frac{\mathrm{d}}{\mathrm{d}\tau} f_{\gamma}(t,t+\tau)|_{\tau=0} =\frac{\mathrm{d}}{\mathrm{d}\tau}f_{\gamma_t}(0,\tau)|_{\tau=0}\\\stackrel{\mathrm{(\ref{42})}}{=}-\alpha_{F,\gamma_t}(p,v)=-A_p(v) =- (\gamma^{*}A)_{t} \left (\frac{\partial}{\partial t} \right)\text{.}
\end{multline}
Then we have
\begin{equation*}
\frac{\mathrm{d}}{\mathrm{d}t} f_\gamma(0,t) \stackrel{\mathrm{(\ref{35})}}{=} \frac{\mathrm{d}}{\mathrm{d}\tau}f_{\gamma}(t,\tau)|_{\tau=t} \cdot f_{\gamma}(0,t) \stackrel{\mathrm{(\ref{39})}}{=}- \mathrm{d}r_{f_{\gamma}(0,t)}|_1 \left((\gamma^{*}A)_{t} \left (\frac{\partial}{\partial t} \right)  \right )\text{.}
\end{equation*} 
Hence, $f_{\gamma}$ solves the initial value problem (\ref{37}). By uniqueness,  $f_{\gamma^{*}A}=f_{\gamma}$
and finally
\begin{equation*}
F(\overline{\gamma})= f_{\gamma}(0,1)=f_{\gamma^{*}A}(0,1)=\fu(A)(\overline{\gamma})\text{.}
\end{equation*} 

\medskip

It remains to show that, conversely, for a given 1-form $A\in\Omega^1(X,\mathfrak{g})$, \begin{equation*}
\fo(\fu(A))=A\text{.}
\end{equation*}
We test the 1-form $\fo(\fu(A))$ at a point $x\in X$ and a tangent vector $v\in T_xX$. Let $\Gamma:\R \to X$ be a curve in $X$ with $x=\Gamma(0)$ and $v=\dot \Gamma(0)$. If we further denote $\gamma_{a} := \Gamma_{*}(0,a)$ we have
\begin{equation}
\label{eq:neu}
-\fo(\fu(A))|_x(v) \stackrel{\mathrm{(\ref{42})}}{=} \left . \frac{\partial f_{\gamma_{a}}}{\partial a} \right|_{(0,0)}\stackrel{\mathrm{(\ref{43})}}{=} \left . \frac{\partial}{\partial a} \right |_{0}\fu(A)(\gamma_{a})= \left . \frac{\partial}{\partial a} \right |_{0}f_{\gamma_{a}^{*}A}(0,1)
\end{equation}
Here, $f_{\gamma_{a}^{*}A}$ is the unique solution of the initial value problem (\ref{37}) for the given 1-form $A$ and the curve $\gamma_{a}$. We calculate
with the product rule
\begin{multline*}
\left . \frac{\mathrm{d}}{\mathrm{d}a} \right |_0 \left . \frac{\mathrm{d}f_{\gamma_{a}^{*}A}(0,\tau)}{\mathrm{d}\tau} \right |_{\tau=t} =- \left . \frac{\mathrm{d}}{\mathrm{d}a} \right |_0 A_{\gamma_a(t)}\left ( \left . \frac{\mathrm{d}\gamma_a(\tau)}{\mathrm{d}\tau} \right|_{\tau=t} \right ) \cdot f_{\gamma_{0}^{*}A}(0,t)\\ -A_{\gamma_0(t)}\left ( \left . \frac{\mathrm{d}\gamma_0(\tau)}{\mathrm{d}\tau} \right|_{\tau=t} \right ) \cdot \left . \frac{\mathrm{d}f_{\gamma_{a}^{*}A}(0,t) }{\mathrm{d}a} \right |_0\text{.} 
\end{multline*}
Since  $\gamma_0=\id_x$, we have $\left . \frac{\mathrm{d}\gamma_0(\tau)}{\mathrm{d}\tau} \right|_{\tau=t}=0$ as well as $f_{\gamma_{0}^{*}A}(0,t)=1$. Hence,
\begin{eqnarray*}
\left . \frac{\mathrm{d}}{\mathrm{d}a} \right |_0 \left . \frac{\mathrm{d}f_{\gamma_{a}^{*}A}(0,\tau)}{\mathrm{d}\tau} \right |_{\tau=t} &=& - \left . \frac{\mathrm{d}}{\mathrm{d}a} \right |_0 A_{\gamma_a(t)}\left ( \left . \frac{\mathrm{d}\gamma_a(\tau)}{\mathrm{d}\tau} \right|_{\tau=t} \right ) 
\\&=& - \lim_{h\rightarrow 0} \frac{1}{h}  A_{\gamma_h(t)}\left ( \left . \frac{\mathrm{d}\gamma_h(\tau)}{\mathrm{d}\tau} \right|_{\tau=t} \right )\text{.} 
\end{eqnarray*}
Within the thin homotopy class of $\gamma_a$ we find a representative that satisfies $\gamma_a(t)=\Gamma(a\varphi(t))$ with $\varphi:[0,1] \to [0,1]$ a surjective smooth map with sitting instants. Using this representative we obtain
\begin{equation*}
\left . \frac{\mathrm{d}\gamma_h(\tau)}{\mathrm{d}\tau} \right|_{\tau=t} = \left . \frac{\mathrm{d}\Gamma}{\mathrm{d}\tau}\right |_{\tau=h\varphi(t)} \cdot h \cdot \left . \frac{\mathrm{d}\varphi(\tau)}{\mathrm{d}\tau} \right |_{\tau=t}\text{.}
\end{equation*}
Thus,
\begin{eqnarray*}
\left . \frac{\mathrm{d}}{\mathrm{d}a} \right |_0 \left . \frac{\mathrm{d}f_{\gamma_{a}^{*}A}(0,\tau)}{\mathrm{d}\tau} \right |_{\tau=t} &=& - \lim_{h\rightarrow 0} \left . \frac{\mathrm{d}\varphi(\tau)}{\mathrm{d}\tau} \right |_{\tau=t}\cdot A_{\gamma_h(t)}\left (  \left . \frac{\mathrm{d}\Gamma}{\mathrm{d}\tau}\right |_{\tau=h\varphi(t)} \right )\\&=&- \left . \frac{\mathrm{d}\varphi(\tau)}{\mathrm{d}\tau} \right |_{\tau=t} A_x(v)\text{.}
\end{eqnarray*}
Integrating $t$ form 0 to 1 gives with $\varphi(1)=1$
\begin{equation*}
\left . \frac{\mathrm{d}f_{\gamma_a}(1)}{\mathrm{d}a} \right |_0 = - A_x(v)\text{.} \end{equation*}
Substituting this in \erf{eq:neu} we obtain the claimed result.

\end{appendix}

\newpage

\section*{Table of Notations}
\addcontentsline{toc}{section}{Table of Notations}

\label{sec13}

\newcommand{\notation}[3]{
\noindent
\begin{minipage}[t]{0.15\textwidth}#1\end{minipage}
\begin{minipage}[t]{0.68\textwidth}#2\vspace{0.3cm}\end{minipage}\hfill
\begin{minipage}[t]{0.105\textwidth}Page \pageref{#3}\end{minipage}
}

\notation{$\upp{M}{\pi}$}{the universal path pushout of a surjective
submersion $\pi:Y \to M$.}{sec5}

\notation{$\ex{\pi}$}{the functor $\ex{\pi}: \loctrivfunct{i}{\mathrm{Gr}}{1} \to \trans{i}{i}{1}$ which
extracts descent data.}{4}

\notation{$f_{*}$}{the functor $f_{*}:\mathcal{P}_1(M) \to \mathcal{P}_1(N)$
of path groupoids induced by a smooth map $f:M \to N$.}{def8}

\notation{$\funct(S,T)$}{the category of functors between categories $S$and $T$.}{p4}

\notation{$\funct^{\infty}(\mathcal{P}_1(X),\mathrm{Gr})$}{\quad\\ the category of smooth functors from $\mathcal{P}_1(X)$ to a Lie groupoid $\mathrm{Gr}$.
}{p3}

\notation{$G\text{-}\mathrm{Tor}$}{the category of smooth principal $G$-spaces
and smooth equivariant maps.}{sec8}

\notation{$\mathrm{Gr}$}{a Lie groupoid}{sec7}

\notation{$i$}{a  functor $i: \mathrm{Gr} \to T$, which relates the typical
fibre $\mathrm{Gr}$ of a functor to the target category $T$.}{sec7}

\notation{$\sm$}{the category of smooth manifolds and smooth maps between
those.}{appB}

\notation{$\pi_1^1(M,x)$}{the thin homotopy group of the smooth manifold
$M$ at the point $x\in M$.}{44}

\notation{$\uppp$}{the projection functor $\uppp:\upp{M}{\pi} \to \mathcal{P}_1(M)$.}{51}

\notation{$PM$}{the set of paths in $M$}{def9}

\notation{$P^1M$}{the set of thin homotopy classes of paths in $M$}{def10}

\notation{$\mathcal{P}_1(M)$}{the path groupoid of the smooth manifold $M$.}{def8}

\notation{$\Pi_1(M)$}{the fundamental groupoid of a smooth manifold $M$}{p1}

\notation{$\mathrm{Rec}_{\pi}$}{the functor $\mathrm{Rec}_{\pi}: \trans{i}{i}{1}
\to \loctrivfunct{i}{\mathrm{Gr}}{1}$ which reconstructs a functor from descent data.}{17}

\notation{$s$}{the section functor $s:\mathcal{P}_1(M) \to \upp{M}{\pi}$
associated to a surjective submersion $\pi:Y \to M$.}{sec5}

\notation{$\mathcal{B}G$}{the category with one object whose set of morphisms
is the Lie group $G$.}{def17}

\notation{$\smsp$}{the category of smooth spaces}{appB}

\notation{$T$}{the target category of transport functors -- the fibres of
a bundle are objects in $T$, and the parallel transport maps are morphisms
in $T$.}{sec7}

\notation{$\trans{i}{i}{1}$}{the category of descent data of $\pi$-locally $i$-trivialized functors.}{def7}

\notation{$\transsmooth{i}{1}{i}$}{the category of smooth descent data of $\pi$-locally $i$-trivialized functors}{def2}

\notation{$\transport{i}{1}{\mathrm{Gr}}{T}$}{\quad\\the category of transport functors with $\mathrm{Gr}$-
structure.}{def16}

\notation{$\loctrivfunct{i}{\mathrm{Gr}}{1}$}{the category of functors $F:\mathcal{P}_1(M) \to T$ together
with $\pi$-local $i$-trivializations}{def6}

\notation{$\loctrivfunctsmooth{i}{\mathrm{Gr}}{\pi}{1}$}{the category of transport functors on $M$ with $\mathrm{Gr}$-structure together
with $\pi$-local $i$-trivializations.}{def3}

\notation{$\mathrm{Vect}(\C_h^n)$}{the category of $n$-dimensional hermitian vector spaces and isometries between those.
}{57}

\notation{$\mathrm{VB}(\C^n_h)^{\nabla}_M$}{the category of hermitian vector
bundles of rank $n$ with unitary connection over $M$.}{p5}

\notation{$\gconn{X}$}{the category of differential cocycles on $X$ with
gauge group $G$.}{def11}

\notation{$\diffco{G}{1}{\pi}$}{the category of differential cocycles of
a surjective submersion $\pi$
with gauge group $G$.}{def12}

\newpage

\kobib{../../bibliothek/tex}

\end{document}